\documentclass[a4paper,11pt]{article}
\usepackage[english]{babel}

\usepackage[T1]{fontenc}

\usepackage[utf8]{inputenc}
\usepackage{amssymb}
\usepackage{amsmath}
\usepackage{amsthm}
\usepackage{mathrsfs}
\usepackage{mathtools}
\usepackage{amsfonts}
\usepackage{anyfontsize}
\usepackage[bb=ams]{mathalpha}
\usepackage[normalem]{ulem}

\usepackage[nottoc,notlot,notlof]{tocbibind}
\usepackage{tikz-cd}
\usepackage[all,cmtip]{xy}
\usepackage{graphicx,tikz-cd,pgf}
\usepackage{enumerate}
\usepackage[margin=2.5cm]{geometry}
\usepackage{leftindex}
\usepackage{url}
\usepackage{stix}
\usepackage{bigints}
\usepackage{cite}
\usepackage{stackengine}
\tikzcdset{arrow style=tikz}
\tikzset{mymatr/.style={every outer matrix/.append style={draw=black, inner xsep=15pt , inner ysep=6pt, rounded corners, very thick}}}

\newcommand{\ltr}[1]{{#1}_{\bot}}
\newcommand{\rtl}[1]{{#1}_{\top}}

\newcommand{\bspace}{\noindent\,\,$\bullet$\,\,\,}

\newcommand{\blspace}{\bigskip}

\newcommand{\twinf}[1]{\sTwin(#1)}
\newcommand{\squaref}[1]{\sSquareF#1}

\newcommand{\csquare}[1]{#1\times{#1}}

\newcommand{\Setg}{U}
\newcommand{\Setuno}{\Setg}
\newcommand{\Setdue}{V}
\newcommand{\Settre}{W}

\newcommand{\ssdomsymbol}{\mathtt{d}}
\newcommand{\sscodomsymbol}{\mathtt{c}}
\newcommand{\ssdom}[1]{{#1}_{\ssdomsymbol}}
\newcommand{\sscodom}[1]{{#1}_{\sscodomsymbol}}

\DeclareMathOperator{\dg}{\Lambda}

\newcommand{\domcod}[1]{\imath_{#1}}

\newcommand{\onedomcod}[1]{{\ssdomsymbol}_{#1}}
\newcommand{\twodomcod}[1]{\sscodomsymbol_{#1}}
\newcommand{\onedomcodof}[2]{
\onedomcod{#1}(#2)}
\newcommand{\twodomcodof}[2]{
\twodomcod{#1}(#2)
}

\newcommand{\sd}[2]{#1\Delta#2}

\newcommand{\R}{\mathbb R}

\makeatletter
\newcommand*{\bigcdot}{}
\DeclareRobustCommand*{\bigcdot}{%
  \mathbin{\mathpalette\bigcdot@{}}%
}
\newcommand*{\bigcdot@scalefactor}{0.75}
\newcommand*{\bigcdot@widthfactor}{1.20}
\newcommand*{\bigcdot@}[2]{%
  \sbox0{$#1\vcenter{}$}
  \sbox2{$#1\cdot\m@th$}%
  \hbox to \bigcdot@widthfactor\wd2{%
    \hfil
    \raise\ht0\hbox{%
      \scalebox{\bigcdot@scalefactor}{%
        \lower\ht0\hbox{$#1\bullet\m@th$}%
      }%
    }%
    \hfil
  }%
}
\makeatother

\newcommand{\vm}[2]{
{#1}{\,\bigcdot\,}{#2}}
\newcommand{\abstractvm}[2]{
{#1}{\,*\,}{#2}}
\newcommand{\cdm}[2]{
{#1}{\,\cdot\,}{#2}}
\newcommand{\abstractcdm}[2]{
{#1}{\,\odot\,}{#2}}

\newcommand{\hm}[2]{
{#1}{\,\hmult\,}{#2}}

\DeclareMathOperator{\sCatObjectSymbol}{\mathtt{obj}}

\DeclareMathOperator{\sCatArrowsSymbol}{\mathtt{arr}}
\DeclareMathOperator{\sCatArrowsSymbolpm}{\boldsymbol{\mathsf{Arr}}}

\newcommand{\sCatO}[1]{\sCatObjectSymbol[#1]}
\newcommand{\sCatA}[1]{\sCatArrowsSymbol[#1]}
\newcommand{\sCatApm}[1]{\sCatArrowsSymbolpm[#1]}

\newcommand{\ttx}{\mathtt{x}}

\DeclareMathOperator{\symbolhom}{\mathsf{hom}}

\newcommand{\functors}[2]{
\mathtt{Fun}(#1,#2)
}

\newcommand{\cathom}[3]{
\symbolhom_{#1}(#2,#3)
}

\newcommand{\hmult}{\circleddash}

\newcommand{\mbf}[1]{
\boldsymbol{#1}
}

\DeclareMathOperator{\twinsymbol}{\boldsymbol{\mathtt{T}}}
\newcommand{\twin}[1]{
\twinsymbol({#1})
}

\newcommand{\arrowws}[3]{
{#1}{\,\,}
\stackon[-3.0pt]{$\xrightarrow{\hspace*{0.7cm}}$}{$#3$}{{\,}#2}
}

\newcommand{\twoarrowws}[5]{
{#1}{\,\,}
\stackon[-4.0pt]{$\xrightarrow{\hspace*{0.5cm}}$}{$#4$}
{{\,}#2}
\,\,
\stackon[-4.0pt]{$\xrightarrow{\hspace*{0.5cm}}$}{$#5$}
{{\,}#3}
}

\newcommand{\ch}[2]{
\mathtt{hom}_{\scCat}(#1,#2)
}

\newcommand{\mybar}[3]{%
    \mathrlap{\hspace{#2}\overline{\scalebox{#1}[1]{\phantom{\ensuremath{#3}}}}}\ensuremath{#3}
}

\DeclareMathOperator{\erl}
{\mybar{0.9}{-0.1pt}{\RR}}

\DeclareMathOperator{\scOne}{\mathbf{1}}
\DeclareMathOperator{\scTwo}{\mathbf{2}}
\DeclareMathOperator{\scThree}{\mathbf{3}}
\DeclareMathOperator{\scSQ}{\mathbf{Q}}
\DeclareMathOperator{\scMatrices}{\mathbf{Matr}}

\newcommand{\dueuno}{\begin{bmatrix}
1\\
0\\ 
\end{bmatrix}}

\newcommand{\trequattro}{\begin{bmatrix}
1&0&0&0\\
0&1&0&0\\
0&0&1&0 
\end{bmatrix}}

\newcommand{\quattrouno}{\begin{bmatrix}
1\\
0\\
0\\
0\\ 
\end{bmatrix}}

\newcommand{\tredue}{\begin{bmatrix}
1&0\\
0&1\\
0&0 
\end{bmatrix}}

\newcommand{\treuno}{\begin{bmatrix}
1\\
0\\
0 
\end{bmatrix}}

\DeclareMathOperator{\runitssymbol}{\texttt{U}}

\newcommand{\runits}[1]{\runitssymbol_{#1}}

\newcommand{\aidentitydueuno}{\aidentity{2}{1}}
\newcommand{\aidentitytredue}{\aidentity{3}{2}}
\newcommand{\aidentitytreuno}{\aidentity{3}{1}}
\newcommand{\aidentityquattrouno}{\aidentity{4}{1}}
\newcommand{\aidentitytrequattro}{\aidentity{3}{4}}

\newcommand{\aidentity}[2]{
\matrixone^{\scalebox{0.5}{#1}}_{\scalebox{0.5}{#2}}
}

\DeclareMathOperator{\pseudotopol}{\mathtt{N}}

\DeclareMathOperator{\KKint}{
\mathtt{I}
}

\DeclareMathOperator{\wlocname}{\mathsf{wloc}}

\newcommand{\mrep}[1]{
\Measurable[#1]
}

\DeclareMathOperator{\cotaNAME}{\mathsf{c}}
\newcommand{\cota}[2]{
\cotaNAME_{#1}(#2)
}

\DeclareMathOperator{\fataNAME}{\mathsf{s}}
\newcommand{\fata}[2]{
\fataNAME_{#1}(#2)
}

\DeclareMathOperator{\afccname}{\mathsf{FCC}}
\newcommand{\afcc}[1]{
\afccname({#1})
}

\DeclareMathOperator{\Lebtrealsymbol}{\mathcal{L}}
\DeclareMathOperator{\Lebtclasssymbol}{\mathsf{L}}

\DeclareMathOperator{\Lebtb}{{\Lebtclasssymbol}}
\DeclareMathOperator{\Lebt}{\Lebtrealsymbol}

\newcommand{\shom}[2]{
\hset{#1}{#2}
}

\newcommand{\fccp}[2]{
{#2}\times\hset{#1}{#2}
}

\newcommand{\nseqofsets}{\boldsymbol{q}}
\newcommand{\seqofsets}{\boldsymbol{\bsubset}}

\newcommand{\classoffunc}{\boldsymbol{\mathcal{R}}}
\newcommand{\fccf}{\varphi}
\newcommand{\fccg}{\beta}
\newcommand{\fcch}{\gamma}
\newcommand{\fcc}{\mathcal{F}}

\newcommand{\limv}{y}

\newcommand{\dePossel}{\mathcal{V}}

\newcommand{\wloc}[1]{
\wlocname(#1)
}

\DeclareMathOperator{\cofname}{\mathsf{cof}}
\newcommand{\cofi}[1]{
\cofname(#1)
}

\DeclareMathOperator{\areatoname}{Q}
\newcommand{\areato}[1]{
\areatoname(#1)
}

\newcommand{\parama}{
{\boldsymbol{r}}
}

\newcommand{\paramb}{
{\boldsymbol{s}}
}

\newcommand{\paramc}{
{\boldsymbol{p}}
}

\newcommand{\enum}{
{\boldsymbol{\alpha}}
}

\DeclareMathOperator{\directions}{dir}
\newcommand{\dirStolz}{
\invim{\mathsf{Fin}}{\mathsf{S}}
}

\DeclareMathOperator{\finalfiltername}{\mathsf{Fin}}
\newcommand{\fiseof}[1]{
\finalfiltername[{#1}]}

\DeclareMathOperator{\pointw}{w}
\DeclareMathOperator{\pointwbold}{\boldsymbol{\pointw}}
\DeclareMathOperator{\pointz}{z}
\DeclareMathOperator{\pointzbold}{\boldsymbol{\pointz}}

\newcommand{\figebyin}[2]{
{#2}^{\uparrow}\!\left[#1\right]
}

\DeclareMathOperator{\inteop}{
\mathring{\pseudotopol}
}

\DeclareMathOperator{\ds}{
\fmm}

\DeclareMathOperator{\newds}{
D}

\DeclareMathOperator{\sihomsymbol}{
{\ensuremath{\normalfont\boldsymbol{\mathsf{twinh}}}}}

\newcommand{\sihom}[3]{\sihomsymbol_{#1}(#2,#3)}

\newcommand{\fctpm}[1]{
\sCatApm{#1}
}
\newcommand{\sfromctpm}{\sCatArrowsSymbolpm}

\newcommand{\fpmtc}[1]{
\snCat{cat}(#1)
}

\newcommand{\snCat}[1]{\mathsf{#1}}

\DeclareMathOperator{\sCat}{{\mathsf{R}}}
\DeclareMathOperator{\scCat}{{\mathsf{Cat}}}
\DeclareMathOperator{\sCatone}{{\mathsf{R}}}
\DeclareMathOperator{\sCattwo}{{\mathsf{S}}}
\DeclareMathOperator{\sCatthree}{{\mathsf{U}}}
\DeclareMathOperator{\sFunctor}{\mathsf{J}}
\DeclareMathOperator{\sTwin}{\mathsf{Twin}}
\DeclareMathOperator{\sSquareF}{\square}

\DeclareMathOperator{\sFunctorone}{\mathsf{J}}
\DeclareMathOperator{\sFunctortwo}{\mathsf{K}}
\DeclareMathOperator{\sFunctorthree}{\mathsf{L}}

\newcommand{\naturalto}{%
  \mathrel{\vbox{\offinterlineskip
    \mathsurround=0pt
    \ialign{\hfil##\hfil\cr
      \normalfont\scalebox{1.8}{.}\cr
      $\longrightarrow$\cr}
  }}%
}

\DeclareMathOperator{\ffSetsymbol}{\mathsf{Set}}
\newcommand{\ffSet}[1]{{#1}_{\ffSetsymbol}}
\DeclareMathOperator{\snt}{\naturalto}
\DeclareMathOperator{\sPM}{{\mathsf{PM}}}
\DeclareMathOperator{\sUPM}{{\mathsf{uPM}}}
\DeclareMathOperator{\sRPM}{{\mathsf{RPM}}}
\DeclareMathOperator{\sSet}{{\mathsf{Set}}}
\DeclareMathOperator{\sRing}{{\mathsf{Ring}}}
\DeclareMathOperator{\spSet}{\mathtt{PSet}}
\DeclareMathOperator{\srSet}{\mathtt{Rel}}

\DeclareMathOperator{\sTopCH}{\mathtt{CH}}

\DeclareMathOperator{\sBoole}{\mathtt{BA}}

\DeclareMathOperator{\sDset}{\mathtt{DSet}}
\DeclareMathOperator{\sFset}{\mathtt{FSet}}

\DeclareMathOperator{\tailname}{\mathtt{tail}}
\newcommand{\newtail}[2]{
\tailname_{#2}\!\left(#1\right)}

\newcommand{\fisof}{\mathcal{W}}

\DeclareMathOperator{\nameinttop}{{\mathcal{I}}}

\DeclareMathOperator{\homomorphisms}{\mathtt{hom}}

\DeclareMathOperator{\homomorphismsSet}{
\homomorphisms_{\sSet}}

\DeclareMathOperator{\homomorphismspSet}{
\homomorphisms_{\spSet}}

\DeclareMathOperator{\homomorphismsB}{
\homomorphisms_{\mathtt{BA}}}

\DeclareMathOperator{\homomorphismsTop}{
\homomorphisms_{\mathtt{Top}}}

\newcommand{\hset}[2]{
\homomorphismsSet(#1,#2)}

\newcommand{\hpset}[2]{
\homomorphismspSet(#1,#2)}

\DeclareMathOperator{\homomorphismsSetr}{
\homomorphisms_{\srSet}}

\newcommand{\hrset}[2]{
\homomorphismsSetr(#1,#2)}

\newcommand{\hsetnz}[2]{
\homomorphisms^{\ast}_{\sSet}(#1,#2)}
\newcommand{\hBoole}[2]{
\homomorphismsB(#1,#2)}
\newcommand{\hBoolenz}[2]{
\homomorphisms^{\ast}_{\sBoole}(#1,#2)}
\newcommand{\htop}[2]{
\homomorphismsTop(#1,#2)}

\DeclareMathOperator{\nameBAC}{\boldsymbol{\mathtt{Boole}}}

\newcommand{\Boolehnz}[2]{
\nameBAC_{\ast}({#1},{#2})
}

\newcommand{\inttop}[1]{
{\nameinttop}_{#1}}

\newcommand{\inttopset}[2]{
\nameinttop_{#1}(#2)
}

\newcommand{\identifyf}[1]{
{\boldsymbol{\mathsf{I}}}_{#1}
}

\newcommand{\sofm}{
{\Measurable}_{\hmeas}^{\mathsf{F}}
}

\newcommand{\ns}{
\Measurable_{0}
}

\newcommand{\np}{
\pi_{\hmeas}
}

\newcommand{\inp}{
{\rho}_{\hmeas}
}

\newcommand{\complemento}[1]{
\complement{#1}
}

\newcommand{\standard}{(\mm,\Measurable,\hmeas)}

\newcommand{\sdx}{Q}
\newcommand{\sdxb}{R}

\DeclareMathOperator{\matrixone}{\scriptstyle{\mathsf{M}}}
\DeclareMathOperator{\matrixtwo}{\scriptstyle{\mathsf{N}}}

\DeclareMathOperator{\ione}{\scriptstyle{\mathsf{I}}}
\DeclareMathOperator{\ii}{\scriptstyle{\mathsf{II}}}
\DeclareMathOperator{\iii}{\scriptstyle{\mathsf{III}}}
\DeclareMathOperator{\iv}{\scriptstyle{\mathsf{IV}}}

\newcommand{\unit}[1]{
{\mathsf{1}}_{#1}
}

\newcommand{\indic}[1]{
{\textup{I}}_{#1}
}

\newcommand{\herz}[4]{
{#2}_{#3}^{#4}(#1)
}
\newcommand{\herzNP}[3]{
{#1}_{#2}^{#3}
}

\newcommand{\dirset}{\fmm}

\newcommand{\sffm}{\boldsymbol{\diamond}}

\newcommand{\fdirimFS}[2]{
{#1}_{\sffm}(#2)
}

\newcommand{\fdirimF}[1]{
{#1}_{\sffm}}

\newcommand{\dirimf}[1]{
{#1}_{\ast}}

\newcommand{\dirim}[2]{
{#1}_{\ast}(#2)}
\newcommand{\invim}[2]{
{#1}^{\ast}(#2)}
\DeclareMathOperator{\preorder}{
\mathtt{R}}

\DeclareMathOperator{\nametailsof}{\mathtt{Tail}}

\newcommand{\trailof}[2]{
{#1}_{\bullet}[#2]
}

\newcommand{\tailsof}[2]{
\nametailsof_{#1}[#2]
}

\newcommand{\gseq}{g-sequence}
\DeclareMathOperator{\namefgb}{{\mathsf{T}}}
\DeclareMathOperator{\namefgbpoints}{{\mathsf{t}}}
\newcommand{\fgb}[2]{
\namefgb_{#2}[#1]}
\newcommand{\fgbp}[2]{
\namefgbpoints_{#2}[#1]}

\DeclareMathOperator{\namegsgb}{s}
\newcommand{\gsgb}[1]{
\namegsgb_{#1}}
\newcommand{\cti}[2]{{\mathcal{T}}_{#1}(#2)}

\DeclareMathOperator*{\fcclim}{
{\boldsymbol{\mathpzc{lim}}}}
\DeclareMathOperator*{\flimsup}{
{\boldsymbol{\mathsf{limsup}}}}
\DeclareMathOperator*{\fliminf}{
{\boldsymbol{\mathsf{liminf}}}}
\newcommand{\clusters}[2]{
{\boldsymbol{\mathsf{Cluster}}}[#1,#2]}
\newcommand{\pclusters}[2]{
{\boldsymbol{\mathsf{clusterset}}}(#1,#2)}

\newcommand{\gsclusters}[2]{%
\mathtt{ClusterSet}(#1,#2)}
\DeclareMathOperator*{\flim}{
{\boldsymbol{{\mathsf{lim}}}}}
\DeclareMathOperator*{\gslim}{
{\boldsymbol{\mathtt{lim}}\,}}

\newcommand{\gsaccp}[2]{
\clusters{#1}{#2}}

\DeclareMathOperator{\bubicn}{\tboundary{\ubicn}}
\newcommand{\aperture}{j}

\DeclareMathOperator{\ubicn}{
\mathbb{B}}
\DeclareMathOperator{\eventuallydisjoint}{
(eventually \, disjoint)}

\DeclareMathOperator{\frequentlyoutsideang}{
(frequently \, outside \, the \, angular \, filter)}
\DeclareMathOperator{\aproperty}{
(asymptotic)}
\DeclareMathOperator{\Fatouproperty}{
(Fatou)}
\DeclareMathOperator{\rotationalinvariance}{
(rotational \, invariance)}
\DeclareMathOperator{\PrivalovSetSymbol}{\mathsf{Privalov}}
\newcommand{\PrivalovSet}[1]{\PrivalovSetSymbol(#1)}
\DeclareMathOperator{\LusinSetSymbol}{\mathsf{Lusin}}
\newcommand{\LusinSet}[2]{\LusinSetSymbol_{#2}(#1)}
\newcommand{\LusinSetni}[1]{\LusinSetSymbol(#1)}
\newcommand{\limvalue}{\xi}
\newcommand{\shadownawithp}[1]{
\overline{(#1)}}
\newcommand{\shadow}[2]{
\shadowna{#1}(#2)}
\newcommand{\shadowwithp}[2]{
\shadownawithp{#1}(#2)}
\newcommand{\shadowna}[1]{
\overline{#1}}
\newcommand{\shadowset}[2]{
\shadowna{#1}[#2]}

\makeatletter
\newsavebox\myboxA
\newsavebox\myboxB
\newlength\mylenA
\newcommand*\xoverline[2][0.75]{%
    \sbox{\myboxA}{$\m@th#2$}%
    \setbox\myboxB\null
    \ht\myboxB=\ht\myboxA%
    \dp\myboxB=\dp\myboxA%
    \wd\myboxB=#1\wd\myboxA
    \sbox\myboxB{$\m@th\overline{\copy\myboxB}$}
    \setlength\mylenA{\the\wd\myboxA}
    \addtolength\mylenA{-\the\wd\myboxB}%
    \ifdim\wd\myboxB<\wd\myboxA%
       \rlap{\hskip 0.5\mylenA\usebox\myboxB}{\usebox\myboxA}%
    \else
        \hskip -0.5\mylenA\rlap{\usebox\myboxA}{\hskip 0.5\mylenA\usebox\myboxB}%
    \fi}
\makeatother

\makeatletter
\let\save@mathaccent\mathaccent
\newcommand*\if@single[3]{%
  \setbox0\hbox{${\mathaccent"0362{#1}}^H$}%
  \setbox2\hbox{${\mathaccent"0362{\kern0pt#1}}^H$}%
  \ifdim\ht0=\ht2 #3\else #2\fi
  }
\newcommand*\rel@kern[1]{\kern#1\dimexpr\macc@kerna}
\newcommand*\widebar[1]{\@ifnextchar^{{\wide@bar{#1}{0}}}{\wide@bar{#1}{1}}}
\newcommand*\wide@bar[2]{\if@single{#1}{\wide@bar@{#1}{#2}{1}}{\wide@bar@{#1}{#2}{2}}}
\newcommand*\wide@bar@[3]{%
  \begingroup
  \def\mathaccent##1##2{%
    \let\mathaccent\save@mathaccent
    \if#32 \let\macc@nucleus\first@char \fi
    \setbox\z@\hbox{$\macc@style{\macc@nucleus}_{}$}%
    \setbox\tw@\hbox{$\macc@style{\macc@nucleus}{}_{}$}%
    \dimen@\wd\tw@
    \advance\dimen@-\wd\z@
    \divide\dimen@ 3
    \@tempdima\wd\tw@
    \advance\@tempdima-\scriptspace
    \divide\@tempdima 10
    \advance\dimen@-\@tempdima
    \ifdim\dimen@>\z@ \dimen@0pt\fi
    \rel@kern{0.6}\kern-\dimen@
    \if#31
      \overline{\rel@kern{-0.6}\kern\dimen@\macc@nucleus\rel@kern{0.4}\kern\dimen@}%
      \advance\dimen@0.4\dimexpr\macc@kerna
      \let\final@kern#2%
      \ifdim\dimen@<\z@ \let\final@kern1\fi
      \if\final@kern1 \kern-\dimen@\fi
    \else
      \overline{\rel@kern{-0.6}\kern\dimen@#1}%
    \fi
  }%
  \macc@depth\@ne
  \let\math@bgroup\@empty \let\math@egroup\macc@set@skewchar
  \mathsurround\z@ \frozen@everymath{\mathgroup\macc@group\relax}%
  \macc@set@skewchar\relax
  \let\mathaccentV\macc@nested@a
  \if#31
    \macc@nested@a\relax111{#1}%
  \else
    \def\gobble@till@marker##1\endmarker{}%
    \futurelet\first@char\gobble@till@marker#1\endmarker
    \ifcat\noexpand\first@char A\else
      \def\first@char{}%
    \fi
    \macc@nested@a\relax111{\first@char}%
  \fi
  \endgroup
}
\makeatother

\DeclareMathOperator{\Poisson}{P}
\newcommand{\Poissonf}[1]{
\Poisson#1}
\newcommand{\realbf}{f}
\newcommand{\realbftwo}{g}

\newcommand{\classbf}{\boldsymbol{f}}
\newcommand{\bpoint}{x}
\newcommand{\bpointtwo}{r}
\newcommand{\bpointthree}{z}
\newcommand{\bsubset}{Q}
\newcommand{\bsubsettwo}{R} 
\newcommand{\bsubsetthree}{S} 

\newcommand\FatouSet[1]{
\Fatou{(#1)}}
\DeclareMathOperator{\Fatouname}{\mathsf{Fatou}}
\DeclareMathOperator{\Fatou}{\Fatouname}
\newcommand\relFatouSet[2]{
\Fatou(#1;#2)}
\DeclareMathOperator{\Har}{H}
\newcommand{\hmeaso}[1]{\hmeasoname#1}
\newcommand{\hmeasoP}[1]{\hmeasoname(#1)}

\newcommand{\normalmeasureP}[1]{\hmeas{(#1)}}
\DeclareMathOperator{\adbasis}{\boldsymbol{\mathcal{Z}}}
\DeclareMathOperator{\radius}{\rho}
\newcommand{\mofoar}[1]{
\displaystyle{\sup_{#1}}}
\newcommand{\mofoartextstyle}[1]{
\textstyle{\sup_{#1}}
}

\newcommand{\mofoarftextstyle}[2]{
\textstyle{\sup_{#1}{#2}}
}

\newcommand{\mofoarf}[2]{
\mofoar{#1}\,{#2}
}

\newcommand{\mofoarfa}[3]{
\mofoarf{#1}{#2}(#3)
}

\newcommand{\mofoarfaWP}[3]{
\mofoar{#1}{#2}\,(#3)
}

\DeclareMathOperator{\fofibox}{\mathsf{g}}
\DeclareMathOperator{\bofofibox}{{\fofibox}}
\newcommand{\fofibop}[1]{\fofibox(#1)}

\DeclareMathOperator{\gs}{
\mathscr{w}
}

\DeclareMathOperator{\gsb}{
\mathscr{t}
}

\newcommand{\sequencenr}[1]{{\{{#1}_{\aperture}\}}_{{\aperture}}}

\newcommand{\sequence}[1]{
{
\{
{#1}_{\aperture}
\}
}_{\aperture\geq 1}
}

\DeclareMathOperator{\newpartition}{\mathsf{C}}

\DeclareMathOperator{\topologyname}{
\boldsymbol{\mathcal{G}}}

\DeclareMathOperator{\topol}{
\Phi}

\DeclareMathOperator{\topoltwo}{
\Theta}

\newcommand{\topology}[1]{\topologyname(#1)}
\newcommand{\topologyNE}[1]{\topologyname_{\bullet}(#1)}
\newcommand{\topologyfromametric}[2]{
\topologyname(#1,#2)}

\DeclareMathOperator{\LebesgueNamesymbol}{\mathsf{Lebesgue}}
\DeclareMathOperator{\LebesgueName}{
\LebesgueNamesymbol}
\newcommand{\LebesgueSet}[1]{
\LebesgueName[#1]}

\newcommand{\constantsymbol}{c}

\newcommand{\constant}{
\boldsymbol{\constantsymbol}
}

\DeclareMathOperator{\Stolzsymbol}{\mathsf{Stolz}}

\newcommand{\Stolztheta}{
\Stolzsymbol_{\bpoint}
}

\DeclareMathOperator{\symmetric}{s}

\newcommand{\Stolzthetasym}{
\Stolzsymbol_{\bpoint}^{\symmetric}
}

\DeclareMathOperator{\Plessnersymbol}{\mathsf{Plessner}}

\newcommand{\PlessnerSet}[1]{\Plessnersymbol(#1)}

\DeclareMathOperator{\realPlessnersymbol}{\mathsf{realPlessner}}

\newcommand{\realPlessnerSet}[1]{\realPlessnersymbol
[#1]}

\newcommand{\boldav}{
\boldsymbol{|}}

\makeatletter
\newcommand*\dashline{\rotatebox[origin=c]{90}{$\dabar@\dabar@$}}
\makeatother

\newcommand{\ballsopenorclosed}{\mathbf{B}}

\newcommand{\closedballs}[1]{
\ballsopenorclosed
[{#1}]
}

\DeclareMathOperator{\symbolforball}{{B}}

\newcommand{\openballcr}[2]{
\symbolforball_{\mm,\mbdr}
\!
\left(#1,#2\right)
}
\newcommand{\openballcrs}[2]{
\symbolforball_{\mm}
\!
\left(#1,#2\right)
}

\newcommand{\openballcrm}[2]{
\symbolforball_{\mbdr}
\!
\left(#1,#2\right)
}

\newcommand{\setone}{
{\bsubset_{1}}
}

\newcommand{\settwo}{
{\bsubset_{2}}
}

\newcommand{\averagenaf}[1]{
{#1}_{\!\!\hmeas}
}

\DeclareMathOperator{\averageoperator}{
{\widetilde{\hmeas}}
}

\newcommand{\apairing}[2]{
\averageoperator\left(#1|#2\right)
}

\newcommand{\newapairing}[2]{
\hmeas\left(#1|#2\right)
}

\newcommand{\cof}[1]{
\widetilde{#1}
}

\DeclareMathOperator{\BorelName}{\boldsymbol{\mathcal{B}}}
\newcommand{\Borelsets}[1]{
\BorelName(#1)}

\DeclareMathOperator{\frepresname}{R}

\DeclareMathOperator{\frepres}{\frepresname}

\newcommand{\frepresf}[1]{
\frepres(#1)
}

\DeclareMathOperator{\Measurable}{\mathcal{M}}

\DeclareMathOperator{\symbolforpowerset}{{\mathcal{P}}}

\newcommand{\totalpowerset}[1]{
\symbolforpowerset(#1)
}

\newcommand{\powersetnotempty}[1]{
\symbolforpowerset_{\bullet}(#1)
}

\DeclareMathOperator{\bqsubset}{
S}

\newcommand{\levelsymbol}{r}
\newcommand{\level}{\levelsymbol}

\newcommand{\artctabp}[4]{
#4\supset#3\supset#1\to#2
}

\DeclareMathOperator{\edisjointNAME}{
\triangleleft
\triangleright
}

\newcommand{\edisjoint}[2]{
{#1}\edisjointNAME{#2}
}

\makeatletter
\newcommand{\oset}[3][0ex]{%
  \mathrel{\mathop{#3}\limits^{
    \vbox to#1{\kern-2\ex@
    \hbox{$\scriptstyle#2$}\vss}}}}
\makeatother

\newcommand{\arequiv}[2]{
#1\sim_{\!\!\!\!\!\bpoint}\,#2
}

\DeclareMathOperator{\one}{\boldsymbol{(1)}}

\DeclareMathOperator{\two}{\boldsymbol{(2)}}

\DeclareMathOperator{\rms}{\boldsymbol{\mathrm{(S)}}}

\DeclareMathOperator{\sone}{\boldsymbol{\mathrm{(S.1)}}}
\DeclareMathOperator{\stwo}{\boldsymbol{\mathrm{(S.2)}}}
\DeclareMathOperator{\sthree}{\boldsymbol{\mathrm{(S.3)}}}
\DeclareMathOperator{\sfour}{\boldsymbol{\mathrm{(S.4)}}}
\DeclareMathOperator{\sfive}{\boldsymbol{\mathrm{(S.5)}}}

\DeclareMathOperator{\rmcg}{\boldsymbol{\mathrm{(CG)}}}

\DeclareMathOperator{\cgone}{\boldsymbol{\mathrm{(CG.1)}}}
\DeclareMathOperator{\cgtwo}{\boldsymbol{\mathrm{(CG.2)}}}

\DeclareMathOperator{\atoms}{
{{\mathsf{Atoms}}_{\hmeas}
}}

\newcommand{\saequiv}{
\oset{\,\hmeas}{=}
}

\newcommand{\aequiv}[2]{
#1\!\saequiv#2
}

\newcommand{\aesubset}[2]{
#1\subset_{\hmeas}#2
}

\DeclareMathOperator{\scollection}{\mathcal{C}}
\DeclareMathOperator{\scollectiontwo}{
\mathcal{D}
}

\DeclareMathOperator{\bnewpartition}{\boldsymbol{\mathsf{c}}}

\DeclareMathOperator{\newpartitiontwo}{\mathsf{D}}

\DeclareMathOperator{\mtclosureName}{
{CL}
}

\newcommand{\mtclosureof}[2]{
\mtclosureName_{#2}^{\hmeas}\{#1\}
}

\newcommand{\aalgebraof}[1]{
\boldsymbol{\mathsf{A}}(#1)
}

\newcommand{\aringof}[1]{
\boldsymbol{\mathsf{R}}(#1)
}

\newcommand{\ringof}[1]{
#1^*
}

\newcommand{\tilesofp}{
\mathcal{T}
}

\DeclareMathOperator{\eSete}{E}

\DeclareMathOperator{\dsetd}{D}

\DeclareMathOperator{\aregion}{V}

\DeclareMathOperator{\aregioneso}{E}

\DeclareMathOperator{\fmm}{
{A}
}

\DeclareMathOperator{\mm}{
{X}
}

\DeclareMathOperator{\smm}{
\boldsymbol{\mm}
}

\DeclareMathOperator{\mmtwo}{
{Y}
}

\DeclareMathOperator{\conjugate}{Q}

\newcommand{\conjugatef}[1]{
\conjugate{\!}#1
}

\newcommand{\closedsets}[1]{
\boldsymbol{\mathcal{F}}(#1)
}

\newcommand{\openballs}[1]{
\ballsopenorclosed(#1)
}

\DeclareMathOperator{\Ff}{
\mathsf{f}\mathbb{N}}

\DeclareMathOperator{\newafilter}{
\mathsf{Z}}

\DeclareMathOperator{\newbfilter}{
\mathsf{W}}

\DeclareMathOperator{\newcfilter}{
\mathsf{Y}}

\DeclareMathOperator{\newdfilter}{
\mathsf{V}}

\newcommand{\eoa}{
\bpoint}

\DeclareMathOperator{\fm}{
\mathsf{b}}

\DeclareMathOperator{\fmb}{
\mathsf{c}}

\DeclareMathOperator{\fmc}{
\mathsf{d}}

\DeclareMathOperator{\fmd}{
\mathsf{e}}

\DeclareMathOperator{\afilter}{
\mathsf{F}}

\DeclareMathOperator{\ambient}{
\mathbf{W}
}

\newcommand{\angularbv}[2]{
{#1}_{\flat}(#2)}

\newcommand{\angularbvwithp}[2]{
{(#1)}_{\flat}(#2)}

\newcommand{\angularbvna}[1]{
{#1}_{\flat}}

\DeclareMathOperator{\collaregions}{{
{\boldsymbol{\Lambda}}}
}

\newcommand{\newfilter}[2]{
{[#1]}_#2
}

\newcommand{\newfilterWP}[1]{
[#1]
}

\DeclareMathOperator{\Setone}{\domain}

\DeclareMathOperator{\ntambient}{W}

\newcommand{\bsubsetfour}{T}

\newcommand{\tops}[2]{
(#1,#2)
}

\DeclareMathOperator{\ts}{
\mathbf{Y}
}

\DeclareMathOperator{\tsnt}{
{Y}
}

\DeclareMathOperator{\tsnttwo}{
{W}
}

\DeclareMathOperator{\secondts}{
\mathbf{S}
}

\newcommand{\newrestriction}[2]{#1|_{#2}
}

\DeclareMathOperator{\imbedding}{\hookrightarrow}

\DeclareMathOperator{\interval}{
\texttt{int}
}

\DeclareMathOperator{\metricambientNAME}{
{d}
}

\DeclareMathOperator{\macn}{
\boldsymbol{\metricambientNAME}
}

\DeclareMathOperator{\mbdr}{\delta}

\DeclareMathOperator{\indexset}{
{I}
}

\DeclareMathOperator{\adomain}{
\domain
}

\newcommand{\vvnclosure}[2][3]{%
  {}\mkern#1mu\overline{\mkern-#1mu#2}}

\newcommand{\parrow}
{
{
\,\,\normalfont\scalebox{1.3}{$\rightarrowtail$}
\,
}
}

\newcommand{\newclosure}[2][3]{
{}\mkern#1mu\overline{\mkern-#1mu#2}}

\newcommand{\curva}{\boldsymbol{c}}

\newcommand{\shrinksang}{\xrightarrow[nt]{}}

\newcommand{\meanvalue}[2]{
\frac{1}{\absv{#2}}
\int_{#2}\classbf(e^{i\theta})\,d\theta
}

\newcommand{\ameanvalue}[2]{
\frac{1}{\normalmeasureP{#2}}
\int_{#2}#1{\,}d\!\hmeas
}

\newcommand{\absv}[1]{\left|#1\right|}

\newcommand{\holomorphic}[1]{\mathcal{O}(#1)}
\DeclareMathOperator{\complexmeasuresName}{\mathcal{M}}

\newcommand{\cmeasuresh}[1]{{\complexmeasuresName}(#1,\hmeas)}

\newcommand{\setofsuchthat}[2]{
\left\{
#1\!\!:
#2
\right\}
}

\newcommand{\bootstrap}[1]{
{#1}^{\sharp}
}

\DeclareMathOperator{\outermeasurename}{\textsf{Outer}}
\newcommand{\outermeasure}[1]{\outermeasurename[{#1}]}

\DeclareMathOperator{\cfotc}{
\continuous{
\vvnclosure[0]{\adomain}}
}

\DeclareMathOperator{\hdist}{\mbdr_H}

\DeclareMathOperator{\asecondfilter}{\mathsf{\Psi}}

\DeclareMathOperator{\filterbase}{\Psi}

\newcommand{\nf}[2]{
{\mathcal{N}}_{#1}(#2)
}

\DeclareMathOperator{\newnfSPname}{\mathsf{N}}

\newcommand{\nsdi}[2]{
{\newnfSPname}_{#1}(#2)
}

\newcommand{\newnfSP}[3]{
{\newnfSPname}^{#2}(#3)
}

\newcommand{\newnfT}[1]{
{\newnfSPname}_{#1}
}

\DeclareMathOperator{\absdistributionfunction}{
\lambda
}

\newcommand{\newdistributionfunction}[2]{
\absdistributionfunction(#1\dashline\!#2)
}

\newcommand{\newdistributionfunctionA}[3]{
\newdistributionfunction{#1}{#2}(#3)
}

\DeclareMathOperator{\tauto}{\textsf{v}}

\newcommand{\pbpoint}[2]{
\bpoint(#2)
}

\newcommand{\rbpoint}[2]{
r(#2)
}

\newcommand{\diameterwrtmetric}[2]{
\text{diam}_{#2}\left[#1\right]
}

\DeclareMathOperator{\uspacename}{\mathbb{E}}

\newcommand{\uhssn}{
\uspacename
}

\newcommand{\omib}[1]{
{#1}^{\circ}
}

\newcommand{\omibs}[2]{
\omib{#1}{#2}{\,}
}

\DeclareMathOperator{\cofunctions}{
\mathcal{H}}

\DeclareMathOperator{\coafoasname}{\boldsymbol{\mathsf{F}}}

\DeclareMathOperator{\ucoafoasname}{
\boldsymbol{\mathscr{U}}}

\newcommand{\soaf}[1]{
{\spaceofallfilters{#1}}}

\newcommand{\usoaf}[1]{
\ucoafoasname(#1)}

\newcommand{\soaftcags}[2]{
\coafoasname_{#2}(#1)}

\newcommand{\spaceofallfilters}[1]{
\boldsymbol{\mathscr{F}}(#1)}

\DeclareMathOperator{\namespags}{
\mathscr{S}
}

\DeclareMathOperator{\namespagsF}{
\boldsymbol{\circ}}

\newcommand{\sspags}[1]{
\namespags{(#1)}
}

\newcommand{\sspagsF}[1]{
{#1}_{\namespagsF}
}

\newcommand{\spags}[1]{
\namespags(\powersetnotempty{#1})
}

\newcommand{\moadb}[1]{
\mathcal{M}_{#1}
}

\newcommand{\moadbf}[2]{
\moadb{#1}{#2}
}

\newcommand{\moadbfa}[3]{
\moadbf{#1}{#2}(#3)
}

\DeclareMathOperator{\mosht}{M}
\newcommand{\moshtf}[1]{\mosht_{#1}}
\newcommand{\moshtfa}[2]{\mosht_{#1}({#2})}
\newcommand{\moshtfafoar}[3]{\mosht_{#1,#3}({#2})}
\newcommand{\moshtffoar}[2]{\mosht_{#1,#2}}

\DeclareMathOperator{\cmosht}{m}
\newcommand{\cmoshtf}[1]{\cmosht_{#1}}
\newcommand{\cmoshtfa}[2]{\cmosht_{#1}({#2})}

\newcommand{\mitbotud}[1]{
|{#1}|
}

\DeclareMathOperator{\bpath}{w}

\newcommand{\sotpvf}[1]{[0,+\infty]^{#1}}

\newcommand{\dsubset}{Z}

\DeclareMathOperator{\radiusAR}{R}

\newcommand{\radialbvf}[1]{{#1}_{\radius}}
\newcommand{\radialbvfa}[2]{
\radialbvf{#1}(#2)
}
\newcommand{\radialbvfaAR}[2]{
{#1}_{\radiusAR}(#2)
}

\newcommand{\blim}[3]{
\lim_{#2\ni\dpoint\to#3}#1(\dpoint)
}

\newcommand{\blimDT}[3]{
\lim_{#2\ni\bsubset\to#3}#1
}

\DeclareMathOperator{\NN}{\mathbb{N}}

\def\QQ{{\mathbb Q}}

\def\RR{{\mathbb R}}
\def\CC{{\mathbb C}}
\def\ZZ{{\mathbb Z}}

\DeclareMathOperator{\cintervals}{\boldsymbol{\mathcal{I}}}

\DeclareMathOperator{\triangolo}{\textsc{t}}

\DeclareMathOperator{\tenda}{\Delta}
\newcommand{\tent}[2]{
\tenda^{#1}#2
}

\newcommand{\tentna}[1]{
\tenda^{#1}
}

\DeclareMathOperator{\area}{S}

\newcommand{\areaf}[3]{
\area_{#2}#1(#3)}

\newcommand{\areafwithheight}[4]{
\area_{#2,#4}#1(#3)}

\newcommand{\areafna}[2]{
\area_{#2}#1}

\newcommand{\areafnawithheight}[3]{
\area_{#2,#3}#1}

\newcommand{\areafnawithpwithheight}[3]{
\area_{#2,#3}(#1)}

\DeclareMathOperator{\foaregions}{\varphi}

\DeclareMathOperator{\foaregionstan}{\tau}

\DeclareMathOperator{\hmeas}{\omega}
\DeclareMathOperator{\hmeasoname}{\hmeas^{*}}

\DeclareMathOperator{\hsigma}{
{\widehat{\BorelName}}_{\hmeas}
}

\DeclareMathOperator{\Sm}{{\mathcal{S}}}

\DeclareMathOperator{\ZygmundName}{{\mathcal{A}}}

\newcommand{\Zm}{\ZygmundName}

\newcommand{\ZygmundNM}[1]{
\ZygmundName^{*}
({#1})
}

\newcommand{\tinterior}[1]{
{#1}^{\circ}
}

\newcommand{\ginterior}[1]{
\ZygmundName'({#1})
}

\DeclareMathOperator{\udone}{\mathbb{D}}
\DeclareMathOperator{\budone}{\partial\mathbb{D}}

\newcommand{\hmeasure}[1]{\hmeas_{#1}}

\newcommand{\tboundary}[1]{\partial{#1}}

\newcommand{\trboundary}[2]{
\partial_{#2}{(#1)}}

\newcommand{\trboundaryWP}[2]{\partial_{#2}{#1}}

\newcommand{\domain}{\Omega}

\newcommand{\baccdomain}{\partial_{a}\domain}

\newcommand{\bdomain}{\partial\domain}

\newcommand{\classbftwo}{g}

\newcommand{\dpoint}{z}
\newcommand{\dfunction}{\gs}

\newcommand{\averagenafs}[2]{
\averagenaf{#1}[#2]
}

\DeclareMathOperator{\Elle}{L}
\DeclareMathOperator{\quotient}{\Measurable/\ns}

\DeclareMathOperator{\Ellef}{\mathcal{L}} 

\newcommand{\Lspace}[3]{
\Elle^{#1}(#2,#3)
}

\newcommand{\Lspacensa}[2]{
\Elle^{#1}(#2)
}

\DeclareMathOperator{\Loneboundary}{
\Elle^{1}(\bdomain,\hmeas)
}

\DeclareMathOperator{\Loneboundarycv}{
\Elle^{1}(\bdomain,\hsigma,\hmeas)
}

\DeclareMathOperator{\Solution}{V}
\DeclareMathOperator{\Keldych}{K}

\DeclareMathOperator{\PoissonK}{\mathtt{P}}
\DeclareMathOperator{\Poissonk}{\mathtt{p}}

\DeclareMathOperator{\Lebesgue}{L}
\newcommand{\Lebof}[1]{
\Lebesgue{\!}{#1}
}
\newcommand{\Lebofa}[2]{
\Lebesgue{\!}{#1}(#2)}

\newcommand{\Lebofoar}[1]{
{\Lebesgue}^{#1}
}

\newcommand{
\Lebofoarf}[2]{
\Lebofoar{#1}{#2}
}

\newcommand{\Lebofoarfa}[3]{
\Lebofoarf{#1}{#2}(#3)
}

\DeclareMathOperator{\cont}{C}
\newcommand{\continuous}[1]{
\cont(#1)
}

\DeclareMathOperator{\contreg}{\cont_r}

\newcommand{\regular}[1]{
\contreg(#1)
}

\newcommand{\hardy}[2]{
\har^{#1}({#2})
}

\newcommand{\Hardy}[2]{
\Har^{#1}({#2})
}

\DeclareMathOperator{\har}{h}

\newcommand{\harmonic}[1]{
\har(#1)
}

\newtheoremstyle{slthmstyle}  
   {}       
   {}       
   {\slshape}   
   {}        
   {\bfseries}  
   {}   
   {2mm}       
   {}           
\theoremstyle{slthmstyle}
\newtheorem{theorem}{Theorem}[section]

\newtheorem{question}{Question}[section]
\newtheorem{problem}[theorem]{Problem}
\newtheorem{assumption}{Assumption}[section]
\newtheorem{lemma}[theorem]{Lemma}
\newtheorem{corollary}[theorem]{Corollary}
\newtheorem{proposition}[theorem]{Proposition}
\theoremstyle{definition}
\newtheorem{definition}[theorem]{Definition}
\newtheorem{notation}[theorem]{Notation}
\newtheorem{example}[theorem]{Example}
\newtheorem{examples}[theorem]{Examples}

\newtheorem{remark}[theorem]{Remark}
\theoremstyle{definition}

\numberwithin{equation}{section}

\title{On the Differentiation of Integrals \\ 
in Measure Spaces Along Filters: II}

\date{\today}
\author{
{\bf Fausto Di Biase}\\
Dipartimento di Economia,
Universit\`a ``G.\! D'Annunzio'' di Chieti-Pescara,\\
Viale Pindaro, 42, I-65127, Pescara, Italy
email: {\tt fausto.dibiase@unich.it}\\
\\
{\bf Steven G. Krantz}\\
Department of Mathematics, Washington University\\ 
St. Louis MO 63130, USA
email: {\tt sk@wustl.edu}\\
}

\begin{document}

{\tableofcontents}

\newpage

\maketitle

\begin{abstract} 
Let $X$ be a complete measure space of finite measure. 
The Lebesgue transform 
of an integrable function $f$
on $X$ 
encodes the collection of all the mean-values of $f$ on all measurable 
subsets of $X$ of positive measure. In the problem of the 
differentiation of integrals, one seeks to recapture $f$ from 
its Lebesgue transform. In previous work we showed that, 
in all known results, $f$ may be recaputed from its Lebesgue 
transform by means of a limiting process 
associated to an appropriate family of filters defined on the collection 
$\Zm$
of all measurable subsets of $X$ of positive measure. The first result 
of the present work is that the existence of such a limiting process is equivalent 
to the existence of a Von Neumann-Maharam lifting of $X$.

In the second result of this work we provide an independent argument 
that 
shows that 
the recourse to filters is 
a \textit{necessary consequence} of the requirement that 
the process of recapturing $f$ from its mean-values 
is associated to a \textit{natural transformation}, in the sense of category theory.
This result essentially follows  from the
Yoneda lemma.
As far as we know, this is the first instance of a significant interaction between 
category theory and the problem of the differentiation of integrals. 

In the Appendix 
we have proved, in a precise sense, that 
\textit{natural transformations fall within the general
concept of homomorphism}.
As far as we know, this is a novel conclusion: Although 
it is often said that natural transformations are homomorphisms of functors, 
this statement appears to be presented as a mere  analogy, not in 
a precise technical sense.
In order to achieve this result, 
we had to bring to the foreground a notion 
that is implicit in the subject but has remained hidden in the background, i.e., that of \textit{partial magma}.   
\end{abstract}


%


{\footnotesize

}

\section{Introduction}
\label{section:introduction}

The results in this paper are focused on two \textit{existence problems} in the context of a complete measure space $\smm\eqdef\standard$,
where $\hmeas$ us the measure and $\Measurable$ the $\sigma$-algebra of subsets of 
a set $\mm$. 
In order to describe our results, we first have 
to give a brief presentation of these two existence problems.

\subsection{The First Existence Problem}

The first existence problem
 arose in the study of operator algebras and 
is related to the 
\textit{natural projection} 
\begin{equation}
\np: \Measurable\twoheadrightarrow\quotient
\label{eq:urnaturalprojection} 
\end{equation}
where $\ns$ is the $\sigma$-ideal of $\hmeas$-null sets in $\smm$,
$\quotient$,
called 
the \textit{measure algebra of \/$\smm$},
is the quotient Boolean algebra, 
and 
$\np$ is the 
natural 
projection of $\Measurable$ onto its quotient \cite[p.320]{Royden1968bis}. 
Alfr\'ed Haar asked the following question:

\begin{problem}[The First Existence Problem]
Does there exist 
a right-inverse of the natural projection~\eqref{eq:urnaturalprojection}
\begin{equation}
\inp:\quotient\to\Measurable 
\label{e:urlifting}
\end{equation}
which is a
Boolean algebra homomorphism?  
\label{p:first}
\end{problem}

\begin{remark}
The existence of 
a  
right-inverse 
of 
the natural projection~\eqref{eq:urnaturalprojection}
\textit{in the category of sets}
follows at once from the fact that $\np$
is surjective   
 \cite[pp. 6-7]{MacLaneBirkhoffbis}, 
but a 
construction of this sort, 
based on purely set-theoretic grounds, does not 
necessarily 
preserve the Boolean structures involved. 
 What is at stake in this problem 
 is 
the existence of a right-inverse 
of $\np$ that 
\textit{belongs to 
the category of Boolean algebras}. 
\end{remark}

In 1931 John von Neumann gave the first contribution to this problem
showing
that the answer is positive for Lebesgue measure on $\RR$ 
\cite{VonNeumannbis}. In 1958 
Dorothy Maharam gave the second contribution showing 
that 
the answer is positive if $\smm$ is complete and $\sigma$-finite.
See \cite{Fremlin3bis}, 
\cite{IonescuTulceaAC1961}, 
and 
\cite{IonescuTulceaAC1969}
for different proofs, 
more general versions, 
and a wealth of applications. 
The problem is open if 
$\standard$ is not complete, and then 
subtle issues arise \cite{Shelahbis}, \cite{VonNeumannStonebis}.

\subsection{The Second Existence Problem}
\begin{definition}
A real-valued \textit{partial function on}  
 $\mm$ is a real-valued function 
defined on a (possibly empty) subset of $\mm$. 
Instead of the familiar arrow notation, 
the \textit{sparrow tail notation} 
$\mm\parrow\RR$
will be used for added emphasis, and 
we denote by 
$\R^{(\mm)}$  the collection of real-valued 
\textit{partial} functions on $\mm$.
\label{d:partialfunction}
\end{definition}

The following notation is adopted. 

\smallskip

\bspace $\Sm\eqdef\{\sdx\in\Measurable:\hmeas(\sdx)<+\infty\}\subset\Measurable$ 
is the collection 
of \textit{summable} sets in $\smm$.

\blspace

\bspace 
$\Zm\eqdef\{\sdx\in\Sm:0<\hmeas(\sdx)\}\subset\Zm$  is the 
 collection 
of \textit{averageable} sets in $\smm$.

\blspace

\bspace 
$\RR^{\Zm}$
is the 
collection of functions from $\Zm$ to $\RR$.

\blspace

\bspace 
$\Ellef^{0}(\smm)\subset\RR^{(\mm)}$
is the vector space of measurable real-valued 
functions defined $\hmeas$-a.e.\ on $\mm$.

\blspace

\bspace 
$\Ellef^1(\smm)\eqdef
\{
\realbf\in\Ellef^{0}(\smm): 
{|\realbf|} \text{ is integrable}\}$
is the vector space of summable real-valued 
functions defined a.e.\ on $\mm$, while $\Lspacensa{1}\smm$
is the quotient of 
$\Ellef^1(\smm)$ under a.e.\ equality. The natural projection 
\begin{equation}
\Ellef^1(\smm)\twoheadrightarrow\Lspacensa{1}{\smm}
\label{eq:urprojection} 
\end{equation}
maps $\realbf\in\Ellef^1(\smm)$
to the class of functions which are a.e.\ equal to 
$\realbf$: We denote this class 
by $\classbf$ (in bold font), and say that 
$\realbf\in\Ellef^1(\smm)$ is a \textit{representative} of 
$\classbf\in\Lspacensa{1}{\smm}$.

The \textit{Lebesgue map} associated to 
$\smm$ is the 
linear, positive, and injective map 
\begin{equation}
\Lebtb_{\smm}:\Lspacensa{1}{\smm}
\to
\RR^{\Zm}
\label{eq:urnewLebtboldnew} 
\end{equation}
defined as follows:
If 
$\classbf\in\Lspacensa{1}{\smm}$
then 
$\Lebtb_{\smm}(\classbf)\in{\RR}^{\Zm}$
is the 
\textit{Lebesgue transform of} 
$\classbf$
\begin{equation}
\Lebtb_{\smm}(\classbf):\Zm\rightarrow\RR
\label{eq_uraverageonebis}
\end{equation} 
which encodes all the mean-values of $\classbf$:
If $\bsubset\in\Zm$
and
$\realbf\in\Ellef^1(\smm)$ is any representative of 
$\classbf$,
\begin{equation}
\Lebtb_{\smm}(\classbf)(\sdx)
\eqdef\ameanvalue{\realbf}{\bsubset}
\label{eq_urnewaveragebisrealbfnew}
\end{equation} 
Since  
$\Lebtb_{\smm}(\classbf):\Zm\rightarrow\RR$
 is a \textit{bona fide} function defined 
on $\Zm$, while $\classbf$ is an equivalence class of functions, 
the notion of left-inverse of~\eqref{eq:urnewLebtboldnew} is better formulated 
as follows.
\begin{definition}
A \textit{left-inverse of the Lebesgue map} is a function 
$\ell:\RR^{\Zm}\to\RR^{(\mm)}$ such that 
for each $\classbf\in\Lspacensa{1}{\smm}$
\begin{equation}
\ell(\Lebtb_{\smm}(\classbf)) \text{ is a representative of } \classbf
\label{e:mofleftiinverse} 
\end{equation}
\end{definition}

\begin{problem}[Preliminary Version of The Second Existence Problem.]
Does there exist 
a left-inverse
\begin{equation}
\ell:\RR^{\Zm}\to\RR^{(\mm)}
\label{e:tpotdoi}
\end{equation}
of the Lebesgue map which is given by a \textup{limiting process}? 
\label{p:secondpreliminary}
\end{problem}
The precise contours of Problem~\ref{p:secondpreliminary} depends on 
what is intended by a \textit{limiting process}. 
The main goal of this work is to prove that 
the appropriate notion of \textit{limiting process} 
to be employed in the statement of \textup{Problem~\ref{p:secondpreliminary}}
is the one associated to the notion of \textit{filter-kernel}, 
described in \cite{DiBiaseKrantz2023bis}.

\begin{definition}
A  {\it filter-kernel} on $\smm$ is a map 
\begin{equation}
\fofibox:\mm\to\spaceofallfilters{\Zm}
\label{eq:urfamiliesoffiltersbis} 
\end{equation}
where $\spaceofallfilters{\Zm}$ is the collection of all filter on $\Zm$.
\label{d:filterkernel}
\end{definition}
Every filter-kernel on $\smm$ yields a \textit{limiting process}, denoted by 
$\displaystyle{\flim_{\fofibox}}$, as follows. See 
\cite{DiBiaseKrantz2023bis} for background.

\begin{definition}
The \textit{limiting operator}
associated to
a filter-kernel
$\fofibox$ 
on $\smm$
is the 
map
\begin{equation}
\flim_{\fofibox}:
\RR^{\Zm}
\to
\RR^{(\mm)}
\label{eq:urlimitingoperator}
\end{equation}
defined as follows:
If $\lambda\in\RR^{\Zm}$ then 
$\displaystyle{\flim_{\fofibox}\lambda:\mm\parrow\RR}$ 
is the partial function, defined 
by 
$\displaystyle{({\flim_{\fofibox}\lambda)(x)}\eqdef\flim_{\fofibop{\bpoint}}\lambda}$ 
for 
$\displaystyle{\bpoint\in C(\lambda,\fofibox)\eqdef
\{\bpoint\in\mm:
\text{ the limiting value }
\flim_{\fofibop{\bpoint}}\lambda
\text{ exists and is finite }
\}\subset\mm}$.
\end{definition}
We are now ready to present a more precise version of 
Problem~\ref{p:secondpreliminary}.

\begin{problem}[The Second Existence Problem]
Does there exist 
a left-inverse~\eqref{e:tpotdoi} 
of the Lebesgue map which equal to the limiting operator associated to some 
filter-kernel on $\smm$? 
\label{p:second}
\end{problem}

\begin{remark}
The existence of 
a  
left-inverse 
of 
the Lebesgue map
\textit{in the category of sets}
follows at once from the fact that it 
is injective   
 \cite[pp. 6-7]{MacLaneBirkhoffbis}, 
but a 
construction of this sort, 
based on purely set-theoretic grounds, does not  
necessarily fall under the notion of a \textit{limiting process}, unless this notion 
is not given any contours, i.e., it 
is set to encompass 
\textit{any} map~\eqref{e:tpotdoi}. 
\end{remark}

The first contribution to Problem~\ref{p:secondpreliminary} 
was given by Lebesgue, who proved that if $\mm=\RR$ with Lebesgue measure then
the answer is positive, since  
the limiting value of the mean values of 
$\classbf$ 
on intervals shrinking at $\bpoint$
exists for a.e.\ $\bpoint\in\mm$ and yield a representative of $\classbf$.
Lebesgue later extended this result to $\mm=\RR^n$ and proved the existence of 
the limiting values of the mean-values of $\classbf$
over balls which do not necessarily contain $\bpoint$, 
provided the balls shrink to $\bpoint$ in an appropriate “nontangential” manner; 
see \cite{SteinWeiss1971}. 
In 1979, experts were convinced that the mean-values over 
sequences of balls shrinking to $\bpoint$ in a  “tangential” manner  could 
not recapture $\classbf$ \cite{Rudin1979bis}. However, 
in 1984 Alexander Nagel and Elias M.\ Stein showed 
this belief to be wrong \cite{Nagel--Stein1984bis}. 
See \cite{DiBiase1998bis} for an extension of these results to settings 
which lack the presence of a group acting on the metric space. 
In these results, 
the metric structure of the ambient space 
is used to construct a limiting process which, applied to 
the 
Lebesgue transform of $\classbf$, recaptures (a representative of) each $\classbf$.
Results of this kind are 
known as \textit{theorems on the differentiation of integrals}.
See \cite{Brucknerbis}, \cite{DeGuzman1975bis} and \cite{DeGuzman1981bis} for an overview, 
\cite{ArcozziDiBiaseUrbankebis},
\cite{Chirka1973bis}, 
\cite{DePossel1936bis},
\cite{DiBiase1998bis},
\cite{DiBiase2009bis},
\cite{DiBiaseFischer1998bis},
\cite{DiBiaseKrantz2021bis},
\cite{DiBiaseKrantz2023bis},
\cite{DiBiaseStokolosSvenssonWeiss2006bis},
\cite{Fatou1906bis},
\cite{HauptAumannPaucbis},
\cite{Koranyi1969abis},
\cite{Krantz1991bis},
\cite{Krantz2001bis},
\cite{Krantz2007bis},
\cite{Stein1970bis}, \cite{Stein1972bis},
\cite{SteinWeiss1971}, 
for more results of this kind and related applications, 
and \cite{LucicPasqualetto} for recent results.

In 1936 Ren{\'e} de Possel 
observed that if the space is not endowed with a metric
then the familiar device of emplying shrinking metric balls 
is not available, and posed the problem of devising some kind of 
\textit{limiting process} in this generality, 
to be applied to the Lebesgue transform of $\classbf$, that would recapture 
a representative of each $\classbf$ \cite{DePossel1936bis}. 
All known results,
where a left-inverse of the Lebesgue map is exhibited and is given by a 
\textit{limiting process} of some kind, may be reformulated in terms of a
\textit{filter-kernel} on $\smm$; see \cite{DiBiaseKrantz2021bis}. 
In \cite{DiBiaseKrantz2023bis} we argued that the language of filter-kernels is preferable to 
that of Moore-Smith sequences.

The fact that, in all known results,
the limiting process that yields a left-inverse of the Lebesgue map
is associated to a \textit{filter-kernel}, supports the belief that the 
right notion of limiting process, 
to be employed in Problem~\ref{p:secondpreliminary},
is precisely the one associated to a filter-kernel.
However, this belief is merely based on empirical grounds that do not allow us 
to exclude that some different and more useful notion
may exist, still to be discovered. In \cite{DiBiaseKrantz2023bis} we have provided 
some theoretical justification for our belief. 
In this paper, we will provide two theorems which appear to be 
definitive arguments in support of our belief.

\subsection{The First Goal of This Paper}

The statement of the Second Existence Theorem 
is motivated in part 
by the framework introduced 
in \cite{DiBiaseKrantz2023bis}. The first goal of this paper 
is to provide an independent argument (that is, independent of the theoretical considerations presented in \cite{DiBiaseKrantz2023bis})
that 
shows that 
the 
right notion of limiting process, 
to be employed in Problem~\ref{p:secondpreliminary},
is precisely that of a limiting operator associated to a filter-kernel. 
Indeed, we show that the similarities between 
Problem~\ref{p:first} and Problem~\ref{p:second} 
are not superficial, since 
these two existence problem admit the same answer.

\begin{theorem}
Let $\standard$ be a complete measure space of finite measure. 
Then the following conditions 
are equivalent.
\begin{description}

\item[(a)] There exists 
a left-inverse
of the Lebesgue map that is equal to the limiting operator associated to some 
filter-kernel on $\smm$.

\item[(b)] There exists 
a right-inverse of the natural projection~\eqref{eq:urnaturalprojection}
that is a
Boolean algebra homomorphism.

\end{description}

\label{thm:equivalenceofthetwotheorems}
\end{theorem}

\begin{proof}
We may assume that $\Zm\not=\emptyset$. 
The proof of the fact that 
\textbf{(a)} 
implies  
\textbf{(b)}
will be given in Section~\ref{s:equivalenceAB}. 
The proof of the fact that 
\textbf{(b)} 
implies  
\textbf{(a)}
will be given in Section~\ref{s:equivalenceBA}.
\end{proof}

\subsection{The Second Goal of This Paper}

The second goal of this paper is to provide a second, independent argument 
that 
shows that 
the 
right notion of limiting process, 
to be employed in  Problem~\ref{p:secondpreliminary},
is precisely that of a limiting operator associated to a filter-kernel. 
More precisely, 
we show that 
the recourse to filters is 
a \textit{necessary consequence} of the requirement that 
$\ell$ \textit{is associated to a natural transformation}, in the sense of category theory (see below for a precise definition). 
This result is achieved in Theorem~\ref{thm:Yoneda}. We will give two different proofs 
of this result. The first one is a direct, hand-on proof, which only 
requires the notion of natural transformation. The second one is based on the
Yoneda lemma and 
on the notion of adjoint functor. 
The additional interest of Theorem~\ref{thm:Yoneda} lies in the fact 
that, although the requirement that $\ell$ is associated to 
a natural transformation 
has \textit{a priori} nothing to do with limiting processes, it does indeed imply 
that $\ell$ is the limiting operator associated to a filter-kernel. 
As far as we know, this is the first instance of a significant interaction between 
category theory and the problem of the differentiation of integrals.

\subsection{The Third Goal of This Paper}

Since the readership of this work probably has a primary interest and background 
in mathematical analysis, and our second main result 
is cast in the language of category theory, 
we have provided an appendix with a crash course 
in this subject that sheds light on the notion of 
natural transformation. Indeed, 
in the appendix we have strived to show that 
\textit{natural transformations fall within the general
concept of homomorphism}, which is familiar to all mathematician:
This is the third goal of this paper, 
reached in Theorem~\ref{t:proofofequivalence}. 
As far as we know, this is a novel conclusion: Although 
it is often said that natural transformations are homomorphisms of functors, 
this statement appears to be presented as a mere  analogy, not in 
a precise technical sense.
In order to achieve our third goal, 
we had to adopt an iconoclastic take on the subject, and bring to the foreground a notion 
that is implicit in the subject but has remained hidden in the background, i.e., that of \textit{partial magma}.   
We hope the appendix will be of interest to experts as well. 

\section{\textbf{(a)} implies \textbf{(b)} in Theorem~\ref{thm:equivalenceofthetwotheorems}}\label{s:equivalenceAB}

We assume that there exists 
a left-inverse
of the Lebesgue map that is equal to the limiting operator associated to some 
filter-kernel $\fofibox$ on $\smm$.
Our goal in this section is to show that the natural projection~\eqref{eq:urnaturalprojection} has a right-inverse in the category of Boolean algebras. Some preliminary notation 
will be helpful.

\subsection{The Lebesgue Transform on Functions}

The composition of~\eqref{eq:urnewLebtboldnew} with~\eqref{eq:urprojection} 
is the map 
$$
\Lebt_{\smm}:\Ellef^1(\smm)\to{\RR}^{\Zm}
$$
(also called, with a slight abuse of language, 
\textit{the Lebesgue map associated to $\smm$}). 
It follows that, if 
$\realbf\in \Ellef^1(\smm)$
is a \textit{representative} of 
$\classbf\in\Lspacensa{1}{\smm}$,  
then 
$\Lebt_{\smm}(\realbf)\in{\RR}^{\Zm}$
is the function 
\begin{equation}
\Lebt_{\smm}(\realbf):\Zm\rightarrow\RR
\label{eq:Lebesguetransform}
\end{equation}
defined by 
\begin{equation}
\Lebt_{\smm}(\realbf)(\sdx)\eqdef\Lebtb_{\smm}(\classbf)(\sdx) \, .
\label{eq:newmv}
\end{equation}
The map~\eqref{eq:Lebesguetransform}
is also called
\textit{the Lebesgue transform of $\realbf\in\Ellef^1(\smm)$}.

\subsection{Conditional Probability}

If
  $\sdx\in\Measurable$
and
$\sdx'\in\Zm$ 
then  
it is  convenient to borrow 
from 
 probability theory
the vertical bar notation   
to denote the following quantity:
\begin{equation}
\newapairing{\sdx}{\sdx'}\eqdef
\frac{\hmeas(\sdx\cap\sdx')}{\hmeas(\sdx')} \, .
\label{eq:conditional} 
\end{equation}
Recall that the \textit{indicator function of} 
$\sdx\subset\mm$ is the function  
$\indic{\sdx}:\mm\to\RR$ 
defined by 
$\indic{\sdx}(x)=1$ 
if $x\in{}\sdx$ and 
$\indic{\sdx}(x)=0$
if $x\not\in{}\sdx$.
If 
$\indic{\sdx}\in\Ellef^1(\smm)$
and
$\bsubset'\in\Zm$
then
\begin{equation}
\Lebt(\indic{\sdx})(\sdx')=
\newapairing{\sdx}{\sdx'} \, .
\label{eq:uguaglianza} 
\end{equation}

Observe that $\Lebt_{\smm}(\realbf):\Zm\rightarrow\RR$
 is a \textit{bona fide} function defined 
on $\Zm$, while $\classbf$ is an equivalence class of functions, and that 
the values  $\realbf(x)$ of a representative of $\classbf$ may be recovered only up to 
a set of measure zero (called the \textit{exceptional set of} 
$\classbf$). In particular, if $\sdx\in\Measurable$
then $\indic{\sdx}\in\Ellef^1(\smm)$ and 
\begin{equation}
\flim_{\fofibop{\bpoint}}\Lebt_{\smm}(\indic{\sdx})=\indic{\sdx}(x)
\quad\text{$\hmeas$ - a.e. }\bpoint\in\mm \, .
\label{e:almostevery} 
\end{equation}
\subsection{Lower Densities}
Recall that the \textit{symmetric difference} of two sets $\bsubset,\bsubsettwo\subset\mm$ is defined by $\bsubset\Delta\bsubsettwo\eqdef
(\bsubset\setminus\bsubsettwo)\cup(\bsubsettwo\setminus\bsubset)$.
If $\sdx,\sdxb\in\Measurable$, we say that $\sdx$ and 
$\sdxb$ are \textit{almost everywhere equal}, and write 
$\aequiv{\bsubset}{\bsubsettwo}$, 
if $\hmeas(\bsubset\Delta\bsubsettwo)=0$. The binary relation
  ``$\boldsymbol{\saequiv}$'' is an equivalence relation on $\Measurable$.

\begin{definition}
We define the map
$\lambda:\Measurable\to\Measurable$ 
as follows: If $\sdx\in\Measurable$, then 
\begin{equation}
\lambda(\sdx)\eqdef\left\{
\bpoint\in\mm:
\flim_{\fofibop{\bpoint}}\Lebt_{\smm}(\indic{\sdx})=1
\right\} \, .
\label{e:lowerdensity}
\end{equation}
\label{d:lowerdensity}
\end{definition}

\begin{lemma}
The set in~\eqref{e:lowerdensity} is measurable. 
\end{lemma}
\begin{proof}
Observe that~\eqref{e:almostevery}  
implies that $\lambda(\sdx)$
is equal a.e. to $\sdx$, and the conclusion follows 
from the fact that $\smm$ is complete.
\end{proof}
In the following result, we apply Definition~1.12 in \cite{DiBiaseKrantz2023bis}.
\begin{lemma}
If $\alpha,\beta\in\RR^{\Zm}$, $\newafilter\in\soaf{\Zm}$, 
$\alpha(q)\leq\beta(q)\leq1$ for all $q\in\Zm$, and 
$\displaystyle{\flim_{\newafilter}\alpha=1}$,
then
$\displaystyle{\flim_{\newafilter}\beta=1}$.
\label{l:tddc}
\end{lemma}
\begin{proof}
Let $\epsilon>0$. If $q\in\Zm$ and $\alpha(q)\in(1-\epsilon,1+\epsilon)$, 
then $\beta(q)\in(1-\epsilon,1+\epsilon)$. 
Hence 
\begin{equation}
\{q\in\Zm:\alpha(q)\in(1-\epsilon,1+\epsilon)\}
\subset
\{q\in\Zm:\beta(q)\in(1-\epsilon,1+\epsilon)\} \, .
\label{e:larger} 
\end{equation}
Since the smaller set in~\eqref{e:larger} belongs to $\newafilter$, it follows 
that the larger set in~\eqref{e:larger} also belongs to $\newafilter$.
\end{proof}
The following notions will be useful.
\begin{definition}
A map 
$
\varphi:\Measurable\to
\totalpowerset{\mm}
$
may have one or more of the following properties.
\begin{description}
\item[(pms)] ($\boldsymbol{\varphi}$ preserves measurable sets)
$\varphi(\sdx)\in\Measurable$ for all $\sdx\in\Measurable$.

\item[(pas)] ($\boldsymbol{\varphi}$ preserves the ambient space)
$\varphi(\mm)=\mm$.

\item[(pfi)] ($\boldsymbol{\varphi}$ preserves finite intersections)
$\varphi(\bsubset\cap\bsubsettwo)=
\varphi(\bsubset)\cap
\varphi(\bsubsettwo)$ for all $\bsubset,\bsubsettwo\in
\Measurable$.

\item[(aei)] ($\boldsymbol{\varphi}$ is an almost everywhere identity)
$\aequiv{\varphi(\bsubset)}{\bsubset}$
for each $\bsubset\in\Measurable$.

\item[(pes)] ($\boldsymbol{\varphi}$ preserves the empty set)
$\varphi(\emptyset)=\emptyset$.

\item[(spmc)] ($\boldsymbol{\varphi}$ strongly preserves the measure class)
$\aequiv{\bsubset}{\sdxb}$
$\implies$
 $\varphi(\bsubset)=\varphi(\sdxb)$
for all $\bsubset,\sdxb\in\Measurable$.

\item[(cwtc)] ($\boldsymbol{\varphi}$ commutes with the complement)
$\varphi(\complement{\sdx})=\complement{(\varphi(\sdx))}$
for each $\bsubset\in
\Measurable$.

\item[(pfu)] ($\boldsymbol{\varphi}$ preserves finite unions)
$\varphi(\bsubset\cup\bsubsettwo)=
\varphi(\bsubset)\cup
\varphi(\bsubsettwo)$ for all $\bsubset,\bsubsettwo\in
\Measurable$.

\item[(spmcns)] ($\boldsymbol{\varphi}$ strongly preserves the measure class 
of null sets)
$\aequiv{\bsubset}{\emptyset}$
$\implies$
 $\varphi(\bsubset)=\varphi(\emptyset)$
for all $\bsubset\in\Measurable$.

\end{description}
\end{definition}
\vspace*{.15in}

\begin{definition}
If $\varphi:\Measurable\to\totalpowerset{\mm}$ and 
$\lambda:\Measurable\to\totalpowerset{\mm}$, then we say that 
$\lambda$ is \textit{subordinate to} $\varphi$ if 
$$
\varphi(\bsubset)\subset\lambda(\sdx)
\subset\complement{(\varphi(\complement{\sdx}))}
$$
for each $\bsubset\in
\Measurable$.
\end{definition}
\begin{lemma}
Observe that if $\varphi:\Measurable\to\totalpowerset{\mm}$
then 
\begin{equation}
\protect{\textup{\textbf{(spmcns)}}}
\,\,
\&
\,\,
\protect{\textup{\textbf{(cwtc)}}}
\,\,
\&
\,\,
\protect{\textup{\textbf{(pfi)}}}
\,\,
\&
\,\,
\protect{\textup{\textbf{(pes)}}}
\,\,
\&
\protect{\textup{\textbf{(pfu)}}}
\,\,
\Rightarrow
\protect{\textup{\textbf{(spmc)}}};
\label{e_useful}
\end{equation}
and
\begin{equation}
\protect{\textup{\textbf{(cwtc)}}}
\,\,
\&
\protect{\textup{\textbf{(pfi)}}}
\,\,
\Rightarrow
\protect{\textup{\textbf{(pfu)}}};
\label{e_usefulbis}
\end{equation}
\end{lemma}
\begin{proof}
If $\aequiv{\sdx}{\sdxb}$ then 
\begin{equation*}
\begin{split}
\varphi(\sdx)\vartriangle\varphi(\sdxb)
&=
[\varphi(\sdx)\cap\complement{(\varphi(\sdxb))}]
\cup
[\varphi(\sdxb)\cap\complement{(\varphi(\sdx))}]
=
[\varphi(\sdx)\cap(\varphi(\complement{\sdxb}))]
\cup
[\varphi(\sdxb)\cap(\varphi(\complement{\sdx}))]
\\
&=
[\varphi(\sdx\cap\complement{\sdxb})]
\cup
[\varphi(\sdxb\cap\complement{\sdx})]
=
\varphi(\sdx\setminus\sdxb)
\cup
\varphi(\sdxb\setminus\sdx)
=\varphi(\sdx\vartriangle\sdxb)=\varphi(\emptyset)=\emptyset
\end{split}
\end{equation*} 
hence $\varphi(\sdx)=\varphi(\sdxb)$ and~\eqref{e_useful} is proved.
If $\sdx,\sdxb\in\Measurable$ then 
\begin{equation*}
\begin{split}
\complement{(\varphi(\sdx\cup\sdxb))}
&=\varphi(\complement{(\sdx\cup\sdxb)})=
\varphi(
(\complement{\sdx})
\cap
(\complement{\sdxb}))
=
\varphi(\complement{(\sdx)})
\cap
\varphi(\complement{(\sdxb)})\\
&=
\complement{(\varphi(\sdx))}
\cap
\complement{(\varphi(\sdxb))}
=
\complement{(
\varphi(\sdx)
\cup
\varphi(\sdxb)
)}
\end{split}
\end{equation*}
hence~\eqref{e_usefulbis} is proved.
\end{proof}

The main idea of the following result can be found in \cite{VonNeumannbis}.

\begin{lemma}
If $\standard$ is a complete measure space of finite measure 
and $\fofibox$ is a  filter-kernel $\fofibox$ on $\smm$
which differentiates $\Lspacensa{1}{\smm}$, then the 
map $\lambda:\Measurable\to\Measurable$ described in~\eqref{e:lowerdensity} 
preserves the ambient space, preserves finite intersections, is an almost everywhere identity, preserves the empty set, and strongly preserves the measure class.
\label{l:lowerdensity}
\end{lemma}
\begin{proof}
First we prove
that $\lambda$ strongly preserves the measure class.
Observe that if 
$\bsubset,\bsubsettwo\in\Measurable$ and 
$\hmeas(\sd{\bsubset}{\bsubsettwo})=0$
then $\Lebt_{\smm}(\indic{\sdx})=\Lebt_{\smm}(\indic{\sdxb})$, 
hence $\lambda(\sdx)=\lambda(\sdxb)$. 
In order to prove that $\lambda$ is an almost everywhere identity,
observe that 
$\indic{\sdx}\in\Ellef^1(\smm)$ and $\fofibox$ 
differentiates $\Lspacensa{1}{\smm}$.
Hence 
\begin{equation}
\begin{split}
\lambda(\sdx)=\left\{
\bpoint\in\mm:
\flim_{\fofibop{\bpoint}}\Lebt_{\smm}(\indic{\sdx})=1
\right\} &{\saequiv}\left\{
\bpoint\in\mm:
\flim_{\fofibop{\bpoint}}\Lebt_{\smm}(\indic{\sdx})=1=\indic{\sdx}(x)
\right\}\\
 &{\saequiv}\left\{
\bpoint\in\mm:
\indic{\sdx}(x)=1
\right\}\\
&=\sdx \, .
\end{split}
\end{equation}
In order to prove
that 
$\lambda$
preserves the ambient space and the empty set,
observe that 
$\Lebt_{\smm}(\indic{\emptyset})\equiv0$ on $\Zm$ and 
$\Lebt_{\smm}(\indic{\mm})\equiv1$ therein. 
We finally prove that 
$\lambda$
preserves finite intersections. 
Observe that 
$\sdx\cap\sdxb\subset\sdx$ and hence 
$$
\newapairing{\sdx\cap\sdxb}{\sdx'}
\leq
\newapairing{\sdx}{\sdx'}
$$
for each $\sdx'\in\Zm$. In other words, 
$$
\Lebt_{\smm}(\indic{\sdx\cap\sdxb})\leq\Lebt_{\smm}(\indic{\sdx})\leq1
$$ 
on 
$\Zm$. Lemma~\ref{l:tddc} implies that 
$\lambda(\sdx\cap\sdxb)\subset\lambda(\sdx)$.
One shows, along similar lines, that 
$\lambda(\sdx\cap\sdxb)\subset\lambda(\sdxb)$, hence 
$\lambda(\sdx\cap\sdxb)\subset\lambda(\sdx)\cap\lambda(\sdxb)$. Now observe that, for each 
$\sdx'\in\Zm$,
$$
\newapairing{\sdx\cap\sdxb}{\sdx'}
=
\newapairing{\sdx}{\sdx'}+\newapairing{\sdxb}{\sdx'}
-
\newapairing{\sdx\cup\sdxb}{\sdx'}
\geq
\newapairing{\sdx}{\sdx'}+\newapairing{\sdxb}{\sdx'}
-
1 \, .
$$
In other words, 
$$
\Lebt_{\smm}(\indic{\sdx})+\Lebt_{\smm}(\indic{\sdxb})-1
\leq
\Lebt_{\smm}(\indic{\sdx\cap\sdxb})\leq 1
$$ 
on $\Zm$, hence Lemma~\ref{l:tddc} implies that
$\lambda(\sdx)\cap\lambda(\sdxb)\subset\lambda(\sdx\cap\sdxb)$,  
and the proof is concluded.
\end{proof}
The following two notions play a central role in Problem~\ref{e:urlifting}.
\begin{definition}
A function $\lambda:\Measurable\to\Measurable$ 
which preserves the ambient space, the empty set, finite intersections, 
and which is an almost everywhere identity and strongly preserves the measure class, 
is called a \textit{lower density} 
on $\standard$.  
\end{definition}
\begin{definition}
A function $\lambda:\Measurable\to\Measurable$ which
is a lower density on $\smm$ and, moreover, preserves finite unions, 
is called a \textit{lifting} 
of $\smm$.  
\end{definition}
We denote by $\totalpowerset{\Zm}$
the collection of all subsets of $\Zm$ and, for any sets $S_1$ and $S_2$ 
by 
$$
\cathom{\sSet}{S_1}{S_2}
$$
the collection of all functions from 
$S_1$ to $S_2$.

The role played by the notion of lower density is explained by the 
following result.
\begin{lemma}[\cite{Traynorbis}]
If there exists a lower density on a complete measure space of finite measure 
$\smm$ then there exists a  lifting on $\smm$.
\label{l:fldtl}
\end{lemma}
\begin{proof}
Observe that there exists a natural isomorphism 
\begin{equation}
\cathom{\sSet}{\Measurable}{\totalpowerset{\mm}}
\cong
\cathom{\sSet}{\mm}{\totalpowerset{\Measurable}}
\label{e:naturaliso} 
\end{equation}
that can be described thusly: If $\lambda$
belongs to the set in the left-hand side of~\eqref{e:naturaliso}
then we define $\ltr{\lambda}\in\cathom{\sSet}{\mm}{\totalpowerset{\Measurable}}$
as follows: 
$$
\ltr{\lambda}(x)\eqdef\{\sdx\in\totalpowerset{\Measurable}:
x\in\lambda(\sdx)\},\quad(x\in\mm).
$$
Moreover, if $\lambda$ belongs to 
the set in the right-hand side of~\eqref{e:naturaliso}
then we define
$\rtl{\lambda}\in\cathom{\sSet}{\Measurable}{\totalpowerset{\mm}}$ by
$$
\rtl{\lambda}(\sdx)\eqdef\{x\in\mm:\sdx\in\lambda(x)\},\quad(\sdx\in\Measurable).
$$
Observe that $\lambda=\rtl{(\ltr{\lambda})}$ for each $\lambda\in\cathom{\sSet}{\Measurable}{\totalpowerset{\mm}}$ and 
$\lambda=\ltr{(\rtl{\lambda})}$
for each $\lambda\in \cathom{\sSet}{\mm}{\totalpowerset{\Measurable}}$.
We denote by $\mathcal{B}(\mm)$ the collection 
$$
\mathcal{B}(\mm)\eqdef
\{Z:Z\in\totalpowerset{\totalpowerset{\mm}},
Z\not=\emptyset,
\text{ and }
\sdx,\sdxb\in Z \text{ implies that }
\sdx\cap\sdxb\in Z \, .
\}
$$
If we recall 
\cite[Lemma 4.11 and Theorem 4.25]{DiBiaseKrantz2023bis}
we may fix a map 
$$
\natural:\mathcal{B}(\mm)
\to
\usoaf{\mm}
$$
which assigns to each 
$Z\in\mathcal{B}(\mm)$ an ultrafilter ${Z}^{\natural}$
on $\mm$ which contains $Z$. Now let $\lambda$ be a lower density on $\smm$. With a slight abuse of language, $\lambda$ may be seen as an element of 
$\cathom{\sSet}{\Measurable}{\totalpowerset{\mm}}$
which preserves measurable sets. 
The fact that $\lambda$ preserves the ambient space and finite intersections implies that 
\begin{equation}
\ltr{\lambda}(x)\in\mathcal{B}(\mm) 
\quad
\text{ for all } x\in\mm \, .
\end{equation}
Hence we may define 
$\beta:\mm\to\usoaf{\mm}$
by setting
$\beta(x)\eqdef{(\ltr{\lambda}(x))}^{\natural}$
for all $x\in\mm$. Thus for all $x\in\mm$
\begin{equation}
\ltr{\lambda}(x)\subset\beta(x) \, .
\label{e:containment}
\end{equation}
Now define 
$\Lambda\eqdef\rtl{\beta}$, 
hence 
$\Lambda:\Measurable\to\totalpowerset{\mm}$.  
Observe that the fact that $\beta(x)$ is a filter on $\mm$ implies that 
$\Lambda$ preserves finite intersections, that it preserves the ambient space, and 
that
\begin{equation}
\Lambda(\complement{\sdx})\subset\complement{\Lambda(\sdx)}
\quad
\text{ for all }
\sdx\in\Measurable \, .
\label{e:oneside} 
\end{equation}
The fact that $\beta(x)$ is an ultrafilter for all $x\in\mm$, coupled with 
\cite[Lemma 4.26]{DiBiaseKrantz2023bis}, implies that 
$\Lambda(\complement{\sdx})\supset\complement{\Lambda(\sdx)}$
and hence $\Lambda$ commutes with the complement, by~\eqref{e:oneside}. 
Moreover, 
\eqref{e:containment} implies that 
\begin{equation}
\lambda(\sdx)\subset\Lambda(\sdx) 
\quad
\text{ for all }
\sdx\in\Measurable \, .
\label{e:onepart}
\end{equation}
The fact that $\beta(x)$ for all $x\in\mm$ is a filter and~\eqref{e:onepart} imply that 
\begin{equation}
\Lambda(\sdx)\subset\complement{(\lambda(\complement{\sdx}))} 
\quad
\text{ for all }
\sdx\in\Measurable \, .
\label{e:theotherpart}
\end{equation}
Hence~\eqref{e:onepart} and~\eqref{e:theotherpart} imply that 
$\Lambda$ is subordinate to $\lambda$, and this fact implies, coupled with the fact 
that $\lambda$ is an almost everywhere identity, that 
$\Lambda$ is also an almost everywhere identity and $\Lambda$ preserves measurable sets.
Since $\Lambda$ is subordinate to $\lambda$ and 
$\lambda$ preserves the empty set and the ambient space, it follows that 
$\Lambda$ preserves the empty set. Observe that $\Lambda$ preserves 
finite unions by~\eqref{e_usefulbis}. 
In order to conclude the proof, \eqref{e_useful} implies that it suffices to show that 
$\Lambda$ strongly preserves the measure class of null sets, and this follows 
from the fact that it is subordinate to $\lambda$, which strongly preserve the measure class. 
\end{proof}

Lemma~\ref{l:lowerdensity} and Lemma~\ref{l:fldtl} imply that there exists a lifting 
on $\smm$. If $\lambda:\Measurable\to\Measurable$ is a lifting on $\smm$, we define 
$$
\rho_{\lambda}:\quotient\to\Measurable 
$$
by 
$$
\rho_{\lambda}(\np(\sdx))\eqdef \lambda(\sdx).
$$
Then $\rho_{\lambda}$ is a right-inverse of $\np$ 
and is a Boolean algebra homomorphism.

\section{\textbf{(b)} implies \textbf{(a)} in Theorem~\ref{thm:equivalenceofthetwotheorems}}\label{s:equivalenceBA}

Assume that there exists 
a right-inverse $\inp$ \eqref{e:urlifting} of the natural projection~\eqref{eq:urnaturalprojection}
which is a
Boolean algebra homomorphism. Observe that the function 
$$
\lambda:\Measurable\to\Measurable
$$
defined by 
$$
\lambda\eqdef\inp\circ\np
$$
is a lifting on $\smm$. 
\begin{lemma}
If $\lambda:\Measurable\to\Measurable$ is a lifting on $\smm$, then 
for each $\sdx\in\Measurable$ $\lambda(\lambda(\sdx))=\lambda(\sdx)$.
\label{l:idempotent} 
\end{lemma}
\begin{proof}
If $\sdx\in\Measurable$ then 
$\aequiv{\sdx}{\lambda(\sdx)}$ and hence 
$\lambda(\lambda(\sdx))=\lambda(\sdx)$.
\end{proof}

\subsection{The Busemann-Feller Map}

\begin{definition}
The Busemann-Feller map on $\smm$ is the function 
$$
\totalpowerset{\Zm}
\to
\cathom{\sSet}{\mm}{\totalpowerset{\Zm}}
$$ 
which maps $\Lambda\in\totalpowerset{\Zm}$
to $\Zm_{\Lambda}:\mm\to\totalpowerset{\Zm}$, defined as follows: If 
 $\bpoint\in\mm$ then 
$$
\Zm_{\Lambda}(\bpoint)\eqdef\{\sdx:\sdx\in\Lambda, \bpoint\in\sdx\} \, .
$$
\end{definition}

\subsection{Differentiation Bases}
\begin{definition}
A \textit{differentiation basis} on $\smm$ is a collection $\dg\subset\Zm$
such that the following properties hold.
\begin{description}
\item[(1)] For each $\bpoint\in\mm$, if the set 
$\Zm_{\Lambda}(\bpoint)\subset\Zm$ is not empty, 
then it is a directed set under reverse inclusion; see  
\cite[Definition 3.13]{DiBiaseKrantz2023bis} for background on the notion of directed set.

\item[(2)] The set $\mm_{\Lambda}\eqdef\{\bpoint\in\mm:\Zm_{\Lambda}(\bpoint)\not=\emptyset\}$
belongs to $\Zm$ and has full measure in $\smm$.
\end{description}
\end{definition}

\subsection{The Differentiation Basis Associated to a Lifting}
\begin{definition}[\cite{Koelzowbis}]
The collection 
\begin{equation}
\Zm_{\lambda}\eqdef\{\lambda(\sdx):\sdx=\lambda(\sdx)\}
\label{e:dbatal}
\end{equation}
is called the
\textit{differentiation basis} associated to a lifting 
$\lambda$ on $\smm$. 
\end{definition}
The following result is a first example of the significance of the notion of lifting.
\begin{proposition}\textup{\cite[p.\ 57]{Koelzowbis}}
If $\lambda$ is a lifting on $\smm$ then the collection $\Zm_\lambda$ 
in~\eqref{e:dbatal} is a differentiation basis on $\smm$.
\label{p:K1}
\end{proposition}
\begin{proof}
Observe that the set 
$\Zm_{\lambda}$ is not empty. Indeed, since we are assuming that 
$\Zm\not=\emptyset$, if $\sdx\in\Zm$ then $\lambda(\sdx)=\lambda(\lambda(\sdx))$ 
by Lemma~\ref{l:idempotent}
and since $\aequiv{\lambda(\sdx)}{\sdx}$ it follows that 
$\lambda(\sdx)\in\Zm$.
Let $\sdx,\sdxb\in\Zm_{\lambda}(\bpoint)$ and assume that 
$\hmeas(\sdx\cap\sdxb)=0$, i.e.,
$\sdx\cap\sdxb\saequiv\emptyset$.
Then $\lambda(\sdx\cap\sdxb)=\lambda(\emptyset)=\emptyset$. 
Hence $\sdx=(\sdx\cap\sdxb)\cup(\sdx\setminus\sdxb)\saequiv\sdx\setminus\sdxb$ and 
$\lambda(\sdx)=\lambda(\sdx\cap\sdxb)\cup\lambda(\sdx\setminus\sdxb)=
\lambda(\sdx\setminus\sdxb)$ and, symmetrically, 
$\lambda(\sdxb)=\lambda(\sdxb\setminus\sdx)$. 
It follows that 
$$
\bpoint\in\sdx\cap\sdxb=
\lambda(\sdx)\cap\lambda(\sdxb)
=[\lambda(\sdx\setminus\sdxb)]\cap[\lambda(\sdxb\setminus\sdx)]
=\lambda[(\sdx\setminus\sdxb)\cap(\sdxb\setminus\sdx)]
=\lambda[\emptyset]=\emptyset \, ,
$$
a contradiction. Hence $\sdx\cap\sdxb\in\Zm$. Moreover,
$\lambda(\sdx\cap\sdxb)=\lambda(\sdx)\cap\lambda(\sdxb)=\sdx\cap\sdxb$
and $\bpoint\in\sdx\cap\sdxb$, it follows that $\sdx\cap\sdxb\in\Zm_{\lambda}(\bpoint)$, 
and the proof of $\boldsymbol{(1)}$ is complete. Now let $N\eqdef\mm\setminus\mm_{\lambda}$. 
In order to show that $\hmeas(N)=0$ it suffices to show 
that its outer measure is zero. Assume the contrary. Then there exists 
a set $\sdx\in\Zm$ such that
$N\subset\sdx$ and $\hmeas^*(N)=\hmeas(\sdx)$, where $\hmeas^*$ 
is the outer measure induced by $\hmeas$ (see e.g.\ 
\cite[Definition 1.20]{DiBiaseKrantz2023bis}).
Since $\aequiv{\sdx}{\lambda(\sdx)}$, it follows that 
there exists at least one point $\bpoint\in N$ such that 
$\bpoint\in\lambda(\sdx)$, and this is a contradiction. It follows 
that $N$ is a null set, hence $\boldsymbol{(2)}$ is proved. 
\end{proof}
Recall that Proposition~\ref{p:K1} shows that for almost 
every $\bpoint\in\mm$ (indeed, for every $\bpoint\in\mm_{\lambda}$), 
the set $\Zm_\lambda(\bpoint)$ is a directed set under downward inclusion, and 
that $\Zm_\lambda(\bpoint)\subset\Zm$.
For every $\bpoint\in\mm_{\lambda}$, 
we denote by $\imath_\bpoint$ the natural injection of $\Zm_\lambda(\bpoint)$ inside 
$\Zm$. Recall also that the function $\Lebt_{\smm}(\realbf)$
is defined on $\Zm$ and hence the composition 
$\Lebt_{\smm}(\realbf)\circ\imath_x$ is the restriction of 
$\Lebt_{\smm}(\realbf)$ to $\Zm_{\lambda}(\bpoint)$. 
The following commutative diagram illustrates these constructions.
\begin{equation}
\begin{tikzcd}[row sep=large, column sep=huge]
\Zm_\lambda(\bpoint)
\arrow[r,"\imath_x"]
\arrow[dr, "\Lebt_{\smm}(\realbf)\circ\imath_\bpoint" ']
 & 
\Zm
\arrow[d, "\Lebt_{\smm}(\realbf)"]
\\ 
\mbox{}
& 
\RR
\end{tikzcd} 
\label{eq:commutesoneDK}
\end{equation}
In his 1967 \textit{Habilitationsschrift}, Dietrich K\"olzow proved that 
the differentiation basis $\Zm_{\lambda}$ satisfies the strong Vitali-type condition
introduced by Ren\'e de Possel 
in \cite[p.\ 407-408]{DePossel1936bis}, which implies 
\cite[p.216]{HauptAumannPaucbis}
that 
for each $\realbf\in \Ellef^1(\smm)$
\begin{equation}
\gslim_{\Zm_{\lambda}(\bpoint)}\Lebt_{\smm}(\realbf)\circ\imath_{\bpoint}=\realbf(\bpoint) \,\text{ a.\ e.\ } \bpoint\in\mm \, .
\label{e:aeconvergence}
\end{equation}
The meaning of~\eqref{e:aeconvergence} is that for a.e.\ $\bpoint\in\mm$,
the restriction of 
$\Lebt_{\smm}(\realbf)$ on the directed set $\Zm_{\lambda}(\bpoint)$
converges to $\realbf(\bpoint)$ in the sense of Moore-Smith convergence. 
Lemma~3.26 in \cite{DiBiaseKrantz2023bis} implies that 
\begin{equation}
\flim_{\newafilter_{\lambda}(\bpoint)}\Lebt_{\smm}(\realbf)\circ\imath_{\bpoint}=\realbf(\bpoint) \,\text{ a.\ e.\ } \bpoint\in\mm
\label{e:aeconvergencefilter}
\end{equation}
where $\newafilter_{\lambda}(\bpoint)\in\soaf{\Zm_{\lambda}(\bpoint)}$ is the filter 
generated on $\Zm_{\lambda}(\bpoint)$ by the final sets in $\Zm_{\lambda}(\bpoint)$
(see \cite[Section 3.6]{DiBiaseKrantz2023bis} for these notions).
It follows that 
\begin{equation}
\flim_{\fofibox_{\lambda}(\bpoint)}\Lebt_{\smm}(\realbf)
=\realbf(\bpoint) \,\text{ a.\ e.\ } \bpoint\in\mm
\label{e:aeconvergencefilterinZm}
\end{equation}
where 
\begin{equation}
\fofibox_{\lambda}(\bpoint)\eqdef \fdirimFS{\imath_{\bpoint}}{
\newafilter_{\lambda}(\bpoint)}
\label{e:fofibox} 
\end{equation}
is the direct image of $\newafilter_{\lambda}(\bpoint)$ by 
$\imath_{\bpoint}$. 
Observe that $\fofibox_{\lambda}(\bpoint)$ is only defined for 
$\bpoint\in\mm_{\lambda}$, but we may define it 
at other points by setting it equal to the trivial filter on $\Zm$. 
Hence we obtain a filter kernel on $\smm$, also denoted by 
$\fofibox_{\lambda}$, and~\eqref{e:aeconvergencefilterinZm} says that 
the limiting operator 
associated to this filter kernel on $\smm$ 
is a left-inverse of the Lebesgue map.

\section{Notation from Category Theory}\label{s:category}

Since we will use some central notions of category theory, it will be useful to recall some notational devices used in this field.
We adapt to our needs the notation from \cite{MacLane1978bis}, where further background
may be found. Two fundamental collections 
are associated to a category $\sCat$, to wit: The collection of 
\textit{objects}  in $\sCat$ and the collection of \textit{arrows}  in $\sCat$, denoted by 
$$
\sCatO{\sCat}
\quad
\text{ and }
\quad
\sCatA{\sCat}
$$
respectively.
We shall identify the objects $u$ of 
$\sCat$ with the associated identity arrows 
$1_u$ 
as in the “arrows-only” definition of a category 
given in \cite[p.9]{MacLane1978bis}. 
If $u,v\in\sCatO{\sCat}$ then 
\begin{equation}
\cathom{\sCat}{u}{v}\eqdef
\{\alpha\in\sCatA{\sCat}:
\,
\alpha 
\text{ is an arrow in $\sCat$ from $u$ to $v$ }
\}
\label{e:hobjects} 
\end{equation}
Recall that in a concrete category, as for example  
in the category $\sSet$ of sets or in the category $\sRing$
of unital commutative rings, 
arrows do correspond to functions, but in general this is not 
true; see \cite[p.26]{MacLane1978bis}. 
If $\tsnt$ is an object in a concrete category $\sCat$, 
its  underlying set 
is denoted by 
$\ffSet{\tsnt}$. This 
useful notational device for the forgetful functor from 
$\sCat$ to 
$\sSet$ is borrowed from \cite{JamneshanTao2022bis}. 
We may drop the subscript if no confusion is likely to arise. For example, 
if $\fmm$ is a set, then we may occasionally write 
$\cathom{\sSet}{\fmm}{\RR}$
instead of 
$\cathom{\sSet}{\fmm}{\ffSet{\RR}}$.

If $\sCatone,\sCattwo$ are categories, 
and 
$\sFunctor$ and $\sFunctortwo$
are functors from $\sCatone$ to $\sCattwo$, then the notation
\begin{equation}
\tau:\sFunctorone\snt\sFunctortwo
\label{eq:ntMacLane} 
\end{equation}
means that $\tau$ is a natural transformation from $\sFunctorone$ to $\sFunctortwo$.
Following \cite{MacLane1978bis}, 
we omit unnecessary parentheses 
and write, for example, 
$\sFunctorone x$ instead of $\sFunctorone(x)$
whenever possible. 

The fact that 
functors are \textit{homomorphisms of categories} 
is clear at first sight, but the claim that 
natural transformations between functors 
may be conceived as 
\textit{homomorphisms of homomorphisms}, 
may appear to be little more than a mere insubstantial analogy, 
unless one presents this notion within a background where it 
is given a precise meaning. This task is achieved in the 
Appendix, 
which we hope 
will be of interest to 
experts and non-experts alike.

\section{Preliminary Observations}

In order to formulate what it means that the map $\ell$ in~\eqref{e:tpotdoi}
is associated to a natural transformation, observe that we may also consider the Lebesgue map as a function of the following form:
\begin{equation}
\Lebtb_{\smm}:\Lspacensa{1}{\smm}
\to
\cathom{\sSet}{\Zm}{\ffSet{\erl}}
\quad
\label{eq:compact} 
\end{equation}
where $\erl$ is the extended real line $[-\infty,\infty]$, a compact Hausdorff 
topological space; 
see \cite{Bourbaki1995bis}. 
\begin{definition}
A 
left-inverse of the Lebesgue map 
\textit{in the category of sets} is 
a function 
\begin{equation}
\ell:
\cathom{\sSet}{\Zm}{\ffSet{\erl}}
\to
\cathom{\sSet}{\mm}{\ffSet{\erl}}
\label{eq:function} 
\end{equation}
such that, for each  
$\classbf\in\Lspacensa{1}{\smm}$,
\begin{equation}
\ell
\circ
 \Lebtb_{\smm}(\classbf)
\text{ is a representative of }
\classbf
\label{eq:inverse} 
\end{equation}
\label{d:leftinverseinSets}
\end{definition}

\subsection{Two Key Ideas}

In order to understand 
what it means that the map $\ell$ in~\eqref{e:tpotdoi}
is associated to a natural transformation, the first key idea 
is that we should try to 
obtain a map similar to~\eqref{eq:function},   
where  
the extended real line $\erl$ 
is replaced by 
a generic compact Hausdorff space $\tsnt$. 
Hence 
we should seek a map 
\begin{equation}
\tau(\tsnt):
\cathom{\sSet}{\Zm}{\ffSet{\tsnt}}
\to
\cathom{\sSet}{\mm}{\ffSet{\tsnt}}
\label{eq:morphism} 
\end{equation}
for each $\tsnt$,  
seen as a 
parameter that varies in the category $\sTopCH$ of compact Hausdorff topological spaces. 
The value $\tau(\erl)$  of such a map at $\erl$ 
should then give us the desired left-inverse, hence 
\begin{equation}
\ell=\tau(\erl)
\label{e:associatedto} 
\end{equation}
However, 
it would be hopeless to expect that  
every possible assignment
\begin{equation}
\tsnt\mapsto\tau(\tsnt)
\label{eq:association} 
\end{equation}
(where $\tau(\tsnt)$ is as in~\eqref{eq:morphism}
and $\tsnt$ varies among the objects of $\sTopCH$)
will be of some interest. Some further restriction on the map~\eqref{eq:association} 
must be imposed.

The second key idea is 
 that 
the map $\tau(\tsnt)$   \textit{should change in a consistent way} as 
$\tsnt$ varies in $\sTopCH$. The precise meaning of this consistency condition is that
the following condition holds:
\begin{equation}
\varphi\,\circ\,\tau(\tsnttwo)(L)=\tau(\tsnt)(\varphi\circ L)
\label{eq:natural}
\end{equation}
for each pair of compact Hausdorff topological spaces 
$\tsnt, \tsnttwo$,  each continuous function $\varphi:\tsnttwo\to\tsnt$, 
and each function $L:\Zm\to\tsnttwo$. 
Observe that a map $\tau$ as in~\eqref{eq:association} with property~\eqref{eq:natural} 
is a natural transformation 
between the two 
hom-functors $\sTopCH\to\sSet$ that appear in~\eqref{eq:morphism}, i.e., 
as in~\eqref{eq:ntMacLane}, 
\begin{equation}
\tau:
\cathom{\sSet}{\Zm}{
\ffSet{(\,\cdot\,)}}
\snt
\cathom{\sSet}{\mm}{\ffSet{(\,\cdot\,)}} \, .
\label{eq:ntbhf} 
\end{equation}
The meaning of~\eqref{eq:natural}, as MacLane puts it, is that 
$\tau$
is \textit{defined in the same way for all objects in 
the given category}; see \cite[p. 79]{MacLane1978bis}.
Indeed, \eqref{eq:natural} may be seen as a \textit{commutativity condition} 
and hence it indicates that a homomorphism is at work in the background. 
In the Appendix 
we will substantiate these comments and present 
from scratch 
the key ideas of category theory in the context of \textit{partial magmas}, 
that are the appropriate variant of the notion of 
\textit{magmas}, due to Bourbaki.  We hope our presentation will be of interest 
to experts and non experts alike. These two key ideas boil down to the following definition.
\begin{definition}
We say that s left-inverse $\ell$
of the Lebesgue map in the category of sets 
is associated to a natural transformation
if there exists a natural transformation~\eqref{eq:ntbhf} 
such that~\eqref{e:associatedto} holds. 
\end{definition}

\subsection{A Representation Theorem for Certain Natural Transformations}

At this point, our task is to
identify the structure of the natural transformations~\eqref{eq:ntbhf}.
The collection of all filters on $\Zm$ plays a role in the description of this structure,
and we will denote it by $\soaf{\Zm}$, as in~\cite{DiBiaseKrantz2023bis}.
Recall that an {\it ultrafilter} is a filter which is maximal 
in the natural order structure of the powerset in which it is contained. 
See for example \cite[Section 4]{DiBiaseKrantz2023bis}.
\begin{theorem}
Each natural transformation 
$\tau$ in~\eqref{eq:ntbhf}
is induced by some map 
\begin{equation}
\fofibox:\mm\to\spaceofallfilters{\Zm}
\label{eq:familiesoffiltersbis} 
\end{equation}
such that 
\begin{equation}
\text{for each $\bpoint\in\mm$,}\,
\fofibop{\bpoint}
\text{ is an ultrafilter on $\Zm$}
\label{eq:ultrafilters} 
\end{equation}
as follows: Given a map~\eqref{eq:familiesoffiltersbis} such that~\eqref{eq:ultrafilters} holds, 
for each object $\tsnt$ of $\sTopCH$,
we define 
\begin{equation}
\tau_{\fofibox}(\tsnt):
\hset{\Zm}{\ffSet{\tsnt}}
\to
\hset{\mm}{\ffSet{\tsnt}} \, .
\label{eq:nat_transf} 
\end{equation}
Here, for 
$L\in \hset{\Zm}{\tsnt_{\sSet}}$
and
$\bpoint\in\mm$,  
\begin{equation}
\tau_{\fofibox}(\tsnt)(L)(\bpoint)\eqdef 
\flim_{\fofibop{\bpoint}}L\in\tsnt
\label{eq:limitingprocess}
\end{equation}
for $\displaystyle{\flim_{\fofibop{\bpoint}}L\in\tsnt}$
the limiting value of $L:\Zm\to\tsnt$ along the filter $\fofibop{\bpoint}$
on $\Zm$.
Then $\tau_{\fofibox}$ is a natural transformation~\eqref{eq:ntbhf}, 
and each natural transformation~~\eqref{eq:ntbhf} has this form.
\label{thm:Yoneda}
\end{theorem}
\begin{proof}
The proof will be given in Section~\ref{s:Yoneda}.  
\end{proof}
Comprehensive information about 
the notion of limiting value of a function along a filter, due to Henri Cartan, 
which appears in~\eqref{eq:limitingprocess}, 
may be found in \cite{DiBiaseKrantz2023bis}. 

The interest of Theorem~\ref{thm:Yoneda} is twofold. 
On the one hand, it shows that 
maps as in~\eqref{eq:familiesoffiltersbis}
appear in a natural way out of the requirement 
that~\eqref{eq:morphism}  should yield a natural transformation~\eqref{eq:ntbhf}. 
On the other hand, it shows that limiting processes arising from filters, 
as in~\eqref{eq:limitingprocess}, also appear in a natural way because of the same requirement. 
This conclusion confirms on theoretical grounds a fact that can be verified on empirical grounds: 
All known formulations and variants of the Lebesgue Differentiation Theorem are 
based on a limiting process which is expressed in terms of a map 
as in~\eqref{eq:familiesoffiltersbis}, as shown in~\cite{DiBiaseKrantz2021bis}.

If we drop~\eqref{eq:ultrafilters} and the compactness  requirement, we are led to the category $\spSet$
of \textit{partial functions} between sets. Recall that, if 
$f\in\hpset{\mm}{\RR_{\sSet}}$, then $f$ is a real-valued function 
defined on a (possibly empty) subset of $\mm$, and 
instead of the familiar arrow notation 
$f:\mm\to\RR$, the \textit{sparrow tail notation} will be used for added emphasis: 
$f:\mm\parrow\RR$.

\begin{definition}
If $\fofibox$ is a filter-kernel
on $\smm$, then the map 
\begin{equation}
\flim_{\fofibox}:\hset{\Zm}{\ffSet{\RR}}
\to
\hpset{\mm}{\ffSet{\RR}} 
\label{eq:limitingoperator}
\end{equation}
is defined as follows: 
If $L\in\hset{\Zm}{\ffSet{\RR}}$, then 
$\displaystyle{\flim_{\fofibox}L\in\hpset{\mm}{\ffSet{\RR}}}$, i.e., 
$$
\flim_{\fofibox}L:\mm\parrow\RR
$$ 
is a partial function defined on the set of points $\bpoint\in\mm$ for which the limiting value
$\displaystyle{\flim_{\fofibop{\bpoint}}L}$
exists, and 
$\displaystyle{({\flim_{\fofibox}L)(x)}\eqdef\flim_{\fofibop{\bpoint}}L}$ at such points.
The map~\eqref{eq:limitingoperator} is the \textit{limiting operator 
of the filter-kernel} $\fofibox$.  
\end{definition}
\begin{example}
If $\fofibop{\bpoint}$ is the trivial filter on $\Zm$
and 
$L\in\hset{\Zm}{\ffSet{\RR}}$ is not constant, then 
$\displaystyle{\flim_{\fofibox}}L$ is defined nowhere. 
If $L$ is a constant, then  
$\displaystyle{\flim_{\fofibox}}L$ is defined everywhere and is everywhere equal to that constant.
\end{example}

\section{Proof of Theorem~\ref{thm:Yoneda}}\label{s:Yoneda}

We find it instructive to offer two renditions of the proof. 
The first one, presented in Section~\ref{s:firstproof},  
has been obtained having in mind the Yoneda lemma, but \textit{before} 
we could actually deduce it from this general result. 
Hence it is direct and detailed but longer and a bit unpredicable, unless one has in mind 
the general philosophy of the Yoneda lemma. The second one, presented in Section~\ref{s:secondproof}, 
is formally deduced from the Yoneda lemma, with the aid of a pair of adjoint functors. 
If we may draw an analogy with computer programming, 
we could say that 
then the first one
may be considered a 
\textit{low level proof}, akin to 
 \textit{machine language instructions}, while the second one is 
 similar to 
\textit{higher-level languages}. 

Since the first proof of Theorem~\ref{thm:Yoneda} was obtained having in mind 
the general philosophy of the Yoneda lemma, the following result, pointed out by Terence Tao in~\cite{Tao2023}, may be helpful to readers with no preliminary knowledge of this result. 
\begin{proposition}
Let us assume that, for each commutative ring $R$, a function
$$
\tau(R):R\to R
$$
is given with the following property: For each 
ring homomorphism $f:S\to U$ between any given commutative rings
$S$ and $U$, 
the following diagram is commutative:
\begin{equation}
\begin{tikzcd}[row sep=large, column sep=huge]
S
\arrow[r,"f" ']
\arrow[d, "\tau(R)" ']
 & 
U
\arrow[d, "\tau(S)" ']
\\ 
S
\arrow[r, "f" ']
& 
U
\end{tikzcd} 
\label{eq:urcommutesone}
\end{equation}
Then there exists a polynomial $p$ with integer coefficients in one indeterminate, 
say 
$p=\sum_{k=0}^{N}p_k x^k$, 
such that, for every commutative ring $R$
and each $a\in R$, 
\begin{equation}
\tau(T)(a)=\sum_{k=0}^{N}p_k a^k  \, .
\label{eq:YonedaTao} 
\end{equation}
\label{p:TaoYonedapolynomials}
\end{proposition}
\begin{proof}
(See \cite{Tao2023}.) Let $\ZZ[\ttx])$ 
be the ring of polynomials with integer coefficients 
in the indeterminate $\ttx$, and observe that
$\tau(\ZZ[\ttx]):\ZZ[\ttx]\to\ZZ[\ttx]$
is a function
from $\ZZ[\ttx]$ to $\ZZ[\ttx]$, and hence its value at $\ttx$
is a polynomial, say
$$
p\eqdef\tau(\ZZ[\ttx])(\ttx)=\sum_{k=0}^n p_k\ttx^k \, .
$$ 
If $R$ is a commutative unital ring 
and $a\in R$, let 
$f_{a}:\ZZ[\ttx]\to R$ be the unique unital ring homomorphism such that 
$f_a(\ttx)=a$, and consider the commutative diagram
obtained from~\eqref{eq:urcommutesone} 
by replacing $S$ with $\ZZ[\ttx]$,
$U$ with $R$, 
and $f$ with $f_a$.
\begin{equation}
\begin{tikzcd}[row sep=large, column sep=large]
\ZZ[\ttx]
\arrow[r,"f_a" ']
\arrow[d, "\tau(\ZZ{[}\ttx{]})" ']
 & 
R
\arrow[d, "\tau(R)" ]
&\ttx\arrow[r,mapsto]\arrow[d,mapsto] &f_a(\ttx)=a
\arrow[d,mapsto]
\\ 
\ZZ[\ttx]
\arrow[r, "f_a" ']
& 
R
&p\arrow[r, mapsto] &f_a(p)=\tau(R)(a)
\end{tikzcd} 
\label{eq:commutesoneTao}
\end{equation}
Then 
\begin{equation}
\begin{split}
\tau(R)(a)=
\tau(R)(f_a(\ttx))
=
\tau(R)\,\circ\, f_a\,(\ttx)
&=
f_a
\,\circ\,
\tau(\ZZ[\ttx])
\,
(\ttx)
\\
=
f_a\left(\tau(\ZZ[\ttx])(\ttx)\right)
=
f_a(p)
&=
f_a\left(\sum_{k=0}^n p_k\ttx^k\right)
=
\sum_{k=0}^n p_k a^k
\end{split}
\end{equation}
and this means that 
$\tau(R)(a)=\sum_{k=0}^n p_k a^k$, which is precisely~\eqref{eq:YonedaTao}.
\end{proof}
Observe that 
Proposition~\ref{p:TaoYonedapolynomials} 
says that each 
natural transformation from the forgetful functor 
\begin{equation}
\ffSet{(\,\cdot\,)}:\sRing\to\sSet 
\label{e:forgetfulfunctor}
\end{equation}
to itself is given by 
an element of the ring $\ZZ[x]$ of polynomial with integer coefficients.
This result may also be obtained by an application of the Yoneda lemma, 
once 
we observe that, for each unital commutative ring $R$, 
$$
\cathom{\sRing}{\ZZ[\ttx]}{R}\cong \ffSet{R}
$$
and hence the forgetful functor~\eqref{e:forgetfulfunctor}
may be identified with the hom-set functor determined by  
$\ZZ[\ttx]$.

The proof of Theorem~\ref{thm:Yoneda}, 
presented below in Section~\ref{s:firstproof}, 
has been obtained having in mind the direct proof of 
Proposition~\ref{p:TaoYonedapolynomials} 
given in~\cite{Tao2023}, reproduced above, 
which conveys the general philosophy behind the Yoneda lemma.

\subsection{Preliminary Facts}\label{s:prelfacts}

Recall that if $\fmm$ is a nonempty set then 
a filter on $\fmm$ is a nonempty collection of nonempty subsets 
of $\fmm$ which is closed under finite intersection and has the property 
that if it contains a subset it also contains every set which contains that subset. 
It is useful to keep in mind, as a motivating example, 
the collection of all neighborhoods of a given point in a topological space. 
The notion of filter may be thought of as a substitute of this collection 
when no topology is given.

The 
collection $\soaf{\fmm}$ of all filters on $\fmm$ is endowed with a natural topology, 
described in 
\cite{DiBiaseKrantz2023bis}, where further background material on filters may be found.
Under this topology, $\soaf{\fmm}$ is compact but not Hausdorff, unless 
$\fmm$ is a singleton. 
The collection $\usoaf{\fmm}$
of all ultrafilters on $\fmm$, endowed with the induced topology, is compact and Hausdorff, 
but not closed in $\soaf{\fmm}$ unless $\fmm$ is a singleton. Indeed, if 
$\newafilter\in\usoaf{\fmm}$ and we denote by 
$\overline{\{\newafilter\}}$
the closure of the “singleton” set $\{\newafilter\}$ in the natural topology of $\soaf{\newafilter}$, 
then 
$$
\overline{\{\newafilter\}}
=
\{
\newbfilter\in\soaf{\fmm}:
\newbfilter\subset\newafilter
\}
$$
as shown in \cite[Lemma 16.33]{DiBiaseKrantz2023bis}.

In~\eqref{eq:limitingprocess} we use the
notation described 
in~\cite{DiBiaseKrantz2023bis}. 
Hence the statement that 
\begin{equation}
\flim_{\fofibop{\bpoint}}L=y\in\tsnt
\label{eq:lone} 
\end{equation}
is equivalent to the statement that 
\begin{equation}
\fdirimFS{L}{\fofibop{\bpoint}}\supset\nsdi{\tsnt}{y}
\label{eq:ltwo} 
\end{equation}
where 
$\fdirimFS{L}{\fofibop{\bpoint}}$ is the image of $\fofibop{\bpoint}$ 
along $L$ and $\nsdi{\tsnt}{y}$ is the filter of neighborhoods of $y$ 
in $\tsnt$. In~\cite{DiBiaseKrantz2023bis} we showed that~\eqref{eq:ltwo} is equivalent 
to the following statement:
 \begin{equation}
\text{$\fdirimFS{L}{\fofibop{\bpoint}}$ converges to $\nsdi{\tsnt}{y}$ in the natural topology of 
$\soaf{\tsnt}$}
\label{eq:lthree} 
\end{equation}
which means that the constant Moore-Smith sequence identically 
equal to $\fdirimFS{L}{\fofibop{\bpoint}}$ converges to $\nsdi{\tsnt}{y}$ in the natural topology of 
$\soaf{\tsnt}$. It will be useful to put these ideas together in a coherent fashion. 
We need the following preliminary observation. 
\begin{lemma}
For each object  $\tsnt$ of\/ $\sTopCH$ 
and each $\newafilter\in\usoaf{\tsnt}$, 
there exists one and only one element
$y$ of $\tsnt$ such that 
the constant Moore-Smith sequence identically 
equal to $Z$ converges to $\nsdi{\tsnt}{y}$ in the natural topology of 
$\soaf{\tsnt}$, i.e., 
\begin{equation}
\newafilter\supset\nsdi{\tsnt}{y} \, .
\label{eq:limitingvalue} 
\end{equation}
\label{l:limitingvalue}
\end{lemma}
\begin{proof}
It suffices to apply Corollary 4.35, Proposition~10.6, 
and Lemma 10.8 in \cite{DiBiaseKrantz2023bis}. 
 \end{proof}
Lemma~\ref{l:limitingvalue} enables us to introduce the following definition.
\begin{definition}
For each object  $\tsnt$ of\/ $\sTopCH$ the map 
\begin{equation}
\gslim_{\soaf{\tsnt_{\sSet}}}:\usoaf{\tsnt_{\sSet}}\to\tsnt 
\label{eq:limitingvalueofafilter}
\end{equation}
is defined by setting 
$$
\gslim_{\soaf{\tsnt_{\sSet}}}\newafilter=y
$$ 
where $y\in\tsnt$ is the unique element of $\tsnt$ such 
that~\eqref{eq:limitingvalue} holds. 
\label{d:limitingvalue}
\end{definition}
Recall that the natural topology of $\soaf{\tsnt_{\sSet}}$, described in~\cite{DiBiaseKrantz2023bis}, 
does not depend on the topology of $\tsnt$.  
\begin{lemma}
For each object  $\tsnt$ of\/ $\sTopCH$ the map~\eqref{eq:limitingvalueofafilter}
is continuous, where $\usoaf{\tsnt_{\sSet}}$ has the topology induced by 
the natural topology of $\soaf{\tsnt_{\sSet}}$, and, 
for each object $\tsnttwo$ in $\sTopCH$ and each  
each $\varphi\in\htop{\tsnttwo}{\tsnt}$ 
the following diagram commutes:
\begin{equation}
\begin{tikzcd}[row sep=large, column sep=huge]
\usoaf{\tsnttwo_{\sSet}}
\arrow[r,"\displaystyle{\gslim_{\soaf{\tsnttwo_{\sSet}}}}" ']
\arrow[d, "\fdirimF{\varphi}" ']
 & 
\tsnttwo\arrow[d, "\varphi" ]
\\ 
\usoaf{\tsnt_{\sSet}}
\arrow[r, "\displaystyle{\gslim_{\soaf{\tsnt_{\sSet}}}}" ']
& 
\tsnt
\end{tikzcd} 
\label{eq:commutesone}
\end{equation}
\label{l:naturalone}
\end{lemma}
\begin{proof}
First, we demonstrate the continuity of~\eqref{eq:limitingvalueofafilter}.
Let $\newafilter\in\usoaf{\tsnt_{\sSet}}$ and assume that 
$\displaystyle{\lim_{\soaf{\tsnt_{\sSet}}}\newafilter=y}$.  
Hence 
$\newafilter\supset\nsdi{\tsnt}{y}$.
Let $\fm$  be an open set which contains $y$. 
Then $\fm\in\newafilter$.
Recall from \cite{DiBiaseKrantz2023bis} that 
$$
\soaftcags{\tsnt}{\fm}
\eqdef
\setofsuchthat{\newbfilter}{
\newbfilter\in\usoaf{\fmm},
\fm\in\newbfilter
}
\subset\soaf{\tsnt} 
$$
is a neighborhood of $\newafilter$ in the natural topology of $\usoaf{\tsnt_{\sSet}}$.
Let $\newbfilter\in\soaftcags{\tsnt}{\fm}$. We claim that 
$$
\gslim_{\soaf{\tsnt_{\sSet}}}\newbfilter\in\fm \, .
$$ 
Indeed, let $\displaystyle{y'=\gslim_{\soaf{\tsnt_{\sSet}}}\newbfilter}$ and assume that 
$y'\not\in\fm$. Since $\tsnt$ is Hausdorff, there exists $\fmb\in\nsdi{\tsnt}{y'}$
such that $\fmb\cap\fm=\emptyset$. The fact that 
$\displaystyle{y'=\gslim_{\soaf{\tsnt_{\sSet}}}\newbfilter}$
means that 
$\newbfilter\supset\nsdi{\tsnt}{y'}$, hence 
$\fmb\in\newbfilter$. The fact that 
$\newbfilter\in\soaftcags{\tsnt}{\fm}$
means that 
$\fm\in\newbfilter$, but $\fm$ and $\fmb$ cannot be disjoint, since they 
belong to the same filter, hence 
$\displaystyle{\gslim_{\soaf{\tsnt_{\sSet}}}\newbfilter\in\fm}$. 
Hence we have shown that for each open set $\fm$ which contains 
$\displaystyle{y=\gslim_{\soaf{\tsnt_{\sSet}}}\newafilter}$
there exists an open set in $\usoaf{\tsnt_{\sSet}}$ 
that contains $\newafilter$
whose image through 
$\displaystyle{\gslim_{\soaf{\tsnt_{\sSet}}}}$ belongs to $\fm$, and the proof of continuity is completed.

Secondly, we show that the diagram~\eqref{eq:commutesone} commutes.
Let $\newafilter\in\usoaf{\tsnttwo_{\sSet}}$, $y\in\tsnttwo$, and 
assume that $\displaystyle{y=\gslim_{\soaf{\tsnttwo_{\sSet}}}}\newafilter$. 
Hence $\newafilter\supset\nsdi{\tsnttwo}{y}$. Now let 
$\tsnt$ be a compact Hausdorff topological space and assume that 
$\varphi\in
\htop{\tsnttwo}{\tsnt}$. Observe that the continuity of $\varphi$ means that 
$$
\fdirimFS{\varphi}{\nsdi{\tsnttwo}{y}}\supset\nsdi{\tsnt}{\varphi(y)}
$$
Corollary 11.36 in~\cite{DiBiaseKrantz2023bis} implies that 
$\fdirimFS{\varphi}{\newafilter}\supset\fdirimFS{\varphi}{\nsdi{\tsnttwo}{y}}$, 
hence 
$\fdirimFS{\varphi}{\newafilter}\supset\nsdi{\tsnt}{\varphi(y)}$, and this means that 
$$
\gslim_{\soaf{\tsnt_{\sSet}}}\fdirimFS{\varphi}{\newafilter}=\varphi(y) \, .
$$ 
Since $\displaystyle{y=\gslim_{\soaf{\tsnttwo_{\sSet}}}\newafilter}$, we have proved that the diagram~\eqref{eq:commutesone} commutes.
\end{proof}

\begin{remark}
The notion of natural transformation will be defined in the Appendix, 
since it would be inconvenient to interrupt this presentation to define it here. 
The Appendix also contains a scratch course on category theory 
where the subject is recast in terms of the notion of partial magma. 
We believe the presentation given in the Appendix will make it more palatable to 
readers who have a primary training in analysis. 
\end{remark}
\begin{remark}
Recall that the peculiar notation $\displaystyle{\lim_{\soaf{\tsnt_{\sSet}}}}$ 
has been introduced  in
Definition~\ref{d:limitingvalue}. This notation is peculiar because 
one usually takes the limit \textit{of} something. The point here is that 
the object we are taking the limit \textit{of} is the variable of a natural transformation.
\end{remark}
\begin{corollary}
The assignment 
$$
\displaystyle{\tsnt\mapsto\lim_{\soaf{\tsnt_{\sSet}}}}
$$ 
is a natural transformation 
from the functor $\sTopCH\to\sTopCH$ which maps $\tsnt$ to 
$\usoaf{\tsnt_{\sSet}}$ and $\varphi:\tsnttwo\to\tsnt$ to 
 $\fdirimF{\varphi}:\usoaf{\tsnttwo_{\sSet}}\to\usoaf{\tsnt_{\sSet}}$, to the identity functor from $\sTopCH$ to itself. 
\end{corollary}
\begin{proof}
The statement follows immediately from Lemma~\ref{l:naturalone}.  
\end{proof}

\subsection{The Direct Proof of Theorem~\ref{thm:Yoneda}}\label{s:firstproof}

The first presentation will be divided into two parts:

\begin{description}
\item[(a)] Each filter-kernel~\eqref{eq:familiesoffiltersbis} such that~\eqref{eq:ultrafilters} holds induces a natural transformation~\eqref{eq:ntbhf}. 
\item[(b)] Each natural transformation~\eqref{eq:ntbhf}
is induced by a 
filter-kernel~\eqref{eq:familiesoffiltersbis} such that~\eqref{eq:ultrafilters} holds. 
\end{description}

\subsubsection{Proof of \textbf{(a)}}

If $\fofibox$ is a 
filter-kernel~\eqref{eq:familiesoffiltersbis} such that~\eqref{eq:ultrafilters} holds, 
we define, for each object $\tsnt$ of $\sTopCH$,
the map 
$\tau_{\fofibox}(\tsnt)$ in~\eqref{eq:nat_transf} where, for 
$L\in \hset{\Zm}{\tsnt_{\sSet}}$, the value of the function 
$$
\tau_{\fofibox}(\tsnt)(L):\mm\to\tsnt
$$
at $\bpoint\in\mm$
is given by~\eqref{eq:limitingprocess}. 
Recall  that 
$\tau_{\fofibox}(\tsnt)(L)(\bpoint)\eqdef 
\flim_{\fofibop{\bpoint}}L\in\tsnt$
means, in the notation of Section~\ref{s:prelfacts},  that
\begin{equation}
\tau_{\fofibox}(\tsnt)(L)(\bpoint)=\lim_{\soaf{\tsnt_{\sSet}}}\fdirimFS{L}{\fofibop{\bpoint}} \, .
\label{eq:ipotesi}
\end{equation}
In order to show that the assignment 
$$
\tsnt\mapsto\tau_{\fofibox}(\tsnt)
$$
so defined is a natural transformation~\eqref{eq:ntbhf}, recall that, for each fixed set 
$B$, we define the covariant hom-functor $\sTopCH\to\sSet$, 
denoted by
$$
\hset{B}{{(\cdot)}_{\sSet}} \, .
$$ 
This functor maps the object 
$\tsnt$ of $\sTopCH$ to the set 
$\hset{B}{{\tsnt}_{\sSet}}$  
and the morphism $\varphi\in\htop{\tsnt}{\tsnttwo}$ to 
the function
$$
\hset{B}{{\varphi}}:\hset{B}{{\tsnt}_{\sSet}}\to\hset{B}{{\tsnttwo}_{\sSet}}
$$
defined by 
$$
\hset{B}{{\varphi}}(L)\eqdef\varphi\circ L\in\hset{B}{{\tsnttwo}_{\sSet}}
$$
for each $L\in\hset{B}{{\tsnt}_{\sSet}}$. See \cite{MacLane1978bis} for background. 
Now let $\tsnt,\tsnttwo$ be objects in $\sTopCH$ and let 
$\varphi\in\htop{\tsnttwo}{\tsnt}$. We claim that the following diagram commutes:
\begin{equation}
\begin{tikzcd}[row sep=large, column sep=huge]
\hset{\Zm}{\tsnttwo_{\sSet}}
\arrow[r,"\tau_{\fofibox}(\tsnttwo)" ']
\arrow[d, "\hset{\Zm}{\varphi}" ']
 & 
 \hset{\mm}{\tsnttwo_{\sSet}}\arrow[d, "\hset{\mm}{\varphi}" ]
\\ 
\hset{\Zm}{\tsnt_{\sSet}}
\arrow[r, "\tau_{\fofibox}(\tsnt)" ']
& 
\hset{\mm}{\tsnt_{\sSet}}
\end{tikzcd} 
\label{eq:commutes}
\end{equation}
Indeed, if $L\in\hset{\Zm}{\tsnttwo_{\sSet}}$ then the value of 
$$
\hset{\mm}{\varphi}\,\circ\,\tau_{\fofibox}(\tsnttwo)(L)
$$
at $x\in\mm$ is given by 
\begin{equation}
\varphi(\lim_{\soaf{\tsnttwo_{\sSet}}}\fdirimFS{L}{\fofibop{\bpoint}}) \, .
\label{eq:one}
\end{equation}
On the other hand, 
the value of 
$$
\tau_{\fofibox}(\tsnttwo)\,\circ
\hset{\Zm}{\varphi}\,(L)
$$
at $x\in\mm$ is given by 
\begin{equation}
\lim_{\soaf{\tsnt_{\sSet}}}\fdirimFS{(\varphi\circ\,L)}{\fofibop{\bpoint}} \, .
\label{eq:two}
\end{equation}
Now recall from \cite[Lemma 11.36]{DiBiaseKrantz2023bis} that 
$$
\fdirimFS{(\varphi\,\circ\,L)}{\fofibop{\bpoint}}
=
\fdirimFS{\varphi}{\fdirimFS{L}{\fofibop{\bpoint}}}
$$
Hence the conclusion follows from Lemma~\ref{l:naturalone}.
\hfill{\qed}

\subsubsection{Proof of \textbf{(b)}}

Now assume that $\tau$ is a natural transformation~\eqref{eq:ntbhf}. 
In particular, for each $\varphi\in\htop{\tsnttwo}{\tsnt}$, 
where $\tsnt$ and $\tsnttwo$ are objects in $\sTopCH$, the following diagram commutes:
\begin{equation}
\begin{tikzcd}[row sep=large, column sep=huge]
\hset{\Zm}{\tsnttwo_{\sSet}}
\arrow[r,"\tau(\tsnttwo)" ']
\arrow[d, "\hset{\Zm}{\varphi}" ']
 & 
 \hset{\mm}{\tsnttwo_{\sSet}}\arrow[d, "\hset{\mm}{\varphi}" ]
\\ 
\hset{\Zm}{\tsnt_{\sSet}}
\arrow[r, "\tau(\tsnt)" ]
& 
\hset{\mm}{\tsnt_{\sSet}}
\end{tikzcd} 
\label{eq:commutestwo}
\end{equation}
Recall from~\cite{DiBiaseKrantz2023bis} that 
$\usoaf{\Zm}$ is a compact Hausdorff space, and that 
the natural injection
\begin{equation}
\delta:\Zm\to\usoaf{\Zm}
\label{eq:naturalin} 
\end{equation}
is defined by 
$\delta(\sdx)\eqdef\{U\subset\Zm:\sdx\in{}U\}$. Hence  
$\delta(\sdx)$ is the principal ultrafilter on $\Zm$ 
generated by $\sdx$; see \cite{DiBiaseKrantz2023bis}. 
Now observe that 
$\tau(\usoaf{\Zm})(\delta)$
belongs to $\hset{\mm}{\usoaf{\Zm}}$, i.e.,
$\tau(\usoaf{\Zm})(\delta):\mm\to\usoaf{\Zm}$
and define
$$
\fofibox\eqdef\tau(\usoaf{\Zm})(\delta) \, .
$$
Hence $\fofibox$ is a filter-kernel on $\mm$. We claim that 
\begin{equation}
\tau=\tau_{\fofibox} \, .
\label{eq:conclusione}
\end{equation}
Now let 
$\tsnt$ be a compact Hausdorff space 
and let $L\in\hset{\Zm}{\tsnt_{\sSet}}$. Then 
the function 
$$
\fdirimF{L}:\usoaf{\Zm}\to\usoaf{\tsnt_{\sSet}}
$$ 
is continuous, where 
$\usoaf{\Zm}$
and 
$\usoaf{\tsnt_{\sSet}}$ are endowed with the natural topology;
see \cite[Lemma 16.31]{DiBiaseKrantz2023bis}. 
Indeed, this result does not depend on the fact that $\tsnt$ is an object 
of $\sTopCH$. 
Recall also that Lemma~\ref{l:naturalone} implies that 
$\displaystyle{\lim_{\soaf{\tsnt_{\sSet}}}:\usoaf{\tsnt_{\sSet}}\to\tsnt}$ is continuous. 
It follows that the function 
\begin{equation}
\varphi\eqdef\lim_{\soaf{\tsnt_{\sSet}}}\circ\,\fdirimF{L}:\usoaf{\Zm}\to\tsnt
\label{eq:definitionofvarphi} 
\end{equation}
is continuous. 
We claim that 
\begin{equation}
\varphi\,\circ\,\delta=L  \, ,
\label{eq:usefullemmaone} 
\end{equation}
where $\delta$ is the function defined in~\eqref{eq:naturalin}. 
Indeed, if 
$\sdx\in\Zm$ then 
$$
\varphi\,\circ\,\delta(\sdx)
=
\lim_{\soaf{\tsnt_{\sSet}}}\circ\,\fdirimF{L}
\,\circ\,\delta(\sdx)
=
\lim_{\soaf{\tsnt_{\sSet}}}
\fdirimFS{L}{\delta(\sdx)}
=
L(\sdx)
$$
since
$$
\fdirimFS{L}{\delta(\sdx)}
\supset
\nsdi{\tsnt}{L(\sdx)} \, .
$$
In point of fact, if $\fm\in\nsdi{\tsnt}{L(\sdx)}$ then, in particular, ${L(\sdx)}\in\fm$, 
hence $\sdx\in\invim{L}{\fm}$, thus 
$\invim{L}{\fm}\in\delta(\sdx)$, hence 
$\fm\in\fdirimFS{L}{\delta(\sdx)}$. 
Therefore the proof of~\eqref{eq:usefullemmaone} is complete.\hfill{\qed}
\smallskip \\

We have proved that \
$$
\varphi\in\htop{\usoaf{\tsnt_{\sSet}}}{\tsnt}
$$
and 
$$
\hset{\Zm}{\varphi}(\delta)=L \, .
$$
Hence we may specialize~\eqref{eq:commutestwo} to the following commutative 
diagram:
\begin{equation}
\begin{tikzcd}[row sep=large, column sep=huge]
\hset{\Zm}{\usoaf{\tsnt_{\sSet}}}
\arrow[r,"\tau(\usoaf{\tsnt_{\sSet}})" ']
\arrow[d, "\hset{\Zm}{\varphi}" ']
 & 
 \hset{\mm}{\usoaf{\tsnt_{\sSet}}}\arrow[d, "\hset{\mm}{\varphi}" ]
\\ 
\hset{\Zm}{\tsnt_{\sSet}}
\arrow[r, "\tau(\tsnt)" ]
& 
\hset{\mm}{\tsnt_{\sSet}}
\end{tikzcd} 
\label{eq:commutestwobis}
\end{equation}
and do “arrow chasing” starting from $\delta\in\hset{\Zm}{\usoaf{\tsnt_{\sSet}}}$. 
We obtain the following diagram: 
\begin{equation}
\begin{tikzcd}[arrows=mapsto, row sep=large, column sep=huge]
\delta
\arrow[r,"\tau(\usoaf{\tsnt_{\sSet}})" ']
\arrow[d, "\hset{\Zm}{\varphi}" ']
 & 
 \fofibox
 \arrow[d, "\hset{\mm}{\varphi}" ]
\\ 
L
\arrow[r, "\tau(\tsnt)" ]
& 
\tau(\tsnt)(L)=\varphi\circ\fofibox
\end{tikzcd} 
\label{eq:commutestwoter}
\end{equation}
We now apply~\eqref{eq:definitionofvarphi} and conclude that  
\begin{equation}
\tau(\tsnt)(L)(\bpoint)=
\varphi\circ\fofibox(\bpoint)
=
\varphi(\fofibox(\bpoint))
= 
\lim_{\soaf{\tsnt_{\sSet}}}
(\fdirimF{L}
(\fofibox(\bpoint))) \, .
\label{eq:arrived}
\end{equation}
If we compare~\eqref{eq:arrived} to~\eqref{eq:ipotesi} we conclude 
that~\eqref{eq:conclusione} holds.
\hfill{\qed}

\subsection{Theorem~\ref{thm:Yoneda} Follows From the Yoneda Lemma}\label{s:secondproof}

The direct proof of Theorem~\ref{thm:Yoneda}
presented in Section~\ref{s:firstproof}
implicitly contains 
a detailed proof of the following statements, hence we limit 
ourselves to give an outline of the results as seen from the 
perspective of the Yoneda lemma. 
The notion of adjoint functor may be found in \cite{MacLane1978bis}. 

The first ingredient is the fact that the forgetful functor 
$\sTopCH\to\sSet$
is right-adjoint to the ultrafilter functor 
$$
\sSet\to\sTopCH
$$
that maps a set $\fmm$ to the compact Hausdorff topological space 
$\usoaf{\fmm}$, and maps a function $f:\fmm\to\newds$ 
between sets $\fmm, \newds$
to the continuous function 
$\fdirimF{f}:\usoaf{\fmm}\to\usoaf{\newds}$. 
A proof of this statement is implicit in the direct proof given in~Section~\ref{s:firstproof}. 
The meaning of this statement is that, for each set $\fmm$
there exist a natural isomorphism
$$
\hset{\fmm}{\tsnt_{\sSet}}\cong
\cathom{\sTopCH}{\usoaf{\fmm}}{\tsnt} \, .
$$
In particular, for the same reason, 
for each set $\mm$
there exists a natural isomorphism 
\begin{equation}
\hset{\mm}{\tsnt_{\sSet}}
\cong
\cathom{\sTopCH}{\usoaf{\mm}}{\tsnt} \, .
\label{e:adjointfunctors} 
\end{equation}
This result implies that a natural transformation 
\begin{equation}
\tau:\hset{\fmm}{{(\,\cdot\,)}_{\sSet}}\snt
\hset{\mm}{
{(\,\cdot\,)}_{\sSet}
}
\label{eq:ntbhfnew} 
\end{equation}
amounts to a natural transformation
\begin{equation}
\cathom{\sTopCH}{\usoaf{\fmm}}{{(\,\cdot\,)}}
\snt
\cathom{\sTopCH}{\usoaf{\mm}}{{(\,\cdot\,)}} \, .
\label{eq:ntbhfnewnew} 
\end{equation}
The Yoneda lemma implies that a natural 
transformation as in~\eqref{eq:ntbhfnewnew} amounts to an element of 
\begin{equation}
\cathom{\sTopCH}{\usoaf{\mm}}{\usoaf{\fmm}}
\label{e:eureka} 
\end{equation}
and~\eqref{e:adjointfunctors} says that the set in~\eqref{e:eureka} 
amounts 
(via a natural isomorphism)
to 
\begin{equation}
\cathom{\sSet}{{\mm}}{{\usoaf{\fmm}}_{\sSet}}  \, .
\label{e:anewproof} 
\end{equation}
Hence there is a natural correspondence between natural 
 transformations~\eqref{eq:ntbhfnew}
 and elements of~\eqref{e:anewproof}, and the elements of~\eqref{e:anewproof} are precisely 
 the kernel-filters.

\section{Appendix on Natural Transformations}
\label{a:nt}

The first goal of this section is to offer speedy access to 
the notion of a \textit{natural transformation} 
to 
readers who are not familiar 
with this language. A functor is readily seen to be 
a homomorphism between categories, 
and sometimes in the literature a natural transformation between functors 
is described as a 
\textit{homomorphism of homomorphisms}. The second goal of this section 
is to show that this description is more than a mere insubstantial analogy.
In order to achieve these goals, we first will have to  
show how the notions 
of \textit{category}, 
 \textit{functor},
 and
  \textit{natural transformation}
  fit 
 within the more general concept of \textit{partial magma}. 
We hope this section 
will be of interest to 
experts and non-experts alike. 

\subsection{Magmas and Partial Magmas}

The notion of \textit{partial magma}, introduced  in Definition~\ref{d:pm},
may be placed within a hierarchy of
axiomatic theories, expressed in first-order languages, 
which allow for 
partial operations; see 
 \cite{AdamekEtAlbis}, \cite[Ch. 2]{Graetzer}, 
 \cite{Poythress}. 

\subsubsection{Partial Algebraic Structures}
In the context of axiomatic theories, 
expressed in first-order languages, the most general notion, 
which thus plays a special role, is that of 
\textit{partial metamagma}, 
defined by \textit{the axiom of the existence of 
a partial binary operation}. This axiom may be  
expressed by the existence of 
a \textit{functional} ternary relation, as in \cite[p. 217]{McLean}.
Richer theories are obtained by adding other axioms, such as 
\textit{associativity}, \textit{regularity}, and so on
(see below). 
However, we will not dwell on the syntactical description of these axiomatic systems, 
but instead concentrate on their models, i.e., on the \textit{interpretation within set theory 
of the partial metamagma axiom} given above, as in \cite[p.10]{MacLane1978bis} (see Definition~\ref{d:pm}).
Hence our emphasis is placed on these models,  
and 
we will ignore the distinction, made in mathematical logic, between 
axiomatic systems (described in appropriate languages)  
and their models: In so doing, we 
cling to the regular but unsound practice 
adopted by the working mathematicians
who are usually unaware of the distinction 
\cite{Bourbaki1949}, 
\cite{DiBiase2009bis}. 
For background on these ideas, see 
\cite{Johnstone},
\cite{Manin} and 
\cite{VanDalen}.
The diagram~\eqref{eq:moregeneral} describes the 
hierarchy between the structures defined by specializing the notion of 
partial magma, introduced below.  
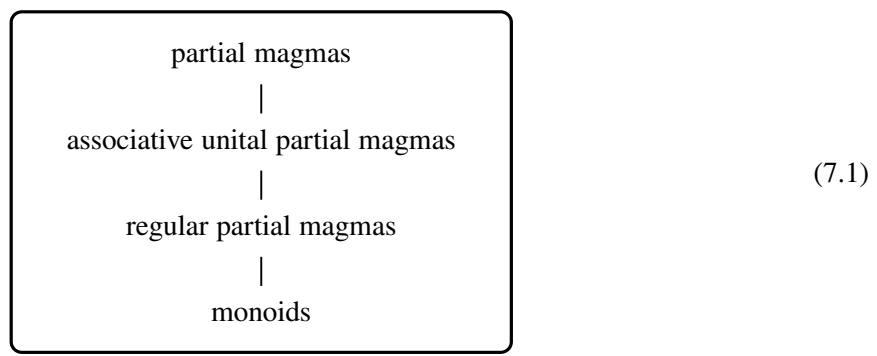
\begin{equation}
\begin{tikzcd}
[mymatr, arrow style=tikz, 
>={Triangle[length=2mm]}, 
cells={nodes={inner sep=2mm}}, 
row sep=0.4cm,column sep=1cm,
every arrow/.append style={dash}]
\text{partial magmas}
\arrow[d]
\\ 
\text{associative unital partial magmas}
\arrow[d]
\\
\text{regular partial magmas}
\arrow[d]
\\
\text{monoids}
\end{tikzcd} 
\label{eq:moregeneral}
\end{equation}
\subsubsection{Caveat}
In dealing with categories, 
one is faced with the so-called “size issues”, where a distinction 
is made between “small sets” and “large sets”; 
see \cite{Blass} for an overview 
and~\cite[Ch. 1]{MacLane1978bis} for an accessible account. 
For this reason, 
some authors, instead of  \textit{set}
use the term \textit{aggregate}, 
as a trigger word to  alert the reader  that
one may be dealing with sets that are “too large”; 
see e.g. \cite[p.237]{EilenbergMacLane}.
These problems  overlap with the study of the 
foundations of mathematics and its set-theoretic 
interpretation, and the rightful desire to stay at large from the 
well-known antinomies of «ordinary intuitive \textit{Mengenlehre}» 
\cite[p.246]{EilenbergMacLane}.
We will not delve deeply into these issues, 
and, in order to avoid these technical problems, 
we will restrict part of our treatment 
to the realm of 
the so-called “small categories”, leaving the task of 
handling more general categories to another occasion; 
see \cite[Ch. 1]{MacLane1978bis}. 
This restriction will not undermine the main goal of this section: 
To show that 
natural transformations are, in a precise sense, 
homomorphisms of homomorphisms.

\subsection{The Notion of Magma}\label{s:caveat}

It will be convenient 
to first recall the notion of \textit{magma}. 
This term was apparently introduced by Jean-Pierre Serre, in his 1964 Harvard lectures 
\cite[Ch. 4]{Serre1964}, published in English in 1965. 
The term was later popularized by Bourbaki in the 1970 edition of his 
\textit{Alg{\`e}bre} 
\cite{Bourbakibis}.

\begin{definition}
A \textit{magma} $B$ is a set 
endowed with a binary operation $B\times B\to B$.
\label{d:magma}
\end{definition}

\subsection{Partial Magmas}\label{s:partialmagmas}

\textit{Partial magmas} belong to the realm of 
\textit{partial algebras}, where operations are 
\textit{partial} functions, i.e., functions defined on a possibly proper subset 
of their “domain”;
see \cite[Ch. 2]{Graetzer}; see also 
\cite{Poythress}.

\begin{definition}
A \textit{partial magma} $B$ is a set of objects 
endowed with a partial 
binary operation 
\begin{equation}
B\times B\parrow B
\label{eq:pbo}
\end{equation}
usually written multiplicatively, 
where the sparrow tail notation in~\eqref{eq:pbo} means, 
as in Definition~\eqref{d:partialfunction}, that the map is defined on 
a (possibly proper) subset of $B\times B$. 
\label{d:pm}
\end{definition}

The notation presented in the following section will be used systematically. 

\subsection{Notation}\label{s:compositionismultiplication}
If $\Setuno$ and $\Setdue$ are sets then the elements of the product
$\Setuno\times\Setdue$ are denoted in bold font: If 
$\mbf{x}\in\Setuno\times\Setdue$ then 
$$
\mbf{x}=(x(1),x(2)), \quad 
x(1)\in\Setuno,
\,
x(2)\in\Setdue.
$$
The following terminology will be very useful, as we will see momentarily.
\begin{definition}
If $\Setuno$ is a set, a
\textit{twin in} $\Setuno$
is an element of $\csquare{\Setuno}$.
If $\mbf{x}$ is a twin in $\Setuno$ then 
$\mbf{x}=(x(1),x(2))$ with $x(k)\in\Setuno$ for $k=1,2$.
\end{definition}

Arrows in \textit{concrete} categories are functions, 
and their “composition” is actual composition of functions, 
but in general this interpretation fails. 
Since “composition” of arrows in a category behaves  like 
multiplication in a partial magma (see Lemma~\ref{l:afapm}), 
it will be written as multiplication, as in~\cite[p. 237]{EilenbergMacLane}. 
Indeed, unless confusion is likely to arise, 
the binary operation in a partial magma will be mostly 
 written multiplicatively, and
the partial binary operation in different partial magmas 
will also be denoted by the same (multiplication) sign. 
The only exception is given by the so-called \textit{horizontal multiplication} in
$\csquare{\Setuno}$, denoted by $\mbf{x}\hmult\mbf{y}$
(Definition~\ref{d:hm}), defined for \textit{any} nonempty set ${\Setg}$.
If $B$ is a magma, then 
$\csquare{B}$ is also endowed with the so-called 
\textit{vertical multiplication}, 
denoted by $\vm{\mbf{x}}{\mbf{y}}$ (see Section~\ref{s:vm}).

\subsubsection{The Square Functor}\label{s:squarefunctor}
\begin{definition}
The \textit{square} of a set $\Setuno$ is the set 
$\squaref{\Setuno}\eqdef\csquare{\Setuno}$. The \textit{square} of a function
 $f:\Setuno\to\Setdue$ is the function 
$\squaref{f}:\csquare{\Setuno}\to\csquare{\Setdue}$
defined, for $\mbf{x}\in\csquare{\Setuno}$, by 
$$
(\squaref{f})(\mbf{x)}\eqdef(f(x(1)),f(x(2))).
$$
\label{d:squarefunctor}
\end{definition}
\begin{lemma}
The assignments in Definition~\ref{d:squarefunctor} determine a functor 
\begin{equation}
\sSquareF:\sSet\to\sSet.
\label{e:twinfunctor} 
\end{equation} 
\label{l:squarefunctor}
\end{lemma}
\begin{proof}
It suffices to show that if 
$g\cdot f$
is
the composition 
of the functions $f:\Setuno\to\Setdue$ and
$g:\Setdue\to\Settre$, then the composition 
$(\squaref{g})\cdot (\squaref{f})$
is defined and is equal to 
$\squaref{(g\cdot f)}$, and $\squaref{1_{\Setuno}}=1_{\squaref{\Setdue}}$.
Indeed, if $\mbf{x}\in\csquare{\Setuno}$ then 
$\squaref{(g\cdot f)}(\mbf{x})
=((g\cdot f)(x(1)),(g\cdot f)(x(2)))=(g(f(x(1))),g(f(x(2))))
=\squaref{g}(\squaref{f}
(\mbf{x}))$, and hence 
$\squaref{(g\cdot f)}(\mbf{x})=(\squaref{g}\cdot\squaref{f})(\mbf{x})$. Moreover, 
$\squaref{1_{\Setuno}}(\mbf{x})
=(1_{\Setuno}(x(1)),1_{\Setuno}(x(2)))=(x(1),x(2))=\mbf{x}$.
\end{proof}
\begin{definition}
The functor in~\eqref{e:twinfunctor} is called the \textit{square functor}.
\end{definition}

\subsection{Examples of Partial Magmas}
Natural numbers form a magma under addition. A more interesting example is the following one. 
\begin{example}
The natural numbers $\NN\eqdef\{0,1,2,\ldots\}$ are a {partial} magma under {subtraction}. 
\end{example}

The following example of a partial magma  seemingly comes out of thin air.

\subsubsection{Horizontal Multiplication: The Twin Partial Magma}\label{s:hm}

\begin{definition}
If $\Setg$ is a nonempty set then
the \textit{horizontal multiplication} 
$$
\hmult:\csquare{(\csquare{\Setg})} \parrow \csquare{\Setg}
$$ 
is defined as follows: If 
$\mbf{x},\mbf{y}\in\csquare{\Setg}$ then 
$\hm{\mbf{x}}{\mbf{y}}$ is defined only if $y(2)=x(1)$, and then 
$$
\hm{\mbf{x}}{\mbf{y}}\eqdef(y(1),x(2)).
$$
\label{d:hm}
\end{definition}
It is convenient to read the terms in the product $\hm{\mbf{x}}{\mbf{y}}$  \textit{from right to left}, 
as if we were dealing with functions. Then we obtain precisely 
$y_1,y_2,x_1,x_2$, and horizontal multiplication consists in erasing the middle terms, if they are equal, and leaving the product undefined otherwise, as in the following diagram.
\begin{equation}
\begin{tikzcd}[mymatr, arrow style=tikz, 
>={Triangle[length=2mm]}, 
cells={nodes={inner sep=2mm}}, 
row sep=0.1cm,column sep=1cm] 
\circ\arrow[d,"y(1)" ', phantom]		
&
\diamond\arrow[d,"y(2)" ', phantom]	
&
\diamond\arrow[d,"x(1)" ', phantom]	
&
\bullet
\arrow[d,"x(2)" ', phantom]
\\	
\mbox{}		
&
\mbox{}	
&
\mbox{}	
&
\mbox{}	
\mbox{}	
\end{tikzcd}
\parrow
\begin{tikzcd}[mymatr,  
arrow style=tikz, 
>={Triangle[length=1mm]}, 
cells={nodes={inner sep=2mm}}, 
row sep=0.1cm,column sep=1cm]
\circ\arrow[d,"y(1)" ', phantom]		
&
\bullet
\arrow[d,"x(2)" , phantom]
\\	
\mbox{}	
&
\mbox{}	
\end{tikzcd}
\label{e:hm}
\end{equation}
\begin{lemma}
If $\Setg$ is a nonempty set then $\csquare{\Setg}$ is a partial magma 
under horizontal multiplication.  
\label{l:horizontal}
\end{lemma}
\begin{proof}
Horizontal multiplication is defined for each pair of elements of 
$\csquare{\Setg}$ only if $\Setg$ is a singleton. 
\end{proof}
\begin{notation}
The \textit{set} $\csquare{\Setg}$ endowed with 
horizontal multiplication $\hmult$
is called the 
\textit{twin partial magma generated by} $\Setg$
and is denoted by $\twinf{\Setg}$. Hence 
$$
\twinf{\Setg}
\eqdef
(\csquare{\Setg},\hmult).
$$
Hence \textit{the partial magma}
$\twinf{\Setg}$ is supported on \textit{the set} $\csquare{\Setg}$.
\end{notation}

\subsubsection{Partial Magmas of Matrices}
Matrices provide a natural source of partial magmas.
\begin{example}
The set $\scMatrices_{\ZZ}$ of all matrices 
with integer coefficients
 is a partial magma. 
Indeed, matrix multiplication $\matrixone\cdot\matrixtwo$
is defined only if 
the number of columns of the matrix
$\matrixone$
is equal to 
the number of rows of $\matrixtwo$; 
see \cite[p.11]{MacLane1978bis}.
\label{eg:matrices}
\end{example}
\begin{examples}
Any subset of $\scMatrices_{\ZZ}$ 
is a partial magma, with the same partial operation. 
The $n\times{}k$ matrix that has $1$ on the diagonal and 
zero as off-diagonal entries 
is denoted by 
$\matrixone^{n}_{k}$.  
Hence
\begin{equation}
\aidentity{2}{1}\eqdef \dueuno,\,
\aidentity{3}{2}\eqdef\tredue,\,
\aidentity{3}{1}\eqdef \treuno,\,
\aidentity{4}{1}\eqdef\quattrouno,\,
\aidentity{3}{4}\eqdef\trequattro
\label{e:matrices}
\end{equation}
Moreover, 
we denote by 
$\ione$, $\ii$, $\iii$, and $\iv$
the identity matrices of order one, two, three and four, respectively. 
Here is a list of partial magmas which, because of their simplicity, are 
useful tools in the theory; cf. \cite{MacLane1978bis}.
\begin{itemize}
\item 
The set $\{\ione\}$ is a partial magma. Observe that $\ione\cdot \ione=\ione$.
\item 
The set $\{\ione,\ii\}$ is a partial magma. Observe that 
$\ione\cdot \ione=\ione$ and $\ii\cdot \ii=\ii$,
but $\ione\cdot \ii$ and $\ii\cdot \ione$ are not defined.
\item 
The set
$\{\ione,\ii,\aidentitydueuno\}$ is a partial magma.
Observe that 
$\aidentitydueuno\cdot \ione=\aidentitydueuno=
\ii\cdot \aidentitydueuno$, 
but for example $\ione\cdot \aidentitydueuno$ is not defined.

\item The set 
$\{\ione,\ii,\iii,\aidentitydueuno,\aidentitytredue,
\aidentitytreuno\}$ is a partial magma. Observe that 
$\aidentitytredue\cdot\aidentitydueuno=
\aidentitytreuno$

\item The set 
$\{\ione,\ii,\iii,\iv,
\aidentitydueuno,
\aidentitytredue,
\aidentitytreuno,
\aidentityquattrouno,
\aidentitytrequattro
\}$ is a partial magma. Observe that 
$\aidentitytreuno=
\aidentitytredue\cdot \aidentitydueuno
=\aidentitytrequattro\cdot \aidentityquattrouno$
\end{itemize}
\label{eg:examples} 
\end{examples}

\subsubsection{Small Categories Yield Partial Magmas}

We now observe that the arrows 
of a small category 
yield a partial magma; see \cite[p.237-238]{EilenbergMacLane}.
\begin{lemma}
If \/ $\sCat$ is a small category then  the set $\sCatA{\sCat}$ of arrows in $\sCat$ forms a partial magma, denoted by $\sCatApm{\sCat}$, where multiplication is composition of arrows. 
\label{l:afapm}
\end{lemma}
\begin{proof}
Composition of arrows in 
 $\sCat$
is defined
for each pair of arrows if and only if $\sCat$ has only one object, since 
it is defined 
 only if 
their domain and codomain are related by the familiar ``chain rule.''
\end{proof}

\subsection{The Category of Partial Magmas} 
The following notion of 
\textit{homomorphism of partial magmas} 
 is best suited to our needs; see \cite[Ch.\ 2]{Graetzer}. 
On this basis, we will define the category of partial magmas. 
Let $Q, B$ be partial magmas.
\subsubsection{Homomorphisms of Partial Magmas} 
\begin{definition}
A \textit{homomorphism of partial magmas} 
is a function $\sFunctor:Q\to B$ 
such that
\begin{description}
\item[(H)]
For every $x,y$ in $Q$, 
if $x\cdot y$ is defined then 
$(\sFunctor x)\cdot(\sFunctor y)$ is defined and is equal to 
$\sFunctor(x\cdot y)$.
\end{description}
\end{definition}

\subsubsection{Partial Magmas Form a Category} 
\begin{lemma}
Partial magmas are the objects of a  category, 
denoted by $\sPM$, whose arrows are the partial magma homomorphisms.
\label{l:PM} 
\end{lemma}
\begin{proof}
Let 
$B_1, B_2, B_3$ be partial magmas, and 
 $f:B_1\to B_2$ and 
$g:B_2\to B_3$ partial magma homomorphisms. 
If $x,y\in B_1$ 
and if $x\cdot y$ is defined in $B_1$,
then $f(x)\cdot f(y)$
is defined in $B_2$, since $f$ is a partial magma homomorphism, 
and $f(x\cdot y)=f(x)\cdot f(y)$. Since 
$g$ is a partial magma homomorphism, it follows that 
$g(f(x))\cdot g(f(y))$ is defined in $B_3$, and 
$g(f(x)\cdot f(y))=g(f(x))\cdot g(f(y))$. Hence 
$$
\text{If $x\cdot y$ is defined in $B_1$ then 
$g(f(x))\cdot g(f(y))$ is defined in $B_3$
and 
$g(f(x\cdot y))=g(f(x))\cdot g(f(y))$.
}
$$
and this means that the composition $g\cdot f:B_1\to B_3$ is a partial magma homomorphism. Associativity holds for composition of functions, 
and if $B$ is a partial magma, then $1_B:B\to B$ is a partial magma 
homomorphism which acts as an identity with respect to composition. 
\end{proof}
The set of homomorphisms of partial magmas from $Q$ to $B$
is denoted by 
$$
\cathom{\sPM}{Q}{B} \, .
$$
\subsubsection{Examples of Partial Magmas Homomorphisms} 
In the following result we will use the 
notation from Section~\ref{s:compositionismultiplication} 
and the notions introduced in Section~\ref{s:hm}. 
\begin{lemma}
The assignments 
$\Setg\mapsto \twinf{\Setg}$
and
$f\mapsto \squaref{f}$
described in Lemma~\ref{l:squarefunctor} and in Section~\ref{s:hm} 
determine a functor 
\begin{equation}
\sTwin:\sSet\to\sPM \, .
\label{e:thetwinfunctor} 
\end{equation}  
\label{l:functorSetTwinThisOne}
\end{lemma}
\begin{proof}
It suffices to show that 
$\squaref{f}:\twinf{\Setuno}\to\twinf{\Setdue}$
is a homomorphism of partial magmas, by Lemma~\ref{l:squarefunctor}. 
If $\mbf{x},\mbf{y}\in\twinf{\Setuno}$ 
and
 $\mbf{x}\hmult{\mbf{y}}$ is defined, then $y(2)=x(1)$, hence  
$f(y(2))=f(x(1))$. It follows that 
$\hm{(\squaref{f})\mbf{x}}{(\squaref{f})\mbf{y}}$
is defined, i.e., 
$(f(x(1)),f(x(2)))\hmult (f(y(1)),f(y(2)))$ is defined, and, in this case, 
$\hm{(\squaref{f})\mbf{x}}{(\squaref{f})\mbf{y}}
=(f(x(1)),f(x(2)))\hmult (f(y(1)),f(y(2)))
=(f(y(1)),f(x(2)))=(\squaref{f})(\hm{\mbf{x}}{\mbf{y}})$.
\end{proof}

\subsection{Vertical Multiplication and Finite Products 
in the Category of Partial Magmas}\label{s:vm} 
In the following result we use the notation for the forgetful functor
introduced in Section~\ref{s:category}.
\begin{lemma}
The product $B_1\times B_2$ of partial magmas 
$B_1$ and $B_2$ exists in the category of partial magmas. Its support is 
the cartesian product of $\ffSet{(B_1)}$ and $\ffSet{(B_2)}$. 
The natural projections $p_k:B_1\times B_2\to B_k$ 
defined by 
$p_k(\mbf{x})\eqdef x(k)$, $k=1,2$,  are partial magma homomorphisms, and for each partial magma 
$D$ and every pair of partial magma homomorphisms 
$f_k:D\to B_k$, $k=1,2$, there exists a unique partial magma homomorphism
$h:D\to B_1\times B_2$ which makes the following diagram commutative. 
\begin{equation}
\begin{tikzcd}[row sep=large, column sep=large] 
\mbox{}
&
D\arrow[dl, "f_1" ']\arrow[dr, "f_2"]\arrow[d,"h"]
&
\mbox{}
\\
{B_1}
&
{B_1}\times {B_2}\arrow[l, "p_1"]\arrow[r, "p_2" ']
&
{B_2}
\end{tikzcd} 
\label{e:diagramproduct}
\end{equation}
\label{l:productofpm}
\end{lemma}

\begin{proof}
The product of 
$\mbf{x}$ and $\mbf{z}$, where 
$x(k)\in B_k$ and 
$z(k)\in B_k$, for $k=1,2$,
is defined if and only if the products $x(k)\cdot z(k)$ 
are defined in $B_k$, and then 
\begin{equation}
\fbox{$\vm{\mbf{x}}{\mbf{z}}\eqdef(x(1)\cdot z(1),x(2)\cdot z(2))$}
\label{e:termwise} 
\end{equation}
If $\vm{\mbf{x}}{\mbf{z}}$ is defined in $B_1\times B_2$
then, by definition, $x(k)\cdot z(k)$ is defined in $B_k$, 
and $p_k(\vm{\mbf{x}}{\mbf{z}})=x(k)\cdot z(k)
=
p_k(\mbf{x})\cdot p_k(\mbf{z})$, hence $p_k$ is a partial magma homomorphism.
If $f_k:D\to B_k$ are partial magma homomorphisms, for $k=1,2$, then 
the unique function $h:D\to B_1\times B_2$ which makes 
the diagram~\eqref{e:diagramproduct} commutative is defined by setting 
$h(a)\eqdef(f_1(a),f_2(a))$ for each $a\in D$.
If $a\cdot b$ is defined in $D$, hence 
$f_k(a)\cdot f_k(b)$ is defined  
in $B_k$  and is equal to $f_k(a\cdot b)$, for each $k=1,2$
thus $\vm{h(a)}{h(b)}$ is defined in $B_1\times B_2$ and is equal to 
$h(a\cdot b)$. We have thus proved that $h$ is a partial magma homomorphism.
\end{proof}
\begin{notation}
The termwise product defined in~\eqref{e:termwise} 
is called
\textit{vertical multiplication}. 
\end{notation}

\begin{corollary}
If $B$ is a partial magma then $\csquare{B}$ is a partial magma under 
vertical multiplication.
\label{c:theproductisapartialmagma}
\end{corollary}
It is useful to conceive vertical multiplication $\vm{\mbf{x}}{\mbf{y}}$ 
as arising from placing 
$(y(1),y(2))$ on top of $(x(1),x(2))$ and reading top-down and right-to-left, 
as in the following 
diagram. 
\begin{equation}
\begin{tikzcd}[mymatr, arrow style=tikz, 
>={Triangle[length=3mm]}, 
cells={nodes={inner sep=2mm}}, 
row sep=0.4cm,column sep=1.5cm]
\ast \arrow[d, "y(1)" ',phantom] 
&
\ast \arrow[d, "y(2)", phantom] 
\\
\diamond\arrow[d, "x(1)" ', phantom] 
&
\diamond\arrow[d, "x(2)", phantom ] 
\\
\mbox{}
&
\mbox{}
\end{tikzcd}
\quad
\parrow
\quad
\begin{tikzcd}[font=\small,
mymatr,  
arrow style=tikz, 
>={Triangle[length=3mm]}, 
cells={nodes={inner sep=3mm}}, 
row sep=1cm, 
column sep=1.5cm]
\circ\arrow[d, "x(1)\cdot y(1)" ', phantom]
&
\circ\arrow[d, "x(2)\cdot y(2)", phantom]
\\
\mbox{}
&
\mbox{}
\end{tikzcd}
\label{e:vm}
\end{equation}

\subsubsection{Double Partial Magmas}\label{s:dpm}
In Section~\ref{s:hm} we have seen that if $\Setg$ is a nonempty set then 
the set $\csquare{\Setg}$ is endowed with the structure of a partial magma, under 
horizontal 
multiplication, denoted by $\hmult$, and the resulting partial magma is called the 
\textit{twin} partial magma associated to $\Setg$, denoted by 
$\twinf{\Setg}$. This result, together with Corollary~\ref{c:theproductisapartialmagma}, 
implies that 
if $B$ is a partial magma then 
the partial magma $\csquare{B}$ (under vertical multiplication, 
denoted by $\vm{\mbox{}}{\mbox{}}$)
is also endowed with a second structure of partial magma, 
given by horizontal multiplication. 
The relation between these two partial magma structures, given 
by~\eqref{e:il}, is called the 
\textit{interchange law}.  
\begin{lemma}
If $B$ is a partial magma then the following identity holds, for each 
$\mbf{x},\mbf{z},\mbf{x'},\mbf{z'}$ in $\csquare{B}$, 
provided both sides are defined:
\begin{equation}
\fbox{$\hm{(\vm{\mbf{x'}}{\mbf{z'}})}{(\vm{\mbf{x}}{\mbf{z}})}
=\vm{(\hm{\mbf{x'}}{\mbf{x}})}{(\hm{\mbf{z'}}{\mbf{z}})}$}
\label{e:il}
\end{equation}
\label{l:doublepartialmagma}
\end{lemma}
\begin{proof}
If both sides of~\eqref{e:il} are defined, they 
are equal to $(\cdm{x(1)}{z(1)},\cdm{x'(2)}{z'(2)})$. 
\end{proof}
\begin{definition}
A partial magma $(D,\abstractvm{\mbox{}}{\mbox{}})$
is said to be a 
\textit{double partial magma} 
under
a second partial operation $\abstractcdm{\mbox{}}{\mbox{}}$ 
if, for each $x,x',z,z'\in D$
\begin{equation}
\fbox{$\abstractvm{(\abstractcdm{{x'}}{{z'}})}{(\abstractcdm{{x}}{{z}})}
=\abstractcdm{(\abstractvm{{x'}}{{x}})}{(\abstractvm{{z'}}{{z}})}$}
\label{e:il2}
\end{equation}
provided both sides are defined.
Then  the partial magma $(D,\abstractcdm{\mbox{}}{\mbox{}})$ 
is a double partial magma under $\abstractvm{\mbox{}}{\mbox{}}$.
\end{definition}
\begin{corollary}
If $B$ is a partial magma, then 
$\csquare{B}$ is a double partial magma under horizontal multiplication, and 
$\twinf{\ffSet{B}}$ is a double partial magma under vertical multiplication.
\label{c:doublepartialmagma}
\end{corollary}
\begin{proof}
It suffices to apply Lemma~\ref{l:doublepartialmagma}.
\end{proof}
\begin{corollary}
If $\sCat$ is a small category then $\twinf{\sCatA{\sCat}}$ 
is a double partial magma under vertical multiplication, and 
$\csquare{\sCatApm{\sCat}}$ is a double partial magma under horizontal multiplication. 
\label{c:Arrowsaadpm}
\end{corollary}
\begin{proof}
It suffices to apply Lemma~\ref{l:afapm} and Corollary~\ref{c:doublepartialmagma}. 
\end{proof}

\subsection{Other Examples of Partial Magmas Homomorphisms} 
In the following result we will use the 
notation from Section~\ref{s:compositionismultiplication} 
and the notions introduced in~\ref{s:vm}. 

\begin{lemma}
The assignments 
$B\mapsto \csquare{B}$
and
$f\mapsto \squaref{f}$
described in
Corollary~\ref{c:theproductisapartialmagma}
and in Lemma~\ref{l:squarefunctor}
determine a functor $\sPM\to\sPM$ which, with a slight abuse of language, 
will be denoted with the same symbol used in~\eqref{e:twinfunctor}, i.e., as follows:
\begin{equation}
\square:\sPM\to\sPM
\label{e:BISthetwinfunctor} 
\end{equation}  
\label{l:functorSetTwinVertical}
\end{lemma}
\begin{proof}
Lemma~\ref{l:squarefunctor} implies that it suffices to show that, 
if $f:Q\to B$ is a partial magma homomorphism, then 
$\squaref{f}:\csquare{Q}\to\csquare{B}$ is a partial magma homomorphism.
If $\mbf{x},\mbf{y}\in\csquare{Q}$ and
$\vm{\mbf{x}}{\mbf{y}}$ is defined in $\csquare{Q}$
then 
$x(1)\cdot y(1)$ and  $x(2)\cdot y(2)$ are defined in $Q$, hence 
$f(x(1))\cdot{}f(y(1))$
and
$f(x(2))\cdot{}f(y(2))$ are defined in $B$. Thus
$\vm{(f(x(1)),f(x(2)))}{(f(y(1)),f(y(2)))}$ is defined
in $\csquare{B}$, i.e., 
$\vm{(\squaref{f}(\mbf{x}))}{(\squaref{f}(\mbf{y}))}$ is defined in 
$\csquare{B}$. 
Since 
$f$ is a partial magma homomorphism, 
$f(x(1)\cdot y(1))=f(x(1))\cdot f(y(1))$
and
$f(x(2)\cdot y(2))=f(x(2))\cdot f(y(2))$, and hence  
\begin{equation*}
\begin{split}
\vm{(\squaref{f}(\mbf{x}))}{(\squaref{f}(\mbf{y}))}&= \vm{(f(x(1)),f(x(2)))}{(f(y(1)),f(y(2)))}
=
({f(x(1))\cdot f(y(1))},{f(x(2))\cdot f(y(2))})\\
         &=(f[x(1)\cdot y(1)],f[x(2)\cdot y(2)])
         =\squaref{f}({x(1)\cdot y(1)},{x(2)\cdot y(2)})
         =\squaref{f}(\vm{\mbf{x}}{\mbf{y}}) \, .
\end{split}
\end{equation*}
Hence $\squaref{f}:\csquare{Q}\to \csquare{B}$ is a partial magma homomorphism.
\end{proof}

\subsection{Twin Homomorphisms in Partial Magmas}\label{s:internalhom} 

We now introduce the useful notion of 
\textit{twin homomorphism}. 

\begin{definition}
If $x,y$ belong to a partial magma $B$,
then we define 
\begin{equation}
\sihom{B}{x}{y}
\eqdef
\{
\mbf{z}\in \csquare{B}:
z(2)\cdot x
\text{ and }
y\cdot z(1)
\text{ are defined, and  }
z(2)\cdot x=y\cdot z(1)
\}
\label{e:internalhomomorphism}
\end{equation}
and say that $\mbf{z}$ is a 
\textit{twin homomorphism} from $x$ to $y$.
\end{definition}
\begin{remark}
The notion of \textit{twin homomorphism} 
is a kind of “external” variant of the notion of homomorphism, since it 
hinges on the partial magma structure of the ambient space $B$ 
rather than on some “internal” structure of $x$ and $y$, and it may be used whenever 
the latter is not available.   
We will see that, under certain conditions, the set in~\eqref{e:internalhomomorphism} is indeed the 
set of arrows from $x$ to $y$ in a certain category.
We may represent the notion of twin homomorphism in the following diagram.
See Section~\ref{s:twincategory}.
\begin{equation}
\begin{tikzcd}
[mymatr,  
arrow style=tikz, 
>={Triangle[length=1mm]}, 
cells={nodes={inner sep=2mm}}, 
row sep=0.7cm, 
column sep=0.5cm,
every arrow/.append style = -{Latex[scale=0.5]}]
x \arrow[d, shift left, "z(2)"]
\arrow[d, shift right, "z(1)" ']\\
y 
\end{tikzcd} 
\end{equation}
\label{remark:twin} 
\end{remark}

\subsection{Partial Magmas vs. Magmas: The Role of Unital Elements}
If 
multiplication is defined \textit{for every} pair of elements of $B$,
and hence we actually have a \textit{magma}, as 
in Definition~\ref{d:magma}, then 
we drop the adjective “partial”; 
for extra emphasis, 
we may 
add the adjective “total”.
We will see momentarily that 
there are 
conditions, based on 
the notion of \textit{unit} in a partial magma, that 
\textit{force} a partial magma to be a magma. 
When dealing with units in a partial magmas, we have to use the indefinite 
article, since they are not necessarily unique.

\subsubsection{Unital Partial Magmas}
\begin{definition}
An element $x$ of a 
partial magma
is a \textit{unit}
if $x\cdot x$ is defined and, for every element $y$,
\begin{description}
\item[(U.1)]
If $x\cdot y$ is defined, then 
 $x\cdot y=y$;

\item[(U.2)] 
If  $y\cdot x$ is defined, then 
$y\cdot x=y$.
\end{description}
In particular, $x\cdot x=x$. 
A \textit{unital partial magma} is a partial magma which 
has at least one unit.
We denote 
by $\runits{B}$ the set of units in $B$.
\end{definition}
\begin{lemma}
If a partial magma is a magma, then it contains at most one unit.  
\label{l:basiclemma}
\end{lemma}
\begin{proof}
If $x$ and $y$ are units, then, since $x\cdot y$ is defined, $x=x\cdot y=y$. 
\end{proof}

\subsubsection{Unital Homomorphisms of Partial Magmas}\label{s:unital} 
The notion of 
\textit{unital homomorphism of partial magmas} 
is a natural one.
\begin{definition}
If $B$ and $Q$ are unital partial magmas, a function $\sFunctor:Q\to B$ is a 
\textit{unital homomorphism of partial magmas} if 
it is a homomorphism of partial magmas and if the following condition holds: 
\begin{description}
\item[(U)]
If $x\in\runits{Q}$ 
then $\sFunctor x\in\runits{B}$.
\end{description}
\end{definition}

\subsubsection{The Category of Unital Partial Magmas} 
\begin{lemma}
Unital partial magmas are the objects of a category, 
denoted by $\sUPM$, whose arrows are the unital partial magma homomorphisms.
\end{lemma}
\begin{proof}
On the basis of Lemma~\ref{l:PM}, it suffices to show that  
the identity $1_B:B\to B$ maps units into units, and that the composition 
of unital partial magma homomorphisms is also a unital partial magma homomorphism. 
Both statements follow at once from the definitions.  
\end{proof}
The set of unital homomorphisms of partial magmas from $Q$ to $B$
is denoted by 
$$
\cathom{\sUPM}{Q}{B}
$$

\subsubsection{Examples of Unital Partial Magmas}\label{s:egunital}

\paragraph{(Natural Numbers)} Natural numbers $\{0,1,2,\ldots\}$ 
under subtraction form a unital partial magma.

\paragraph{(Horizontal Multiplication)} If $\Setg$ is a nonempty set then the units in 
$\twinf{\Setg}$
(under horizontal multiplication) are the elements $\mbf{x}$ with $x(1)=x(2)$.

\paragraph{(Vertical Multiplication)} If $B$ is a unital partial magma then 
$\csquare{B}$ 
is a unital partial magma under vertical multiplication. Its units are the elements 
$\mbf{x}\in \csquare{B}$ where $x(1),x(2)\in\runits{B}$.

\paragraph{(Matrices)} The units in the partial magma  $\scMatrices_{\ZZ}$ are the identity matrices. 

\paragraph{(Arrows in a Category)} If $\sCat$ is a small category then 
the 
partial magma $\sCatApm{\sCat}$ is unital: The units in $\sCatApm{\sCat}$
are the identity arrows associated to the objects of $\sCat$. 

\subsubsection{Examples of Unital Partial Magma Homomorphisms}\label{s:egunitalpmh}

In the following two results we will use the 
notation from Section~\ref{s:compositionismultiplication} 
and the notions introduced in Section~\ref{s:hm} and~\ref{s:vm}. 

\begin{lemma}
The functor 
$\sTwin$ 
defined in Lemma~\ref{l:functorSetTwinThisOne}
yields a functor (denoted by the same symbol)
$$
\sTwin:\sSet\to\sUPM
$$
\label{l:unitaltwin} 
\end{lemma}
\begin{proof}
It suffices to show that 
$\squaref{f}:\twinf{\Setuno}\to\twinf{\Setdue}$
is unital whenever 
$f:\Setuno\to\Setdue$ is a function from $\Setuno$
to
$\Setdue$.
If $\mbf{x}$ is a unit in $\twinf{\Setg}$, then
$x(1)=x(2)$, hence
$f(x(1))=f(x(2))$, and this means that 
$\squaref{f}(\mbf{x})=(f(x(1)),f(x(2)))$
is a unit in $\twinf{\Setdue}$.
\end{proof}
\begin{lemma}
The functor $\square$ 
defined in Lemma~\ref{l:functorSetTwinVertical} 
yields a functor (denoted by the same symbol)
$$
\square:\sUPM\to\sUPM \, .
$$
\end{lemma}
\begin{proof}
If $f:Q\to B$ is a unital partial magma homomorphism and  
$\mbf{x}\in\runits{\csquare{Q}}$
then 
$x(1),x(2)\in\runits{Q}$, hence 
$f(x(1))$ and 
$f(x(2))$ belong to $\runits{B}$; 
thus $(\squaref{f})\mbf{x}=(f(x(1)),f(x(2)))$
is a unit in $\csquare{B}$.
\end{proof}

\subsection{Tauberian Conditions}

A partial 
 magma which contains only one unit is not necessarily a magma. 
 For example, $\NN$ is a unital partial magma under subtraction, 
 and contains only one unit, but it is not a magma under subtraction. 
However, we will introduce two \textit{Tauberian} conditions, called \textbf{(A)}
and
\textbf{($\mathbf{\Phi}$)}, 
under which the following result, whose proof will be given momentarily, holds. 
\begin{lemma}
If a unital partial magma $B$ satisfies conditions 
\textup{\textbf{(A)}}
and
\textup{\textbf{($\mathbf{\Phi}$)}}, 
then
$$
\fbox{$B$ has only one unit.}
\Longleftrightarrow
\fbox{$B$ is a magma}
$$ 
\label{l:equivalence}
\end{lemma}
\begin{proof}
The proof will be given in Section~\ref{s:proofoflequivalence}
\end{proof}

\subsubsection{Associative Partial Magmas}

Condition~\textbf{(A)} says that the partial magma is \textit{associative}.

\begin{definition}[Condition ({\bf A})]
An \textit{associative partial magma} is a partial magma 
with  the following properties:
\begin{description}
\item[(A.1)] 
Whenever 
$x,y, z$ are elements of the partial magma, the following conditions are equivalent: 
\begin{enumerate}
\item the product $(x\cdot y)\cdot z$ is defined;
\item the product $x\cdot (y\cdot z)$ is defined;
\item the products $x\cdot y$ and $y\cdot z$ are both defined.
\end{enumerate}
\item[(A.2)] 
If any of the previous conditions holds, then 
$(x\cdot y)\cdot z=x\cdot (y\cdot z)$
\end{description}
As in~\cite[p.9]{MacLane1978bis} and 
\cite[p.237]{EilenbergMacLane}
the \textit{triple product} which appears in \textbf{(A.2)}
will be written $x\cdot y\cdot z$. 
\end{definition}
\subsubsection{Examples of Associative and Non-Associative Partial Magmas}\label{s:associative}

\paragraph{(Natural Numbers)}
The natural numbers under subtraction are a non-associative 
partial magma. 

\paragraph{(Horizontal Multiplication)} If $\Setg\not=\emptyset$ then $\twinf{\Setg}$
is an associative partial magma.

\paragraph{(Vertical Multiplication)} If $B$ is an \textit{associative} partial magma then 
$\csquare{B}$ is
associative. 

\paragraph{(Matrices)} The partial magma  $\scMatrices_{\ZZ}$ 
is associative.

\paragraph{(Arrows of a Category)} If $\sCat$ is a small category then the 
partial magma $\sCatApm{\sCat}$, under composition of arrows, 
is associative. Indeed, this is one of the axioms of the notion of category.

\subsubsection{Fastened Partial Magmas}
Condition~\textbf{($\mathbf{\Phi}$)} says that the partial magma is \textit{fastened}.
See \cite[p. 9]{MacLane1978bis} or \cite[pp.237-238]{EilenbergMacLane}.
\begin{definition}[Condition \boldmath ($\Phi$)]
A unital partial magma is \textit{fastened} if 
the following conditions hold:
\begin{description}
\item[($\mathbf{\Phi}$.1)] For every $x$ in $B$ there exists 
$l\in\runits{B}$
such that 
$l\cdot x$ is defined.
\item[($\mathbf{\Phi}$.2)] For every $x$ in $B$ there exists
$r\in\runits{B}$ 
such that 
$x\cdot r$ is defined.
\end{description}
If \textbf{($\mathbf{\Phi}$.1)} holds, we say that 
$l$ is a \textit{left pin} of $x$, or that
$x$ is \textit{fastened to the left} (by $l$).
If \textbf{($\mathbf{\Phi}$.2)} holds, we say that 
$r$ is a \textit{right pin} of $x$, or that
$x$ is \textit{fastened to the right} (by $r$).
\end{definition}

\subsubsection{Proof of Lemma~\ref{l:equivalence}.}\label{s:proofoflequivalence}
\begin{proof}
Let us assume that $B$ has only one unit, say $u_{B}$. 
Then \textbf{($\mathbf{\Phi}$)} means that, 
for each $x\in B$,
$u_{B}$ 
is both a left pin and a right pin of $x$, 
i.e., $u_{B}\cdot x$ and $x\cdot u_{B}$ are both defined, hence
$u_{B}\cdot x=x=x\cdot u_{B}$. 
If $y$ is another element of $B$ then $x\cdot u_{B}$
and $u_{B}\cdot{}y$ are both defined, hence \textbf{(A)} implies that 
$x\cdot u_{B}\cdot{}y$ is defined, 
and since it is equal to $x\cdot{}y$, the conclusion follows.
The other implication follows at once from Lemma~\ref{l:basiclemma}. 
\end{proof}

\subsubsection{Examples of Fastened and Non-Fastened Partial Magmas}\label{s:fastened}

\paragraph{(Natural Numbers)} In the 
natural numbers $\NN\eqdef\{0,1,2,\ldots\}$ under subtraction 
$0$ is the unique unit, and each $x\in\NN$
is fastened to the right (by $0$), but only $0$ is fastened to the left. 

\paragraph{(Horizontal Multiplication)} If $\Setg\not=\emptyset$ then 
$\twinf{\Setg}$ under horizontal multiplication is a fastened partial magma, since each 
$\mbf{x}\in\twinf{\Setg}$ is fastened to the 
left by $(x(1),x(1))$ and to the right by $(x(2),x(2))$.

\paragraph{(Vertical Multiplication)} If $B$ is a \textit{fastened} partial magma then 
$\csquare{B}$ is \textit{fastened}
with respect to 
vertical multiplication. For example, if $\mbf{x}$ is in $\csquare{B}$ then 
$\mbf{x}$ is fastened to the left by $\mbf{u}$
provided 
$x(1)$ is fastened to the left by $u(1)$
and
$x(2)$ is fastened to the left by $u(2)$.

\paragraph{(Matrices)} The space of matrices $\scMatrices_{\ZZ}$ is fastened. Indeed, any $n\times k$ matrix 
is fastened to the left by the identity $n\times n$ matrix 
and to the right by the identity $k\times k$ matrix.

\paragraph{(Arrows of a Category)} The unital and associative partial magma 
$\sCatApm{\sCat}$ of all the arrows in a small category 
$\sCat$ 
is fastened. Indeed, any arrow in $\sCat$ is fastened to the right 
by the identity arrow associated to its domain, 
and to the left by the identity arrow associated to its codomain.

\subsection{Regular Partial Magmas}\label{s:rpm}

The interest of 
unital, associative, and fastened partial magmas 
is not limited to the property  in 
Lemma~\ref{l:equivalence}, as we will see momentarily. 
\begin{definition}
A \textit{regular partial magma} is a unital, associative and fastened partial magma.
A \textit{monoid} is a regular partial magma which has only one unit.
\end{definition}
Lemma~\ref{l:equivalence} implies that a monoid is a (total) magma. 
Observe that every group is a monoid.

\subsubsection{The Category of Regular Partial Magmas}\label{s:crpm} 
\begin{lemma}
Regular partial magmas are the objects of a category, 
denoted by $\sRPM$, whose arrows are the unital partial magma homomorphisms.
\end{lemma}
\begin{proof}
Regular partial magmas are, in particular, unital, hence 
$\sRPM$ is a full subcategory of $\sUPM$, and  
$$
\cathom{\sRPM}{Q}{B}=\cathom{\sUPM}{Q}{B}
$$
for each regular partial magmas $Q$ and $B$.\end{proof}

\subsubsection{Examples of Regular and Non-Regular Partial Magmas}\label{s:fctrpm}
\paragraph{(Natural Numbers)}
The natural numbers 
$\{0,1,2,\ldots\}$ under subtraction are a unital, non-associative, 
non-fastened, hence non-regular partial magma. 

\paragraph{(Horizontal Multiplication)} If $\Setg\not=\emptyset$ then 
$\csquare{\Setg}$ under horizontal multiplication is a regular partial magma. 

\paragraph{(Vertical Multiplication)} If $B$ is a regular partial magma then 
$\csquare{B}$ is a regular partial magma with respect to 
vertical multiplication. 

\paragraph{(Matrices)} The space of matrices $\scMatrices_{\ZZ}$ is a regular partial magma. 

\paragraph{(Arrows of a Category)} From the results in 
Section~\ref{s:egunital}, Section~\ref{s:associative}, and Section~\ref{s:fastened}, 
it follows that 
if $\sCat$  is a small category then 
$\sCatApm{\sCat}$ 
is a regular partial magma.

\subsubsection{Functors and Unital Partial Magma Homomorphisms}\label{s:CatRPM}
In Section~\ref{s:fctrpm} we have seen 
that each small category   $\sCat$ yields a regular partial magma 
$\sCatApm{\sCat}$, 
whose elements are the arrows of $\sCat$. In this section we show that 
we may express this construction as a functor
\begin{equation}
\sfromctpm:\scCat\to\sRPM
\label{e:fctrpmbis} 
\end{equation}
where $\scCat$ is the category of all small categories, as in \cite[Ch.1]{MacLane1978bis}, and $\sRPM$ is the category of regular partial magmas, introduced in 
Section~\ref{s:crpm}. 
\begin{lemma}
The assignment $\sCat\mapsto\sCatApm{\sCat}$ is the object-function of a 
functor $\sfromctpm:\scCat\to\sRPM$.  
\end{lemma}
\begin{proof}
In order to show that~\eqref{e:fctrpmbis} yields a functor, we first have to 
describe its action on arrows in $\scCat$, i.e., on functors between small categories. 
If $\sCatone$ and $\sCattwo$ are small categories and 
$f:\sCatone\to\sCattwo$ is a functor from $\sCatone$ to $\sCatone$, then 
$f$ is determined by a pair of functions 
$f_1:\sCatO{\sCatone}
\to\sCatO{\sCattwo}$
and
$f_2:\sCatA{\sCatone}\to\sCatA{\sCattwo}$
such that {\bf (i)} for each $u\in\sCatO{\sCatone}$, 
$f_2(1_u)=1_{f_1(u)}$,
and {\bf (ii)}
for each $x,y\in\sCatA{\sCatone}$, 
if $x\cdot y$ is defined in $\sCatone$, then 
$f_2(x)\cdot f_2(y)$ is defined in $\sCattwo$ and 
$f_2(x\cdot y)=f_2(x)\cdot f_2(y)$. We define 
$\sCatApm{f}:\sCatApm{\sCatone}\to\sCatApm{\sCattwo}$
by $\sCatApm{f}(x)\eqdef f_2(x)$ for each $x\in\sCatApm{\sCatone}$. 
Then {\bf (ii)} implies that $\sCatApm{f}$ is a partial magma homomorphism, and 
{\bf (i)} that it is unital. In particular, if $f:\sCat\to\sCat$ is the identity functor, then 
$\sCatApm{f}:\sCatApm{\sCat}\to\sCatApm{\sCat}$ is the identity, hence 
$\sCatApm{1_{\sCat}}=1_{\sCatApm{\sCat}}$. Finally, if 
$f:\sCatone\to\sCattwo$ and $g:\sCattwo\to\sCatthree$ are functors, then 
$\sCatApm{g\cdot f}(x)={(g\cdot f)}_2(x)=g_2(f_2(x))=\sCatApm{g}\cdot
\sCatApm{f}(x)$ for each $x\in\sCatApm{\sCatone}$, hence 
$\sCatApm{g\cdot f}=\sCatApm{g}\cdot\sCatApm{f}$, and the proof that 
$\sfromctpm$ is a functor is thus complete.
\end{proof}

The idea behind~\eqref{e:fctrpmbis} is to look at 
a category as an “arrows only” structure, as MacLane puts it 
in \cite{MacLane1978bis}. 
This idea is already contained in 
\cite[p.238]{EilenbergMacLane}, where the authors write that 
\begin{quote}
It is thus clear that the objects play a secondary role, and could be entirely omitted from the definition of category. However, the manipulation of the applications would be 
slightly less convenient were this done. 
\end{quote}
As we will see, the added convenience of looking 
at a category as an “arrows only” structure (i.e., 
as a regular partial magma) is that, once  categorical 
notions are 
expressed in the language 
of partial magmas, it becomes possible to give 
a precise meaning to the 
idea that a natural transformation is a 
\textit{homomorphism of homomorphisms}
(an idea to which,
to the best of our knowledge, 
no precise meaning beyond that of a suggestive 
but insubstantial 
analogy is given in the literature). 
The reason for this added possibility is precisely that, once we identify 
a category $\sCat$ with the regular partial magma $\fctpm{\sCat}$, 
then the objects in the category  $\sCat$ 
are none other than the units in  $\fctpm{\sCat}$, 
and hence the two-layered structure
 of the original definition of category,
articulated in  
\textit{objects}
and
\textit{arrows}, seen as ontologically distinct, fades away.

We will show that $\sfromctpm$ is indeed an isomorphism of categories. In order to 
prove it, we will need some preliminary results.

\subsection{Properties of Regular Partial Magmas}\label{s:rpmsc}

Firstly, we show that in a regular partial magma the elements $l$ and $r$ in 
\textbf{($\mathbf{\Phi}$)} are unique; see \cite[p.9]{MacLane1978bis}.

\subsubsection{Projections Onto the Unit Elements}\label{s:coherence1}

\begin{lemma}
If $B$ is a regular partial magma, then each 
$x\in B$ has a unique left pin and a unique right pin.
\label{l:uniquefasteningunits}
\end{lemma}

\begin{proof}
If $x$ is fastened to the right by $r$ and by 
$\rho$, then 
the fact that 
${}x{}={}x{}\cdot r={}x{}\cdot\rho$ implies that 
$({}x{}\cdot{}r)\cdot\rho$
is defined (and is equal to ${}x{}$), hence \textbf{(A)} implies that 
$r\cdot\rho$ is defined, hence 
$r=r\cdot\rho=\rho$. The proof of the 
uniqueness of the left pin of $x$ is similar.
\end{proof}
\begin{definition}
If $B$ is a regular partial magma, 
then the functions represented in the following diagram 
\begin{equation}
\begin{tikzcd}[row sep=large, column sep=large] 
\mbox{}
&
B\arrow[dl, "\onedomcod{B}" ']\arrow[dr, "\twodomcod{B}"]\arrow[d, "\domcod{B}"]
&
\mbox{}
\\
{B}
&
\csquare{B}\arrow[l, "p_1"]\arrow[r, "p_2" ']
&
{B}
\end{tikzcd} 
\label{e:domaincodomain}
\end{equation}
are defined as follows.
\begin{itemize}
\item The unique right pin of $x\in B$ is denoted by $\ssdom{x}$
and the function $\onedomcod{B}$ 
is defined by  $\onedomcodof{B}{x}\eqdef\ssdom{x}$.

\item The unique left pin of $x\in B$ is denoted by $\sscodom{x}$
and the function $\twodomcod{B}$ 
is defined by  $\twodomcodof{B}{x}\eqdef\sscodom{x}$.

\item The functions $p_1$ and $p_2$ are
the coordinate projections, hence $p_1(\mbf{x})\eqdef x(1)$
and $p_2(\mbf{x})\eqdef x(2)$.

\item The function $\domcod{B}$ is the unique map which makes the 
diagram~\eqref{e:domaincodomain} commutative. 
\end{itemize}
We call 
$\ssdom{x}$ the \textit{domain} of $x$ \textit{in} $B$.
We call  
$\sscodom{x}$ the \textit{codomain} of $x$ \textit{in} $B$, and write 
\begin{equation}
\arrowws{\ssdom{x}}{\sscodom{x}}{x}\quad\text{ in $B$} 
\label{e:imagery} 
\end{equation}
\label{d:domcod}
\end{definition}
Lemma~\ref{l:uniquefasteningunits}
says that 
$\ssdom{x}$ and $\sscodom{x}$ are uniquely determined,
hence~\eqref{e:imagery} 
enables us to recapture 
the imagery of category theory, although it 
does not imply that $x$ is a \textit{function}
from $\ssdom{x}$ to $\sscodom{x}$.  
We now describe the outcome 
of~\eqref{e:imagery}
when $x$ is a unit in $B$. 
\begin{lemma}
If $B$ is a regular partial magma and $x\in\runits{B}$ then 
$\ssdom{x}=\sscodom{x}=x$, i.e., 
$$
\arrowws{x}{x}{x}\quad\text{ \textup{in} $B$} 
$$ 
\label{l:onto}
\end{lemma}
\begin{proof}
For each $x\in{B}$,
$\sscodom{x}\in\runits{B}$
and $\sscodom{x}\cdot x$ is defined, and if 
$x\in\runits{B}$
then 
$\sscodom{x}=\sscodom{x}\cdot x=x$, thus 
$\sscodom{x}=x$. 
A similar proof shows that $\ssdom{x}=x$.  
\end{proof}
The following result confirms that 
the imagery in~\eqref{e:imagery} is consistent with category theory, 
since it shows that 
the maps
$\onedomcod{B}$ and $\twodomcod{B}$,
described in~\eqref{e:domaincodomain},
are (idempotent) projections onto the unit elements, and 
behave in the familiar way under multiplication. 
\begin{proposition}
If $B$ is a regular partial magma,
then the following results hold.
\begin{description}
\item[(range)] The range of the maps 
 $\onedomcod{B}$
and
$\twodomcod{B}$
in~\eqref{e:domaincodomain}
is $\textup{$\runits{B}$}$.
\item[(idempotent)]
The maps 
 $\onedomcod{B}$
and
$\twodomcod{B}$
are idempotent.
\item[(behavior under product)] If $x\cdot y$ is defined, then 
\begin{equation}
\sscodom{({x}\cdot{y})}=\sscodom{x} 
\text{ and }
\,
\ssdom{({x}\cdot{y})}=\ssdom{y}\, .
\label{e:composition}
\end{equation}
\end{description}
\label{p:thesamenotion}
\end{proposition}
\begin{proof}
The first two statements follow from Lemma~\ref{l:onto}.
If 
${x}\cdot{y}$ is defined, 
${u}$ is a unit and ${u}\cdot{x}$ is defined, then it is equal to ${x}$. 
Hence $({u}\cdot{x})\cdot{y}$ is defined, and 
\textbf{(A)} implies that  ${u}\cdot({x}\cdot{y})$ is defined, and 
the first equality in~\eqref{e:composition} follows. 
The proof of the second one is similar and is left to the reader. 
\end{proof}
\subsubsection{Examples of Projections Onto the Unit Elements}
We now examine the maps~\eqref{e:domaincodomain} in the examples 
of regular partial magmas we have seen so far. 

\paragraph{(Horizontal Multiplication)}
If $\Setg$ is a nonempty set and 
$\mbf{x}\in\twinf{\Setg}$ then 
$$
\arrowws{(x(1),x(1))}{(x(2),x(2))}{\mbf{x}} 
\quad
\text{ in $\twinf{\Setg}$} \, .
$$
\paragraph{(Vertical Multiplication)}
If $B$ is a regular partial magma and 
$\mbf{x}\in\csquare{B}$
then 
$$
\arrowws{(\ssdom{x(1)},\ssdom{x(2)})}{(\sscodom{x(1)},\sscodom{x(2)})}{\mbf{x}}
\quad\text{ in $\csquare{B}$}  \, .
$$

\paragraph{(Matrices)}
If $\matrixone\in\scMatrices_{\ZZ}$
is any 
$n\times k$ matrix, and $\ione_n$ ($\ione_k$) is the identity matrix
of order $n$ ($k$, resp.)  then 
$$
\arrowws{\ione_k}{\ione_n}{\matrixone}
\quad\text{ in $\scMatrices_{\ZZ}$} \, .
$$

\paragraph{(Arrows of a Category)}
If $\sCat$ is a small category and $x$ belongs to the partial magma 
$\sCatApm{\sCat}$, then 
$$
\arrowws{\ssdom{x}}{\sscodom{x}}{x}
\quad\text{ in $\,\sCatApm{\sCat}$}  \, .
$$
if and only if, 
\textit{in the category $\sCat$},
$\ssdom{x}$ is the domain of $x$ 
and $\sscodom{x}$ the codomain of $x$.

\subsubsection{Coherence with Category Theory}\label{s:coherence2}

The following result says  that,  
in a regular partial magma,
the maps described in~\eqref{e:domaincodomain}  
\textit{determine} whether  
multiplication is defined in terms of the \textit{chain rule} that appears 
in Lemma~\ref{l:chainrule}.
Hence the imagery
in~\eqref{e:imagery}
 is 
entirely consistent with category theory. 

\begin{lemma}
If $B$ is  a regular partial magma and 
$x,z\in B$, then 
 the following conditions are equivalent.
\begin{description}
\item[(1)] \textbf{(the product is defined)} The product $x\cdot z$ is defined.
\item[(2)] \textbf{(chain rule)} $\ssdom{x}=\sscodom{z}$.
\item[(3)] \textbf{(chain rule, in a diagram)} There exist units $u,v,w$ in $B$
such that 
\begin{equation}
\twoarrowws{u}{v}{w}{{z}}{{x}}
\label{e:imagery2} 
\end{equation}
\end{description}
\label{l:chainrule}
\end{lemma}
\begin{proof}
It suffices to prove that \textbf{(1)} is equivalent to \textbf{(2)}, since \textbf{(3)} is a diagrammatic 
expression of \textbf{(2)}.
If $x\cdot{z}$ is defined, then 
$x\cdot{z}=x\cdot
(
\sscodom{z}
\cdot{z}
)
$, and 
\textbf{(A)} implies that 
$x\cdot\sscodom{z}$
is defined, hence 
$\sscodom{z}=\ssdom{x}$, by Lemma~\ref{l:uniquefasteningunits}.
Now assume that $\sscodom{z}=\ssdom{x}$
and denote this element by $u$. 
Then $x\cdot u\, (=x)$
and 
$u\cdot{z} \,(={z})$
are defined, hence 
$x\cdot(u\cdot{z})$ is defined,
by \textbf{(A)},
and since it is equal to 
$x\cdot{z}$, the conclusion follows.
 \end{proof}

We now show that the notion 
of \textit{functor} is the categorical version 
of the notion of unital partial magma homomorphism, defined in Section~\ref{s:unital}.   
This result is related to the idea of associating 
a small category to each regular partial magma, thus obtaining a functor 
$\snCat{cat}:\sRPM \to\scCat$, 
which will be introduced momentarily. 

\begin{lemma} Let $Q$ and $B$ be regular partial magmas
and let 
$\sFunctor\in
\cathom{\sRPM}{Q}{B}$.
Then 
for each $x\in Q$
\begin{equation}
\begin{tikzcd}[mymatr, arrow style=tikz, >={Triangle[length=1mm]}, 
cells={nodes={inner sep=2mm}}, row sep=0.7cm, column sep=0.5cm,
every arrow/.append style = -{Latex[scale=0.5]}]
\arrowws{\ssdom{x}}{\sscodom{x}}{x} 
\end{tikzcd} 
\Longrightarrow
\begin{tikzcd}[mymatr, arrow style=tikz, >={Triangle[length=1mm]}, 
cells={nodes={inner sep=2mm}}, row sep=0.7cm, column sep=0.5cm,
every arrow/.append style = -{Latex[scale=0.5]}]
\arrowws{\ssdom{\sFunctor{x}}}{\sscodom{\sFunctor{x}}}{\sFunctor{x}}
\end{tikzcd} 
\label{e:functors}
\end{equation}
\label{l:equivalent}
\end{lemma}
\begin{proof}
In order to prove~\eqref{e:functors}, we have to show that 
\begin{equation} 
\ssdom{(\sFunctor{x})}=
\sFunctor{(\ssdom{x})}
\text{ and }
\sscodom{(\sFunctor{x})}=
\sFunctor{(\sscodom{x})} \, .
\label{e:commutative}
\end{equation} 
Observe that 
$x\cdot \ssdom{x}$ is defined
and
$\ssdom{x}\in\runits{Q}$. 
Since $\sFunctor$ is a unital partial magma homomorphism, it follows that 
$\sFunctor({x})\cdot \sFunctor(\ssdom{x})$ is defined 
and is equal to $\sFunctor(x\cdot \ssdom{x})$, 
and that 
$\sFunctor(\ssdom{x})\in\runits{B}$.
Hence $\sFunctor(\ssdom{x})=
\ssdom{(\sFunctor{x})}$.
The proof of the other equality in~\eqref{e:functors} is similar and is left to the reader. 
\end{proof}
\begin{remark}
In view of~\eqref{e:commutative}, we may write 
$\ssdom{\sFunctor{x}}$
and
$\sscodom{\sFunctor{x}}$
without parentheses, since the two meanings that may be attributed to this expression have the same value. 
\end{remark}

\subsection{Each Regular Partial Magma Yields a Small Category}
In Section~\ref{s:rpm} we have seen 
that each small category   $\sCat$ yields a regular partial magma 
$\sCatApm{\sCat}$, 
whose elements are the arrows of $\sCat$, and that 
this assignment is part of a functor
$\sfromctpm:\scCat\to\sRPM$ 
where $\scCat$ is the category of all small categories, as in \cite[Ch.1]{MacLane1978bis}. 
In order to show that $\sfromctpm$ is an isomorphism of categories, we introduce 
a functor 
\begin{equation}
\snCat{cat}:\sRPM \to\scCat
\label{e:inverse}
\end{equation}
and show that it is a bilateral inverse of $\sfromctpm$.

\begin{lemma}
If $B$ is a regular partial magma, then the small category 
$\fpmtc{B}$ 
is defined as follows:
The objects of $\fpmtc{B}$ are the units of $B$, i.e.,
\begin{equation}
{\sCatO{\fpmtc{B}}\eqdef\runits{B}}
\label{e:definitionofobjects} 
\end{equation}
and the arrows of $\fpmtc{B}$ are the elements of $B$
with its partial magma structure, i.e.,
\begin{equation}
\sCatApm{\fpmtc{B}}\eqdef{}B \, .
\label{e:definitionofarrows} 
\end{equation}
Then $\fpmtc{B}$ is a small category.
\label{l:cfrpm}
\end{lemma}
\begin{proof}
Composition of arrows in $\fpmtc{B}$ amounts to multiplication in $B$, and, 
if $u,v\in\sCatO{\fpmtc{B}}$, then 
\begin{equation}
{\cathom{\fpmtc{B}}{u}{v}\eqdef
\{{x}\in B:
x\cdot u
\,
\text{ is defined and }
v\cdot x
\,
\text{ is defined}
\}} \, .
\label{e:homset} 
\end{equation}
Hence $\cathom{\fpmtc{B}}{u}{v}$
is the collection of $x\in B$ such that 
$\ssdom{x}=u$
and
$\sscodom{x}=v$, as in Definition~\ref{d:domcod}.
Associativity of composition of arrows follows from associativity of multiplication 
in $B$. The other category axioms follow from the results proved in 
Section~\ref{s:rpmsc}. 
\end{proof}
\begin{lemma}
The assignment $B\mapsto\fpmtc{B}$ is the object-function of a 
functor $\snCat{cat}:\sRPM\to\scCat$.  
\end{lemma}
\begin{proof}
In order to show that~\eqref{e:inverse} yields a functor, we first have to 
describe its action on arrows in $\sRPM$, i.e., on unital partial magma homomorphisms.
If $\sFunctor:B\to Q$ is a unital partial magma homomorphism between 
regular partial magmas $B$ and $Q$, then 
$\fpmtc{\sFunctor}:\fpmtc{B}\to\fpmtc{Q}$ is the functor defined by setting 
$\fpmtc{\sFunctor}_{1}:{\sCatO{\fpmtc{B}}}\to{\sCatO{\fpmtc{Q}}}$, where 
$\fpmtc{\sFunctor}_{1}(u)\eqdef \sFunctor(u)$ for each $u\in{\sCatO{\fpmtc{B}}}$, and 
$\fpmtc{\sFunctor}_{2}:\sCatA{\fpmtc{B}}\to\sCatA{\fpmtc{Q}}$, where 
$\fpmtc{\sFunctor}_{2}(x)\eqdef\sFunctor(x)$ for each $x\in\sCatA{\fpmtc{B}}$. 
The fact that $\fpmtc{\sFunctor}:\fpmtc{B}\to\fpmtc{Q}$ is a functor follows at once from the fact that $\sFunctor:B\to Q$ is a unital partial magma homomorphism. 
Observe that if $\sFunctor:B\to B$ is the identity partial magma homomorphism then 
$\fpmtc{\sFunctor}:\fpmtc{B}\to\fpmtc{B}$ is the identify functor, and that if 
$\sFunctorone:B_1\to B_2$ and 
$\sFunctortwo:B_2\to B_3$ are unital partial magma homomorphisms, then 
$\fpmtc{\sFunctortwo\circ\sFunctorone}=\fpmtc{\sFunctortwo}\circ\fpmtc{\sFunctorone}$. 
Hence $\snCat{cat}:\sRPM\to\scCat$ is a functor.
\end{proof}

Recall that  in Section~\ref{s:CatRPM} we have introduced the functor
$\sfromctpm:\scCat\to\sRPM$. 
\begin{proposition}
The functor $\snCat{cat}:\sRPM\to\scCat$ is a bilateral inverse of the functor
$\sfromctpm:\scCat\to\sRPM$.
\label{p:equivalent}
\end{proposition}
\begin{proof}
It suffices to observe that if $B$ is a regular partial magma, then $B=\fctpm{\fpmtc{B}}$
and if $\sCat$ is a category then $\sCat=\fpmtc{\sCatApm{\sCat}}$.
\end{proof}

The following notation will help us draw a useful corollary of 
Proposition~\ref{p:equivalent}.
For given categories $\sCatone$ and $\sCattwo$, 
the space of functors from $\sCatone$ to $\sCattwo$ 
is denoted as follows:
\begin{equation}
\ch{\sCatone}{\sCattwo} \, . 
\label{e:spaceofunctors}
\end{equation}

\begin{corollary}
Let $Q$ and $B$ be regular partial magmas, 
$\sCatone$ and
let $\sCattwo$ be small categories. 
Then 
\begin{equation}
\begin{split}
\ch{\sCatone}{\sCattwo}&\cong\cathom{\sRPM}{\fctpm{\sCatone}}{\fctpm{\sCattwo}} \, . \\
\cathom{\sRPM}{Q}{B}&\cong\ch{\fpmtc{Q}}{\fpmtc{B}} \, .
\end{split}
\label{e:isomorphism}
\end{equation}
\end{corollary}
\begin{proof}
The isomorphisms in~\eqref{e:isomorphism} is given by the functors 
$\sfromctpm$ and $\snCat{cat}$.
\end{proof}
\begin{remark}
With a slight abuse of language, 
in order to avoid too pedantic a notation, 
we will use the same symbol 
to denote a unital partial magma homomorphism 
from $Q$ to $B$
and 
the associated functor from $\fpmtc{Q}$
to
$\fpmtc{B}$, and, by the same token, we will use the same symbol to denote 
a functor from $\sCatone$ to $\sCattwo$ and the corresponding unital partial magma homomorphism from $\fctpm{\sCatone}$ to $\fctpm{\sCattwo}$. 
\label{r:tss} 
\end{remark}

\subsection{Examples of Categories Associated to Partial Magmas}

It is useful to specify the nature of the categories 
$\fpmtc{B}$ when $B$ is one of the regular partial magmas 
described in  
Examples~\ref{eg:examples}. In the following diagrams, 
the identity arrows will not be represented. 
The category $\scOne\eqdef\fpmtc{\{\ione\}}$
has only one object, denoted $\ione$, and no arrows apart from the identity. 
It may 
be represented by the following diagram.
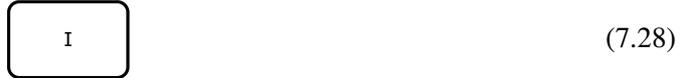
\begin{equation}
\begin{tikzcd}[mymatr, arrow style=tikz, >={Triangle[length=1mm]}, 
cells={nodes={inner sep=2mm}}, row sep=0.7cm, column sep=0.5cm,
every arrow/.append style = -{Latex[scale=0.5]}]
\ione 
\end{tikzcd} 
\end{equation}
The category $\scTwo\eqdef\fpmtc{\{\ione,\ii\}}$
has only two objects and one arrow apart from the identities, 
and it may represented by the following diagram.
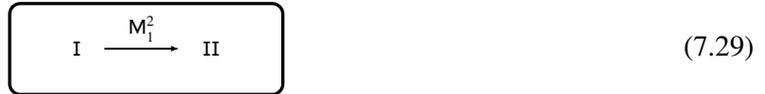
\begin{equation}
\begin{tikzcd}[mymatr, arrow style=tikz, >={Triangle[length=1mm]}, 
cells={nodes={inner sep=3mm}}, row sep=0.8cm, column sep=0.5cm,
every arrow/.append style = -{Latex[scale=0.5]}]
\ione\arrow[rr,"\aidentitydueuno"]&&\ii 
\end{tikzcd} 
\end{equation}
The category 
$\mathbf{II}\eqdef\fpmtc{\{\ione,\ii\}}$
has only two objects and no arrows apart from the identity arrows. 
\begin{equation}
\begin{tikzcd}[mymatr, arrow style=tikz, >={Triangle[length=1mm]}, 
cells={nodes={inner sep=2mm}}, row sep=0.7cm, column sep=0.5cm,
every arrow/.append style = -{Latex[scale=0.5]}]
\ione&&\ii 
\end{tikzcd} 
\end{equation}
The category $\scThree\eqdef\fpmtc{\{\ione,\ii,\iii,\aidentitydueuno,\aidentitytredue,
\aidentitytreuno\}}$
has three objects and three arrows beside the identity arrows. It may be 
represented by the following commutative diagram. 
\begin{equation}
\begin{tikzcd}[mymatr, arrow style=tikz, >={Triangle[length=1mm]}, 
cells={nodes={inner sep=2mm}}, row sep=1cm, column sep=1cm,
every arrow/.append style = -{Latex[scale=0.5]}]
\ione\arrow[rr,"\aidentitydueuno"]
\arrow[drr,"\aidentitytreuno" ']&&\ii\arrow[d,"\aidentitytredue"] \\
\mbox{} &&\iii
\end{tikzcd} 
\end{equation}
The category $\scSQ\eqdef\fpmtc{\{\ione,\ii,\iii,\iv,
\aidentitydueuno,
\aidentitytredue,
\aidentitytreuno,
\aidentityquattrouno,
\aidentitytrequattro
\}}$ has four objects and 
five arrows apart from the identities, and it may be represented 
by the following commutative diagram.
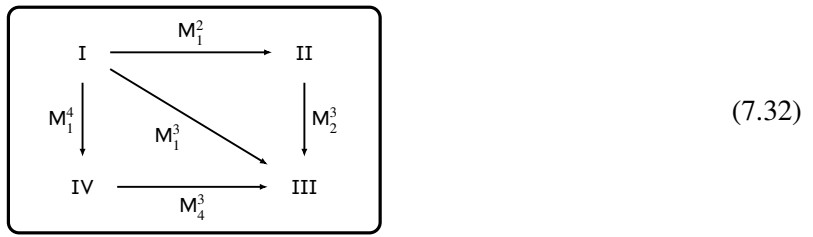
\begin{equation}
\begin{tikzcd}[mymatr, arrow style=tikz, >={Triangle[length=1mm]}, 
cells={nodes={inner sep=3mm}}, row sep=1cm, column sep=1cm,
every arrow/.append style = -{Latex[scale=0.5]}]
\ione\arrow[d,"\aidentityquattrouno" ']\arrow[rr,"\aidentitydueuno"]
\arrow[drr,"\aidentitytreuno" ']&&\ii\arrow[d,"\aidentitytredue"] \\
\iv\arrow[rr,"\aidentitytrequattro" '] &&\iii
\end{tikzcd} 
\end{equation}

\subsection{The Twin Category}\label{s:twincategory}

Recall that if 
$u,v$
are \textit{objects} of a small category $\sCat$,  
then the expression 
\begin{equation}
\cathom{{\sCat}}{x}{y} \, ,
\label{e:ha}
\end{equation}
defined in~\eqref{e:hobjects},
denotes a subset of $\sCatA{\sCat}$, but if $x$ and $y$
are arrows which are not units in $\fctpm{\sCat}$,
then the expression~\eqref{e:ha} 
has no meaning.
We now show that 
the arrows  
in $\sCat$ are 
the \textit{objects} of a new small category, denoted by 
$\twin{\sCat}$, 
called the \textit{twin category} associated to $\sCat$,
where
it makes sense to 
consider 
\begin{equation}
\cathom{\twin{\sCat}}{x}{y}   \, .
\label{e:harrows}
\end{equation}
Indeed, we will see momentarily that 
$\sCatO{\twin{\sCat}}=\sCatA{\sCat}$
and
$\sCatApm{\twin{\sCat}}\subset \csquare{\sCatApm{\sCat}}$.
See Section~\ref{s:internalhom}.
\begin{remark}
The notion of \textit{twin} (or \textit{arrow}) {category} 
is due to F.\ William Lawvere \cite{Lawvere}.  
\end{remark}

\begin{definition}
If $\sCat$ is a small category, then the elements of $\csquare{\sCatApm{\sCat}}$
are called \textit{twin arrows} in $\sCat$. 
If $\mbf{z}$ is a twin arrow in $\sCat$, then 
$\mbf{z}=(z(1),z(2))$, where 
${z}(k)\in\sCatA{\sCat}$ for $k=1,2$.
\end{definition}

\begin{definition}
If $\sCat$ is a small category,
$x,y$ are 
arrows of $\sCat$, 
and
$\mbf{z}$ is a twin arrow in $\sCat$, we say that 
$\mbf{z}$ is a twin arrow \textit{from $x$ to $y$} if the following conditions hold:
\begin{description}
\item[(1)]
$z(2)\cdot x$ 
and
$y\cdot z(1)$ are both defined in $\sCatA{\sCat}$

\item[(2)] $z(2)\cdot x=y\cdot z(1)$
\end{description}
\label{d:dotw}
\end{definition}
\begin{lemma}
The conditions in Definition~\eqref{d:dotw} 
may be expressed by 
the following commutative diagram. 
\begin{equation}
\begin{tikzcd}[row sep=large, column sep=large]
\ssdom{x}
\arrow[r,"x" ']
\arrow[d, "z(1)" ']
 & 
\arrow[d, "z(2)" ]
\sscodom{x}
\\ 
\ssdom{y}
\arrow[r, "y" ']
& 
\sscodom{y}
\end{tikzcd} 
\label{eq:generaldiagram}
\end{equation}
\label{l:meaning}
\end{lemma}
\begin{proof}
Observe that $z(2)\cdot x$
in~\eqref{eq:generaldiagram}
is defined if and only if $\sscodom{x}=
\ssdom{z(2)}$, and 
$y\cdot z(1)$ is defined if and only if 
$\sscodom{z(1)}=\ssdom{y}$. Moreover, if 
$y\cdot z(1)=z(2)\cdot x$, then Proposition~\eqref{p:thesamenotion}
implies that 
$\ssdom{x}=\ssdom{z(1)}$
and $\sscodom{y}=\sscodom{z(2)}$.  
\end{proof}

\begin{proposition}
If $\sCat$ is a small category, then the arrows  
in $\sCat$ are 
the objects of a small category, denoted by $\twin{\sCat}$, 
called the \textit{twin category} associated to $\sCat$,
and defined as follows:
If $x,y$ are 
arrows of $\sCat$ then 
\begin{equation}
\cathom{\twin{\sCat}}{x}{y}
\eqdef
\{\mbf{z}\in\csquare{\sCatApm{\sCat}}:
\mbf{z} \text{ is a twin arrow from $x$ to $y$} 
\}
\label{e:definitionofinternalh}
\end{equation}
\label{p:twin}
\end{proposition}
\begin{proof}
We verify the category axioms described in~\cite[p. 27]{MacLane1978bis}.
The composition
$$
\cathom{\twin{\sCat}}{y}{w}
\times
\cathom{\twin{\sCat}}{x}{y}
\to
\cathom{\twin{\sCat}}{x}{w}
$$
is \textit{vertical} multiplication in $\csquare{\sCatApm{\sCat}}$, induced 
by multiplication in $\sCatApm{\sCat}$; cf. Section~\ref{s:vm}. 
Indeed, if 
$\mbf{z}$ is in 
$\cathom{\twin{\sCat}}{x}{y}$
and
$\mbf{s}$ is in 
$\cathom{\twin{\sCat}}{y}{w}$, then 
the vertical multiplication 
$\vm{\mbf{s}}{\mbf{z}}$ 
is defined and 
belongs to 
$\cathom{\twin{\sCat}}{x}{w}$. 
Now observe that, for each arrow $x$ in $\sCat$, 
$$
\mbf{x}\eqdef(\ssdom{x},\sscodom{x})\in\cathom{\twin{\sCat}}{x}{x}
$$
and if $\mbf{z}\in\cathom{\twin{\sCat}}{x}{y}$ then 
$\vm{\mbf{z}}{\mbf{x}}=
\vm{\mbf{z}}{(\ssdom{x},\sscodom{x})}=
(z(1)\cdot \ssdom{x},z(2)\cdot \sscodom{x})
=(z(1),z(2))=\mbf{z}$,
since $\ssdom{x}$ and $\sscodom{x}$ are units. 
In a similar way one shows that if 
$\mbf{z}\in\cathom{\twin{\sCat}}{y}{x}$
then $\vm{\mbf{x}}{\mbf{z}}=\mbf{z}$.
Associativity for vertical multiplication 
in $\csquare{\sCatApm{\sCat}}$
holds, 
since $\sCatApm{\sCat}$ is associative, 
as observed in Section~\ref{s:associative}.
\end{proof}
Since~\eqref{e:harrows} is defined when $x$ and $y$ 
are arrows of $\sCat$, it also makes sense when $x$
and $y$ are objects of $\sCat$, since objects may be identified with the 
associated identity arrows. 
We now show that, with this identification, 
the hom-sets (between arrows) in $\twin{\sCat}$
 recapture 
the hom-sets (between objects) in $\sCat$.
\begin{lemma}
If $\sCat$ is a small category and $u,v\in\sCatO{\sCat}$, then 
\begin{equation}
\cathom{\sCat}{u}{v}
\cong
\cathom{\twin{\sCat}}{u}{v}
=
\{
(x,x):x\in\cathom{\sCat}{u}{v}
\} 
\end{equation}
\label{l:recaptures}
\end{lemma}
\begin{proof}
Let $\mbf{z}\in\cathom{\twin{\sCat}}{u}{v}$.
Then $\mbf{z}=({z}(1),{z}(2))$
where ${z}(1)$
and
${z}(2)$
belong to $\sCatA{\sCat}$, 
${z}(2)\cdot u$ 
and 
$v\cdot{z}(1)$
are defined and equal. Since 
${z}(2)\cdot u$ 
is defined, 
Lemma~\ref{l:chainrule} implies that  
$\ssdom{z(2)}=\sscodom{u}$. 
Since $u$ is a unit, 
$\sscodom{u}=u$, and, for a similar reason, 
$v=\ssdom{v}=\sscodom{{z}(1)}$. Since $u$ and $v$
are units, the equality
${z}(2)\cdot u= 
v\cdot{z}(1)$ 
implies that 
${z}(1)={z}(2)$.
Let $x\eqdef{z}(1)$. 
Then
$\ssdom{x}=u$, 
$\sscodom{x}=v$, and $\boldsymbol{z}=(x,x)$
with $x\in \cathom{\sCat}{u}{v}$.
\end{proof}
\begin{remark}
The proof of Lemma~\ref{l:recaptures}
may be better appreciated using the following diagram.
\begin{equation}
\begin{tikzcd}[row sep=large, column sep=large]
u
\arrow[r,"u" ']
\arrow[d, "z(1)" ']
 & 
\arrow[d, "z(2)" ]
u
\\ 
v
\arrow[r, "v" ']
& 
v
\end{tikzcd} 
\label{eq:cdfihone}
\end{equation}
\end{remark}
We now examine~\eqref{e:harrows} when only one among 
$x$ and $y$ is an object.
\begin{lemma}
If $\sCat$ is a small category, 
$x=u$ is an object, and $y$ an arrow but not an identity arrow, then 
\begin{equation}
\cathom{\twin{\sCat}}{u}{y} 
=
\{(z(1),z(2)):
z(1)\in\cathom{\sCat}{u}{\ssdom{y}},
\,\,
z(2)\in\cathom{\sCat}{u}{\sscodom{y}},
\,\,
z(2)=y\cdot z(1)
\}    \, .
\label{e:firstcase}
\end{equation}
If $y=v$ is an object and $x$ is an arrow but not an identity arrow then 
\begin{equation}
\cathom{\twin{\sCat}}{x}{v} 
=
\{(z(1),z(2)):
z(1)\in\cathom{\sCat}{\ssdom{x}}{v},
\,\,
z(2)\in\cathom{\sCat}{\sscodom{x}}{v},
\,\,
z(1)=z(2)\cdot x
\}  \, .
\label{e:secondcase}
\end{equation}
\end{lemma}
\begin{proof}
The proof follows the same lines of 
the proof of Lemma~\ref{l:recaptures}. Indeed, it suffices to 
adapt  the diagram~\eqref{eq:generaldiagram} to the case at hand and 
apply Lemma~\ref{l:chainrule}. Hence the proof of~\eqref{e:firstcase}
may be followed along this commutative diagram:
\begin{equation}
\begin{tikzcd}[column sep=small]
& u \arrow[dl, "z(1)" '] \arrow[dr, "z(2)" ] & \\
\ssdom{y} \arrow[rr, "y" '] & & \sscodom{y}
\end{tikzcd} 
\end{equation}
The proof of~\eqref{e:secondcase} may be followed along this
commutative diagram.
\begin{equation}
\begin{tikzcd}[column sep=small]
\ssdom{x}\arrow[rr, "x"] \arrow[dr, "z(1)" ']      &                 & \sscodom{x}\arrow[dl, "z(2)"]\\
		      &   v                                      & 
\end{tikzcd} 
\end{equation}
\end{proof}
\subsection{Horizontal and Vertical Composition of Twin Arrows}
In Corollary~\ref{c:Arrowsaadpm} we have seen that if 
$\sCat$ is a small category then $\twinf{\sCatA{\sCat}}\cong\csquare{\sCatApm{\sCat}}$
is a double partial magma. 
In the proof of Proposition~\eqref{p:twin} 
we have seen that the arrows in the twin category 
associated to a given small category $\sCat$ are the so-called twin arrows in 
$\sCat$, i.e., elements of $\csquare{\sCatApm{\sCat}}$. Vertical multiplication 
in $\csquare{\sCatApm{\sCat}}$
may be represented
in the following diagram, a special case of~\eqref{e:vm},
where for simplicity we only indicate the relevant data and denote 
domains and codomains with small circles.
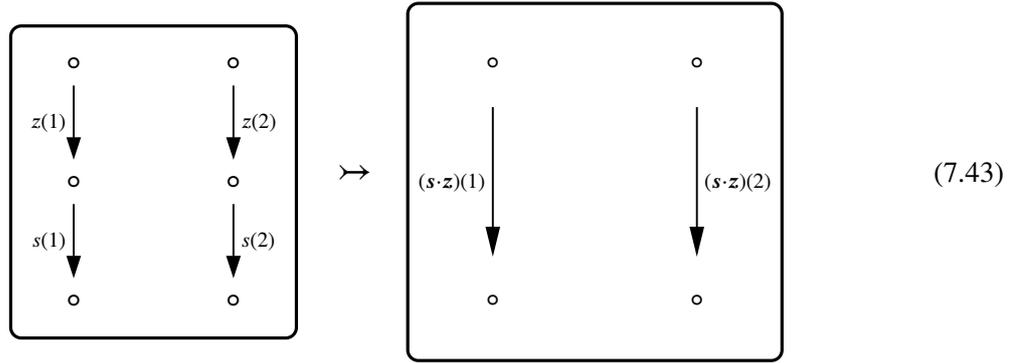
\begin{equation}
\begin{tikzcd}[mymatr, arrow style=tikz, 
>={Triangle[length=3mm]}, 
cells={nodes={inner sep=2mm}}, 
row sep=1cm,column sep=1.5cm]
\circ \arrow[d, "z(1)" '] 
& 
\circ \arrow[d, "z(2)" ] 
\\
\circ\arrow[d, "s(1)" '] 
&
\circ\arrow[d, "s(2)" ] 
\\
\circ  
& 
\circ 
\end{tikzcd}
\quad
\parrow
\quad
\begin{tikzcd}[font=\small, mymatr,  arrow style=tikz, 
>={Triangle[length=4mm]}, cells={nodes={inner sep=5mm}}, 
row sep=1cm, column sep=1.5cm]
\circ\arrow[dd, "(\mbf{s}\cdot\mbf{z})(1)" ']
&
\circ\arrow[dd, "(\mbf{s}\cdot\mbf{z})(2)" ]
\\
&
\\
\circ 
&
\circ\end{tikzcd}
\label{e:vmtwo}
\end{equation}

The horizontal multiplication 
${\boldsymbol{s}}\hmult{\boldsymbol{z}}$
of 
two twin arrows 
  $\mbf{s}$
and
$\boldsymbol{z}$ in ${\sCat}$, seen as elements of 
$\twinf{\sCatA{\sCat}}$, 
may be 
illustrated by the following diagram, 
which is a special case of~\eqref{e:hm}.
On the right side of the diagram there appears the horizontal multiplication of the twin arrows illustrated on the left side.
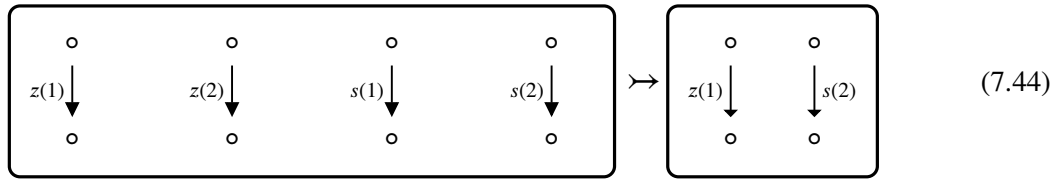
\begin{equation}
\begin{tikzcd}[mymatr, arrow style=tikz, 
>={Triangle[length=2mm]}, 
cells={nodes={inner sep=2mm}}, 
row sep=0.7cm,column sep=1.5cm]
\circ\arrow[d,"z(1)" '] 
&
\circ\arrow[d,"z(2)" '] 	
&
\circ\arrow[d,"s(1)" '] 
&
\circ\arrow[d,"s(2)" '] 	
\\	
\circ	
&
\circ
&
\circ	
&
\circ\end{tikzcd}
\parrow
\begin{tikzcd}[mymatr,  
arrow style=tikz, 
>={Triangle[length=1mm]}, 
cells={nodes={inner sep=2mm}}, 
row sep=0.7cm, 
column sep=0.5cm]
\circ\arrow[d,"z(1)" ']		&\circ\arrow[d,"s(2)" ]\\	
\circ		&\circ	
\end{tikzcd}
\label{e:hmtwo}
\end{equation}

\subsection{Natural Transformations as Homomorphisms}

We are now ready to show in what sense natural transformations 
are indeed homomorphisms 
of functors. It will be useful to keep in mind Remark~\ref{r:tss}.
\begin{definition}
Let $\sCatone$ and $\sCattwo$ be small categories and 
assume that 
$\sFunctor$, $\sFunctortwo$
are functors from $\sCatone$ to $\sCattwo$. 
A \textit{homomorphism} from $\sFunctor$
to
$\sFunctortwo$ is a partial magma homomorphism
\begin{equation}
\mbf{\alpha}:\sCatApm{\sCatone}
\to \twinf{\sCatA{\sCattwo}}
\label{e:urdefinition} 
\end{equation}
(where $\twinf{\sCatA{\sCattwo}}$ is endowed with horizontal multiplication)
such that
\begin{equation}
\mbf{\alpha}(x)\in 
\cathom{\twin{\sCattwo}}{\sFunctor x}{\sFunctortwo x}
\quad
\forall x\in\sCatApm{\sCatone} \, .
\label{e:urpointwise} 
\end{equation}
\label{d:urhof}
\end{definition}
\begin{remark}
A map $\mbf{\alpha}$
as in Definition~\ref{d:urhof} is a homomorphism for two reasons: 
On the one hand, $\mbf{\alpha}$
is a (global) partial magma homomorphism from 
$\sCatApm{\sCatone}$ and 
$\twinf{\sCatA{\sCattwo}}$ (endowed with horizontal multiplication).
On the other hand, \eqref{e:urpointwise}
says that $\mbf{\alpha}$ is a (pointwise) homomorphism, since 
its value 
at each point 
$x\in\sCatApm{\sCatone}$ is a homomorphism between 
$\sFunctor x$
and
$\sFunctortwo x$ in the twin category associated to $\sCattwo$.
\end{remark}
\begin{remark}
The requirement that~\eqref{e:urdefinition} 
is a partial magma homomorphism, where $\twinf{\sCatA{\sCattwo}}$ 
is endowed with horizontal multiplication, means that 
for each $x,y\in\sCatA{\sCatone}$, 
if $x\cdot y$ is defined in $\sCatA{\sCatone}$
then $\hm{\mbf{\alpha}(x)}{\mbf{\alpha}(y)}$
is defined in $\twinf{\sCatA{\sCattwo}}$ and 
\begin{equation}
\mbf{\alpha}(x\cdot y)
=
\hm{\mbf{\alpha}(x)}{\mbf{\alpha}(y)} \, .
\label{e:global} 
\end{equation}
The requirement that~\eqref{e:urpointwise} holds and Lemma~\ref{l:meaning} imply that 
\begin{equation}
\begin{tikzcd}[row sep=large, column sep=large]
\sFunctor \ssdom{x}
\arrow[r,"\sFunctor x"]
\arrow[d, "{{\mbf{\alpha}_1}}({x})" ']
 & 
\arrow[d, "{{\mbf{\alpha}_2}}({x})" ]
\sFunctor \sscodom{x}
\\ 
\sFunctortwo \ssdom{x}
\arrow[r, "\sFunctortwo x" ']
& 
\sFunctortwo \sscodom{x}
\end{tikzcd} 
\label{e:premise}
\end{equation}
\label{r:basicr}
\end{remark}
\begin{lemma}
If $\sCatone$ and $\sCattwo$ are small categories then 
the collection of functors from $\sCat$ to $\sCattwo$, denoted by 
$\ch{\sCatone}{\sCattwo}$, 
is the class of objects of a category, denoted $\functors{\sCatone}{\sCattwo}$, whose arrows are defined as follows: If $\sFunctor,\sFunctortwo\in\ch{\sCatone}{\sCattwo}$
then 
$$
\cathom{\functors{\sCatone}{\sCattwo}}{\sFunctor}{\sFunctortwo}
\eqdef
\{\mbf{\alpha}: \mbf{\alpha}
\text{ is a homomorphism from $\sFunctor$
to
$\sFunctortwo$}
\} \, .
$$
\end{lemma}
\begin{proof}
If $\sFunctor,\sFunctortwo,\sFunctorthree\in\ch{\sCatone}{\sCattwo}$,
$\mbf{\alpha}\in\cathom{\functors{\sCatone}{\sCattwo}}{\sFunctor}{\sFunctortwo}$, 
and 
$\mbf{\beta}\in\cathom{\functors{\sCat}{\sCattwo}}{\sFunctortwo}{\sFunctorthree}$, 
we define 
$$\mbf{\beta}\circ\mbf{\alpha}\in\cathom{\functors{\sCat}{\sCattwo}}{\sFunctor}{\sFunctorthree}$$
by setting $(\mbf{\beta}\circ\mbf{\alpha})(x)\eqdef\vm{\mbf{\beta}(x)}{\mbf{\alpha}(x)}$,
the vertical multiplication of $\mbf{\beta}(x)$ and $\mbf{\alpha}(x)$.
If $x\cdot y$ is defined in $\sCatA{\sCatone}$
then 
\begin{equation}
\begin{split}
(\hm{\mbf{\beta}\circ\mbf{\alpha})(x)}{(\mbf{\beta}\circ\mbf{\alpha})(y)}
&=\hm{
[
\vm
{\mbf{\beta}(x)}
{\mbf{\alpha}(x)}
]
}
{
[
\vm
{\mbf{\beta}(y)}
{\mbf{\alpha}(y)}
]
}\\
&=
\vm
{
[
\hm
{\mbf{\beta}(x)}
{\mbf{\beta}(y)}
]
}
{
[
\hm
{\mbf{\alpha}(x)}
{\mbf{\alpha}(y)}
]
}
\\
&=\vm
{\mbf{\beta}(xy)}
{\mbf{\alpha}(xy)}
=
(\mbf{\beta}\circ\mbf{\alpha})(xy) 
\end{split}
\end{equation}
where in the second equality we applied Lemma~\ref{l:doublepartialmagma}.
Hence $\mbf{\beta}\circ\mbf{\alpha}$ is a partial magma homomorphism from 
$\sCatApm{\sCatone}$
to $\twinf{\sCatA{\sCattwo}}$. Proposition~\ref{p:twin} implies that 
$\mbf{\beta}\circ\mbf{\alpha}(x)\in
\cathom{\twin{\sCattwo}}{\sFunctor x}{\sFunctorthree x}$. The verification 
of the remaining axioms of category is straightforward and is left to the reader.
\end{proof}
\begin{examples}
If $\sCat$ is a category then 
\begin{equation}
{\sCat}\text{ is isomorphic to }\functors{{\scOne}}{\sCat}
\label{e:functorspace1} 
\end{equation}
and 
\begin{equation}
\twinf{\sCat}\text{ is isomorphic to }\functors{{\scTwo}}{\sCat} \, .
\label{e:functorspace2} 
\end{equation}
\end{examples}

\begin{lemma}
If $\mbf{\alpha}$ is a homomorphism from $\sFunctor$ to $\sFunctortwo$, as in Definition~\ref{d:urhof}, then for each object $u$ in $\sCatone$ there exists a unique 
arrow 
\begin{equation}
\tau_{\alpha}(u)\in\cathom{\sCattwo}{\sFunctor u}{\sFunctortwo u}
\label{e:alphatilde} 
\end{equation}
such that
\begin{equation}
\mbf{\alpha}(1_u)=(\tau_{\alpha}(u),\tau_{\alpha}(u)) \, .
\label{e:tte} 
\end{equation}
\label{l:tte}
\end{lemma}
\begin{proof}
If $u$ is an object in 
${\sCatone}$ then 
$\sFunctor u$ and $\sFunctortwo u$ are objects in $\sCattwo$. 
Then~\eqref{e:urpointwise} and 
Lemma~\ref{l:recaptures} imply that there exists a unique
$\tau_{\alpha}(u)\in\cathom{\sCattwo}{\sFunctor u}{\sFunctortwo u}$
such that~\eqref{e:tte} holds. 
\end{proof}
\begin{corollary}
If $\mbf{\alpha}$ is a homomorphism from $\sFunctor$ to $\sFunctortwo$, as in Definition~\ref{d:urhof}, then $\mbf{\alpha}$
is unital as a partial magma homomorphism from $\sCatApm{\sCatone}$ to 
$\twinf{\sCatA{\sCattwo}}$.  
\label{c:tte}
\end{corollary}
\begin{proof}
It suffices to show that if $x$ is a unit in 
$\sCatApm{\sCatone}$
then $\mbf{\alpha}(x)$ is a unit in $\twinf{\sCatA{\sCattwo}}$. 
Recall that 
the units in $\twinf{\sCatA{\sCattwo}}$
(under horizontal multiplication) are the elements $\mbf{w}$ with $w(1)=w(2)$.
Hence Lemma~\ref{l:tte} implies that 
$\mbf{\alpha}(x)$ is a unit in $\twinf{\sCatA{\sCattwo}}$.
\end{proof}

\begin{lemma}
If $\mbf{\alpha}$ is a homomorphism from $\sFunctor$ to $\sFunctortwo$ 
in the sense of Definition~\ref{d:urhof}, and 
\begin{equation}
u\in\sCatO{\sCatone}\mapsto
\tau_{\alpha}(u)\in\cathom{\sCattwo}{\sFunctor u}{\sFunctortwo u}
\label{e:fromhtont} 
\end{equation}
is the map described in Lemma~\ref{l:tte}, then 
for each arrow $x$ in $\sCatone$ the following diagram commutes.
\begin{equation}
\begin{tikzcd}[row sep=large, column sep=large]
\sFunctor \ssdom{x}
\arrow[r,"\sFunctor x"]
\arrow[d, "\tau_{\alpha}(\ssdom{x})" ']
 & 
\arrow[d, "\tau_{\alpha}(\sscodom{x})" ]
\sFunctor \sscodom{x}
\\ 
\sFunctortwo \ssdom{x}
\arrow[r, "\sFunctortwo x" ']
& 
\sFunctortwo \sscodom{x}
\end{tikzcd} 
\label{e:conclusion}
\end{equation}
\label{l:firststeptwo}
\end{lemma}
\begin{proof}
Observe that 
if $x$ is an arrow in $\sCatone$ then 
$x\cdot1_{\ssdom{x}}=x$ and
hence
\begin{equation}
\hm{\mbf{\alpha}(x)}{\mbf{\alpha}(1_{\ssdom{x}})}
\text{ is defined in $\twinf{\sCatA{\sCattwo}}$}
\label{e:definedone} 
\end{equation}
and 
\begin{equation}
\hm{\mbf{\alpha}(x)}{\mbf{\alpha}(1_{\ssdom{x}})}
=\mbf{\alpha}(x\cdot 1_{\ssdom{x}})=\mbf{\alpha}(x) \, .
\label{e:equalitytwo} 
\end{equation}
Lemma~\ref{l:tte} implies that 
$\mbf{\alpha}(1_{\ssdom{x}})
=(\tau_{\alpha}({\ssdom{x}}),\tau_{\alpha}({\ssdom{x}}))$. 
Hence $\tau_{\alpha}({\ssdom{x}})=\mbf{\alpha}_1(x)$. 
Since $1_{\sscodom{x}}\cdot x=x$, a similar argument shows that 
$\tau_{\alpha}({\sscodom{x}})=\mbf{\alpha}_2(x)$.
Hence~\eqref{e:conclusion} follows from~\eqref{e:premise}.
\end{proof}

\begin{theorem}
The notion of homomorphism between two functors $\sFunctor$
and $\sFunctortwo$, 
described in Definition~\ref{d:urhof},
is equivalent to the standard notion of natural transformation from 
 $\sFunctor$
to
$\sFunctortwo$. 
\label{t:proofofequivalence}
\end{theorem}
\begin{proof}
On the one hand, Lemma~\ref{l:firststeptwo} implies that 
if $\mbf{\alpha}$ is a homomorphism from $\sFunctor$ to $\sFunctortwo$ 
in the sense of Definition~\ref{d:urhof}, then the map 
$\tau_{\alpha}$ defined in~\eqref{e:fromhtont}
is a natural transformation from $\sFunctor$ to $\sFunctortwo$.
On the other hand, if 
\begin{equation}
u\in\sCatO{\sCatone}\mapsto
\tau(u)\in\cathom{\sCattwo}{\sFunctor u}{\sFunctortwo u}
\label{e:nt} 
\end{equation}
is a natural transformation from $\sFunctor$ 
to $\sFunctortwo$, and we define $\mbf{\alpha}(x)\eqdef(\tau(\ssdom{x}),\tau(\sscodom{x}))$, then a straightforward verification shows that 
$\mbf{\alpha}$ is 
a homomorphism from $\sFunctor$ to $\sFunctortwo$ 
in the sense of Definition~\ref{d:urhof}.
\end{proof}

\section{Summary of Results}

Let $\smm=\standard$ be a complete measure space of finite measure. 
The Lebesgue transform of an integrable function 
$\classbf$ on $\smm$ 
encodes the collection of all the mean-values of $\classbf$ on all measurable 
subsets of $\mm$ of positive measure. In the problem of the 
differentiation of integrals, one seeks to recapture $\classbf$ from 
its Lebesgue transform. In previous work we showed that, 
in all known results, $\classbf$ may be recaputed from its Lebesgue 
transform by means of a limiting process 
associated to an appropriate family of filters defined on the collection 
of all measurable subsets of $\smm$ of positive measure. The first result 
of the present work is that the existence of such a limiting process is equivalent 
to the existence of a Von Neumann-Maharam lifting of $\smm$.

In the second result of this work we provide an independent argument 
that 
shows that 
the recourse to filters is 
a \textit{necessary consequence} of the requirement that 
the process of recapturing $\classbf$ from its mean-values 
is associated to a \textit{natural transformation}, in the sense of category theory.
This result essentially follows  from the
Yoneda lemma.
As far as we know, this is the first instance of a significant interaction between 
category theory and the problem of the differentiation of integrals. 

In the Appendix 
we have proved, in a precise sense, that 
\textit{natural transformations fall within the general
concept of homomorphism}.
As far as we know, this is a novel conclusion: Although 
it is often said that natural transformations are homomorphisms of functors, 
this statement appears to be presented as a mere  analogy, not in 
a precise technical sense.
In order to achieve this result, 
we had to bring to the foreground a notion 
that is implicit in the subject but has remained hidden in the background, i.e., that of \textit{partial magma}.   

\paragraph{There is no conflict of interest. \\
Useful advice on category theory and logic from Gianluca Amato 
is gratefully acknowledged. \\
Partial support from Indam-Gnampa and Fondi di Ricerca di Ateneo is gratefully acknowledged.}

\end{document}